\documentclass[mathpazo]{cicp}

\usepackage{amsmath}
\usepackage{amsthm}
\usepackage[T1]{fontenc}
\usepackage[utf8]{inputenc}
\usepackage{units}
\usepackage{subfigure,graphicx}
\usepackage{microtype}
\usepackage{zymacros}

\makeatletter

\usepackage{changepage}

\usepackage{url}% simple URL typesetting
\usepackage{booktabs}% professional-quality tables
\usepackage{amsfonts}% blackboard math symbols
\usepackage{nicefrac}% compact symbols for 1/2, etc.
\usepackage{caption}
\usepackage{xcolor}

\@ifundefined{showcaptionsetup}{}{%
 \PassOptionsToPackage{caption=false}{subfig}}
\usepackage{subfig}
\makeatother

\newcommand{\bs}[1]{\boldsymbol{#1}}

\graphicspath{{./figs/}}

\begin{document}

\title{Multi-scale Deep Neural Network (MscaleDNN) Methods for Oscillatory Stokes Flows in Complex Domains}

\author[Wang, B. et.~al.]{Bo Wang\affil{1}, Wenzhong Zhang\affil{2},  Wei Cai\affil{2}\comma\corrauth}
\address{
\affilnum{1}\ LCSM(MOE), School of Mathematics and Statistics, Hunan Normal University, Changsha, Hunan, 410081, P. R. China. \\
\affilnum{2}\ Dept. of Mathematics, Southern Methodist University, Dallas, TX 75275}
\emails{{\tt bowang@hunnu.edu.cn} (Bo Wang), {\tt wenzhongz@smu.edu} (Wenzhong Zhang), {\tt cai@smu.edu} (W.~Cai). Date: September 26, 2020, submitted to CiCP special issue on Machine Learning for Scientific Computing.}

\begin{abstract}
In this paper, we study a multi-scale deep neural network (MscaleDNN) as a meshless numerical method for computing oscillatory Stokes flows in complex domains.
The MscaleDNN employs a multi-scale structure in the design of its DNN using radial scalings to convert the approximation of
high frequency components of the highly oscillatory Stokes solution to one of lower frequencies. The MscaleDNN solution to the Stokes problem
is obtained by minimizing a loss function in terms of $L^2$ norm of the residual of the Stokes equation. Three forms of loss functions
are investigated based on vorticity-velocity-pressure, velocity-stress-pressure, and velocity-gradient of velocity-pressure formulations of the Stokes equation. We first conduct a systematic study of the MscaleDNN methods with various loss functions on the Kovasznay flow in comparison with
normal fully connected DNNs. Then, Stokes flows with highly oscillatory solutions in a 2-D domain with six randomly placed holes are simulated by the MscaleDNN. The results show that MscaleDNN has faster convergence and consistent error decays in the simulation of Kovasznay flow for all four tested loss functions. More importantly, the MscaleDNN is capable of learning highly oscillatory solutions when the normal DNNs fail to converge.
\end{abstract}


\maketitle

\section{Introduction}

Numerical methods for incompressible flow is one of the major topics in computational fluid dynamics, which has been intensively studied over last five decades. Various techniques have been proposed to address the incompressibility condition of the flow, including projection methods \cite{chorin69} \cite{temam69}, Gauge methods \cite{weinan03}, and time splitting methods \cite{gk91}, among others. Finite element and spectral element methods \cite{canuto87} are mostly used to discretize the Navier-stokes equation where special attentions are needed for the approximation spaces of velocity and pressure to satisfy the Babuska and Brezzi inf-sup condition for a saddle point problem \cite{raviart12}. Besides, for large scale engineering applications, body-fitted mesh generations for 3-D objects and efficient linear solvers for the resulting linear systems have been a major issue for computational resources.

The emerging deep neural network (DNN) has found many applications beyond its traditional applications such as image classification and speech recognition. Recent work in extending DNNs to the field of scientific and engineering computing has shown much promise \cite{weinan18}\cite{han18}\cite{gk19}.
DNN based numerical methods are usually formulated as an optimization problem where the loss function could be an energy functional as in a Ritz formulation of a self-adjoint differential equation \cite{weinan18} or simply the least squared mean of the residual of the PDEs \cite{gk20} \cite{cai19}\cite{cai20}.
The DNN technique provides a powerful approximation method to represent solutions of high dimensional variables while the traditional finite element and spectral element methods encounter the well known curse of dimensionality problem. Also, there are several advantages of using DNN to approximate the solution of the incompressible flows. Firstly, the stochastic optimization algorithm employed by DNN based methods relies on loss calculated on randomly sampled points in the computational domain rather than over an unstructured mesh fitting the geometry of the complex objects in the fluid problem. This feature renders the DNN-based methods for solving PDEs a truly mesh-less method. Secondly, due to the capability of the DNN in handling high dimensional functions, the approximation of a time dependent solution can be carried out in the temporal-spatial four dimensional space.
%, thus eliminating the need of time marching schemes with strict stability requirement on the time steps by the traditional CFD algorithms.
Thirdly, boundary conditions for the fluid problems can be simply enforced by introducing penalty terms in the loss function and no need to find and implement appropriate and non-trivial boundary conditions for pressure \cite{orszag86} or vorticity variables in corresponding formulations for the Stokes or Navier-Stokes equations.

Normal fully connected DNNs used for image classification and data science applications have been shown to be ineffective in learning high frequency contents of the solution as illustrated in recent works on DNNs' frequency dependent convergence properties \cite{xu20}. Unfortunately, fluid flow at high Reynolds number will contain many scales, which is the hallmark of the onset of turbulent flow from a laminar one. Therefore, in order to make the DNN based approaches to be  competitive numerical methods, in terms of resolution power, to popular spectral \cite{canuto87} and spectral element methods \cite{gk89}, it is important to develop new classes of DNNs which can represent scales of drastic disparities arising from the study of turbulent flows.  For this purpose, we have recently developed strategies to speed up the convergence of DNNs in learning high frequency content of the solutions of PDEs. Two new DNNs have been proposed: a PhaseDNN \cite{cai19} and a MscaleDNN \cite{cai20}. The PhaseDNN uses a series of phase shifts to convert high frequency contents to a low frequency range before the learning is carried out. This method has been shown to be very effective in simulating high frequency Helmholtz equations in acoustic wave scattering. On the other hand, the MscaleDNN uses a radial scaling technique in the frequency domain (or a corresponding scaling in the physical domain) to convert solution content of a range of higher frequency to a lower frequency one, which will be learned quickly with a small size DNN, and the latter is then scaled back in the physical space to approximate the original solution content. MscaleDNN is more effective to handle higher dimensional PDEs and has already been shown to be superior over traditional fully connected DNNs for solving Poisson-Boltzmann equation in complex and singular domains \cite{cai20}. In this paper, we will extend the MscaleDNN approach to find the solution of the Stokes problem as a first step to develop DNN based numerical methods for time-dependent incompressible Navier-Stokes equations.

The rest of the paper is organized as follows. In section 2, we will present the structure of the MscaleDNN to be used for solving the Stokes problems. Section 3 will propose several loss functions for training, based on three different first order system reformulations of the Stokes equation. A benchmark test on a low frequency Kovasznay flow will be conducted  in section 4 to evaluate the performance of normal fully connected DNN and MscaleDNNs as well as different loss functions. Section 5 will present the numerical tests of highly oscillatory Stokes flows with multiple frequencies in a complex domain. Finally, a conclusion and discussion of future work are given in Section 6.

\section{Multi-scale DNN (MscaleDNN)}

\begin{figure}[htbp]
\centering
\includegraphics[width=0.5\linewidth]{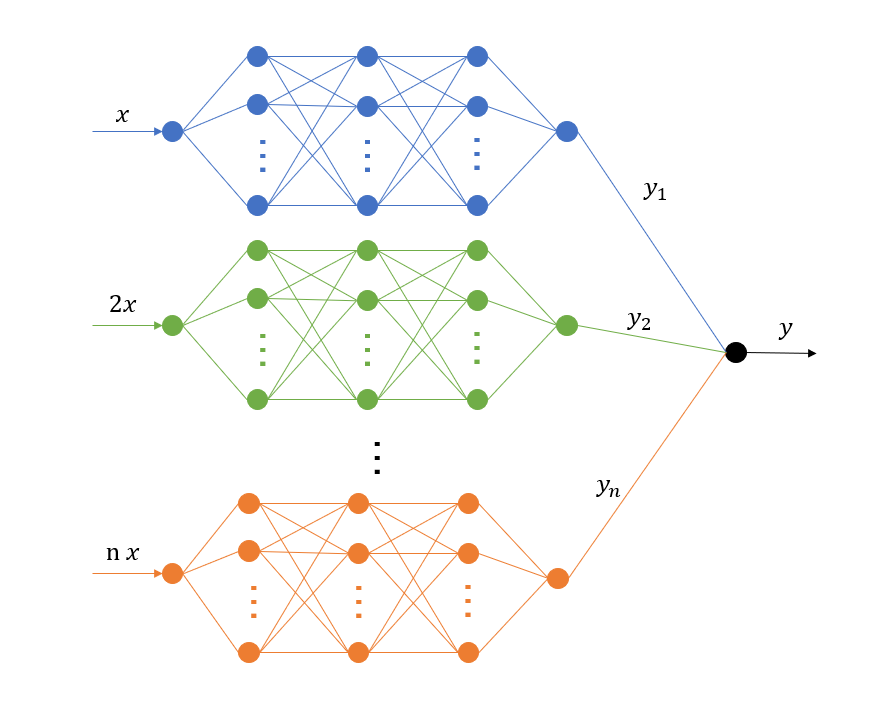}
\caption{Illustration of a MscaleDNN.}
\label{net}
\end{figure}

In a recent work \cite{cai20}, a multi-scale DNN was proposed, which consists of a series of parallel normal sub-neural networks.  Each of the sub-networks will receive a scaled version of the input and their outputs will then be combined to make the final out-put of the MscaleDNN (refer to Fig. \ref{net}). The individual sub-network in the MscaleDNN with a scaled input is designed to approximate a segment of frequency content of the targeted function and the effect of the scaling is to convert a specific high frequency content to a lower frequency range so the learning can
be accomplished much quickly. Recent work \cite{xu20}  on the frequency dependence of the DNN convergence shows that much faster convergence occurs in approximating low frequency function compared with approximating high frequency ones, the MscaleDNN takes advantage of this property.
 In addition, in order to produce scale separation and identification capability for a MscaleDNN, we borrowed the idea of compact mother scaling and wavelet functions from the wavelet theory \cite{daub92}, and found that the activation functions with a localized frequency profile works better than normal activation functions, e.g., ReLU, tanh, etc.

Fig. \ref{net} shows the schematics of a MscaleDNN consisting of $n$ networks. Each scaled input passing through a sub-network can be expressed in the following formula
\begin{equation}
    f_{\vtheta}(\vx) = \vW^{[L-1]} \sigma\circ(\cdots (\mW^{[1]} \sigma\circ(\mW^{[0]} (\vx) + \vb^{[0]} ) + \vb^{[1]} )\cdots)+\vb^{[L-1]}, \label{mscalednn}
\end{equation}
where  $W^{[1]}$ to $W^{[L-1]}$ and $b^{[1]}$ to $b^{[L-1]}$ are the weight matrices and bias unknowns, respectively, to be optimized via the training, $\sigma(x)$ is the activation function. In this work, the following plane wave activation function will be used for its localized frequency property \cite{cai20},
\begin{equation}
\sigma(x)=\sin(x).
\end{equation}

For the input scales, we could select the scale for the $i$-th sub-network to be $i$ (as shown in Fig. \ref{net}) or $2^{i-1}$.
Mathematically, a MscaleDNN solution $f(\vx)$ is represented by the following sum of sub-networks $f_{\theta^{n_{i}}}$ with network parameters denoted by $\theta^{n_{i}}$ (i.e., weight matrices and bias)
\begin{equation}
f(\vx)\sim
{\displaystyle\sum\limits_{i=1}^{M}}
f_{\theta^{n_{i}}}(\alpha_{i}\vx),\label{f_app}%
\end{equation}
where $\alpha_i$ is the chosen scale for the $i$-th sub-network in Fig. \ref{net}. For more details on the design and discussion of the MscaleDNN, please refer to \cite{cai20}.

For comparison studies in this paper, we will refer to a ``{\bf normal}'' network as an one fully connected DNN with the same total number of neurons as the MscaleDNN, but without multi-scale features. We would perform extensive numerical experiments to examine the effectiveness of different settings and select efficient ones to solve complex problems. All DNN models are trained by Adam \cite{adam}.

\section{Loss functions and the MscaleDNN for Stokes problem}
The following two dimensional (2-D) Stokes problem
\begin{eqnarray}
-\nu\triangle{\bs u}+\nabla p={\bs f}, &{\rm in} \quad\Omega,\label{stokeseq}\\
\nabla\cdot{\bs u}=0, & {\rm in} \quad\Omega,\label{incompressible}\\
{\bs u}={\bs g}, & {\rm on}\quad \partial\Omega,\label{dirichletbc}
\end{eqnarray}
will be solved by the MscaleDNN, here $\Omega$ is an open bounded domain in $\mathbb R^2$, and the boundary condition ${\bs g}$ satisfies a compatibility condition
\begin{equation}
\int_{\partial\Omega}{\bs g}\cdot{\bs n}ds=0.
\end{equation}

The MscaleDNN solution will be found as in the traditional least square finite element method \cite{bochev98} where the solution is obtained by minimizing a loss function in terms of the residual of the Stokes problem (\ref{stokeseq}).
To introduce loss functions for the DNN algorithms, we first reformulate \eqref{stokeseq}-\eqref{dirichletbc} into a first order system as in least square finite element methods for solving Stokes problem. There are various possible ways of recasting (\ref{stokeseq}) into a first order system, and we will focus on the following three popular approaches used in the development of least square finite element methods \cite{bochev98}.

\medskip

\noindent{\bf $\bullet$ Vorticity-velocity-pressure ($\bs\omega$VP) formulation:} The first approach introduces a vorticity variable, a scalar quantity for 2-D flows,
\begin{equation}
\bs\omega=\nabla\times\bs u=\partial_x u_y-\partial_y u_x,
\end{equation}
arriving at a vorticity-velocity-pressure ($\omega$VP) system:
\begin{subequations}\label{vvpsystem}	
\begin{align}
\nu\nabla\times\bs \omega+\nabla p=\bs f,  &\quad{\rm in} \quad\Omega,\label{vvp1}\\
\bs \omega=\nabla\times \bs u, & \quad{\rm in} \quad\Omega,\label{vvp2}\\
\nabla\cdot{\bs u}=\bs 0,&\quad {\rm in} \quad\Omega.\label{vvp3}		
\end{align}
\end{subequations}

%For this formulation, a total of three MscaleDNNs will be used: one for the scalar vorticity $\bs\omega$, one for the velocity vector $u$ where the output $y= \bs u$ in Fig. \ref{net}, and one for the scalar pressure $p$.

\medskip
\noindent{\bf $\bullet$ Velocity-stress-pressure (VSP) formulation:} The second approach introduces a stress tensor
\begin{equation}
\bs T=\sqrt{2}\nu(\nabla\bs u+\nabla\bs u^{\top})/2,
\end{equation}
 while a velocity-stress-pressure (VSP) system
%\begin{eqnarray}
%-\sqrt{2\nu}\nabla\cdot\bs T+\nabla p=\bs f,  &{\rm in} \quad\Omega,\label{vsp1}\\
%\bs T=\frac{\sqrt{2\nu}}{2}(\nabla\bs u+\nabla\bs u^{\top}), & {\rm in} \quad\Omega,\label{vsp2}\\
%\nabla\cdot{\bs u}=\bs 0,& {\rm in} \quad\Omega,\label{vsp3}
%\end{eqnarray}
\begin{subequations}\label{vspsystem}
	\begin{align}
	-\sqrt{2\nu}\nabla\cdot\bs T+\nabla p=\bs f,  &\quad{\rm in} \quad\Omega,\label{vsp1}\\
	\bs T=\frac{\sqrt{2\nu}}{2}(\nabla\bs u+\nabla\bs u^{\top}), &\quad {\rm in} \quad\Omega,\label{vsp2}\\
	\nabla\cdot{\bs u}=\bs 0,&\quad {\rm in} \quad\Omega,\label{vsp3}	
	\end{align}
\end{subequations}
is obtained.

%For this formulation, a total of three MscaleDNNs will be used:  one for the velocity vector $u$ where the output $y= \bs u$ in Fig. \ref{net}, one for the stress tensor $\bs T$ where the output $y= \bs T$ in Fig. \ref{net}, and one for the scalar pressure $p$.

\medskip

\noindent{\bf $\bullet$ Velocity-gradient of velocity-pressure (VgVP) formulation:} The third approach simply introduces a variable $\bs U=\nabla\bs u$ (by taking gradient on each component of the velocity field), which leads to a velocity-gradient of velocity-pressure (VgVP) system
%\begin{eqnarray}
%-\nu\nabla\cdot\bs U+\nabla p=\bs f,  &{\rm in} \quad\Omega,\label{vgvp1}\\
%\bs U=\nabla \bs u, \quad \nabla\cdot{\bs u}=\bs 0& {\rm in} \quad\Omega.\label{vgvp2}
%\end{eqnarray}
\begin{subequations}\label{vgvpsystem}
	\begin{align}
		-\nu\nabla\cdot\bs U+\nabla p=\bs f,  &\quad {\rm in} \quad\Omega,\label{vgvp1}\\
		\bs U=\nabla \bs u, &\quad {\rm in} \quad\Omega.\label{vgvp2}\\
		\nabla\cdot{\bs u}=\bs 0&\quad {\rm in} \quad\Omega.\label{vgvp3}	
	\end{align}
\end{subequations}

%For this formulation, a total of three MscaleDNNs will be used:  one for the velocity vector $u$ where the output $y= \bs u$ in Fig. \ref{net}, one for the tensor of the gradient of velocity  $\bs U$ where the output $y= \bs U$ in Fig. \ref{net}, and one for the scalar pressure $p$.

\medskip
It is well known that it is more difficult to compute the pressure than the velocity in computational fluid dynamics. We find that the velocity also converges faster than pressure in the DNN-based methods. In order to take care of the pressure, we take divergence on both sides of the Stokes equation \eqref{stokeseq} to obtain a Poisson equation
\begin{equation}\label{poissonp}
\Delta p=\nabla\cdot \bs f,\quad{\rm in}\quad\Omega.
\end{equation}
The residual of this equation will be an extra term in the loss function, and a tunable weight on the loss due to pressure is introduced. To be consistent with the first order systems above, we also reformulate the Poisson equation \eqref{poissonp} into
%\begin{eqnarray}
%\bs q=\nabla p,\quad{\rm in}\quad\Omega,\\
%\nabla\cdot\bs q=\nabla\cdot\bs f,\quad{\rm in}\quad\Omega.
%\end{eqnarray}
\begin{subequations}\label{poissonsystem}
	\begin{align}
	\bs q=\nabla p,\quad{\rm in}\quad\Omega,\\
	\nabla\cdot\bs q=\nabla\cdot\bs f,\quad{\rm in}\quad\Omega.
	\end{align}
\end{subequations}

Together with the first order systems \eqref{vvpsystem},\eqref{vspsystem} or \eqref{vgvpsystem}, respectively, we can design the MscaleDNN algorithms. In each algorithm, a total of four MscaleDNNs will be used: one for the velocity vector $\bs u$, one for the pressure $p$, one for the gradient of pressure $\bs q$ and one for the vorticity $\bs\omega$, stress $\bs T$ or the gradient of velocity $\bs U$, respectively. The DNN solutions are denoted by $\bs u(\bs x, \theta_{\bs u}), p(\bs x, \theta_p), \bs{\omega}(\bs x,\theta_{\bs{\omega}}), \bs T(\bs x, \theta_{\bs T}), \bs U(\bs x, \theta_{\bs U}), \bs q(\bs x, \theta_{\bs q})$ accordingly. Based on the first order systems, we define loss functions as follows
\begin{equation}\label{firstorderloss}
\begin{split}
\bs L_{\mbox{$\omega$}VP}(\theta_{\bs u}, \theta_p, \theta_{\bs\omega},\theta_{\bs q}):=&\|\nu\nabla\times\bs \omega+\bs q-\bs f\|^2_{\Omega}+\alpha\|\nabla\cdot\bs q-\nabla\cdot\bs f\|_{\Omega}+\|\nabla\times \bs u-\bs \omega\|^2_{\Omega}\\
&+\|\nabla\cdot\bs u\|^2_{\Omega}+\|\nabla p-\bs q\|^2_{\Omega}+\beta\|\bs u-\bs g\|_{\partial\Omega}^2,\\
\bs L_{VSP}(\theta_{\bs u}, \theta_p, \theta_{\bs T},\theta_{\bs q}):=&\|\sqrt{2\nu}\nabla\cdot\bs T-\bs q+\bs f\|^2_{\Omega}+\alpha\|\nabla\cdot\bs q-\nabla\cdot\bs f\|_{\Omega}+\Big\|\frac{\sqrt{2\nu}}{2}(\nabla\bs u+\nabla\bs u^{\top})-\bs T\Big\|^2_{\Omega}\\
&+\|\nabla\cdot\bs u\|^2_{\Omega}+\|\nabla p-\bs q\|^2_{\Omega}+\beta\|\bs u-\bs g\|_{\partial\Omega}^2,\\
\bs L_{VgVP}(\theta_{\bs u}, \theta_p, \theta_{\bs U},\theta_{\bs q}):=&\|\nu\nabla\cdot\bs U-\bs q+\bs f\|^2_{\Omega}+\alpha\|\nabla\cdot\bs q-\nabla\cdot\bs f\|_{\Omega}+\|\nabla \bs u-\bs U\|^2_{\Omega}\\
&+\|\nabla\cdot\bs u\|^2_{\Omega}+\|\nabla p-\bs q\|^2_{\Omega}+\beta\|\bs u-\bs g\|_{\partial\Omega}^2,
\end{split}
\end{equation}
where $\alpha,\beta$ are penalty constants. We emphasize that the Poisson residual
$$\alpha\|\nabla\cdot\bs q-\nabla\cdot\bs f\|_{\Omega}+\|\nabla p-\bs q\|^2_{\Omega},$$
in the loss function is important for the convergence of the pressure as to be shown via numerical results in Section \ref{section_P}.

For the brevity of notations, the loss functions in \eqref{firstorderloss} are named as $\omega$VP-loss, VSP-loss and VgVP-loss, accordingly. In the rest of this paper, these loss functions will be compared with the simple loss function directly obtained from the original Stokes equation:
\begin{equation}\label{secondorderloss}
\bs L_{VP}(\theta_{\bs u}, \theta_p)=\|\nu\Delta \bs u-\nabla p+\bs f\|^2_{\Omega}+\|\nabla\cdot\bs u\|^2_{\Omega}+\beta\|\bs u-\bs g\|_{\partial\Omega}^2,
\end{equation}
which is named as VP-loss. In the DNN algorithms using this loss function, a total of two MscaleDNNs will be used: one for the velocity vector $\bs u$ where the output $y= \bs u$ in Fig. \ref{net}, and one for the scalar pressure $p$.

\section{Kovasznay flow in a square domain}
As a benchmark test, we first consider the Stokes problem in a square domain $\Omega=[0, 2]\times[-0.5, 1.5]$ with an exact solution coinciding with the analytical solution of the incompressible Navier-Stokes equations obtained by Kovasznay \cite{Kovasznay48}, i.e.,
\begin{equation}\label{kovasznayexactsolution}
\begin{split}
u_1=&1-e^{\lambda x_1}\cos(2\pi x_2),\qquad
u_2=\frac{\lambda}{2\pi}e^{\lambda x_1}\sin(2\pi x_2),\qquad
p=\frac{1}{2}e^{2\lambda x_1},
\end{split}
\end{equation}
where
$$\lambda=\frac{Re}{2}-\sqrt{\frac{Re^2}{4}+4\pi^2},\quad Re=\frac{1}{\nu}.$$
The source term $\bs f$ is obtained by substituting the exact solution into the Stokes equation (\ref{stokeseq}).
We set the viscosity $\nu=0.1$ and investigate the performance of algorithms using fully connected and Multi-scale DNNs. In the simulations of this benchmark problem, all MscaleDNNs are set to have six scales: {\bf $\{x, 2x, 4x, 8x, 16x, 32x\}$} and their fully connected sub-networks all have 4 hidden layers and 50 neurons in each hidden layer. On the other hand, a fully connected DNN with 4 hidden layers and 300 neurons in each hidden layer was tested for comparison. Therefore, the total number of neurons in the fully connected DNN and MscaleDNNs are the same. Nevertheless, the fully connected DNN does have more connectivity with more parameters. In the loss functions, we fix $\alpha=1$ and $\beta=100$. We randomly sample 50000 points inside $\Omega$ and 10000 points on the boundary for learning. In the learning process, we set batch size equal to 1000 points inside the domain and randomly pick 400 points on the boundary for each step.
\medskip

\noindent{\bf $\bullet$ Adaptive learning rates.} We have found that reducing learning rate as the training progresses can have a noticeable improvement in the reduction of loss. In our numerical tests, the learning rate of the first $100$ epochs is set to be $0.001$. Then, the learning rate will be reduced by a factor of 10 after each $100$ epochs. The change of learning rate can be seen clearly in the history of losses.

In order to check the accuracy of the algorithms, we define $\ell^2$-errors
\begin{equation}\label{ell2err}
Err({\bs u})=\Big(\frac{1}{N}\sum\limits_{j=1}^N|\bs u^{\rm DNN}(\bs x_{j})-\bs u(\bs x_{j})|^2\Big)^{\frac{1}{2}},\quad Err(p)=\Big(\frac{1}{N}\sum\limits_{j=1}^N|p^{\rm DNN}(\bs x_{j})-p(\bs x_{j})|^2\Big)^{\frac{1}{2}},
\end{equation}
between the DNN solution $\{\bs u^{\rm DNN}(\bs x), p^{\rm DNN}(\bs x)\}$ and the exact solution $\{\bs u(\bs x),\bs p(\bs x)\}$ given in \eqref{kovasznayexactsolution}. Here, $\{\bs x_{j}=(x_1^j, x_2^j)\}_{j=1}^N$ are locations of a uniform $200\times 200$ mesh of the domain $\Omega$.
\begin{figure}[ht!]
	\center
	\subfigure[loss]{\includegraphics[scale=0.3]{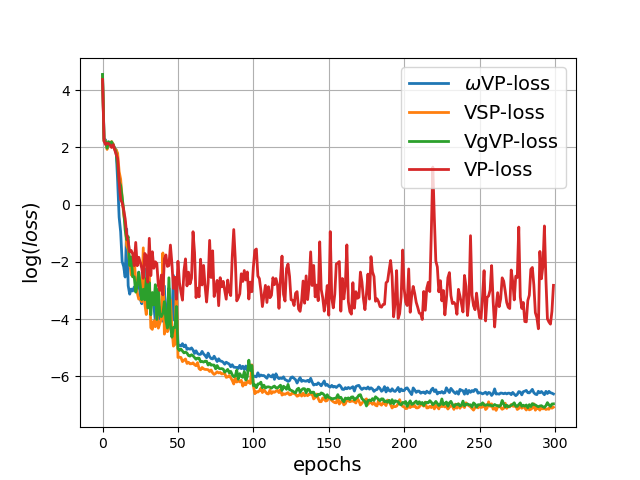}}
	\subfigure[$Err({\bs u})$]{\includegraphics[scale=0.3]{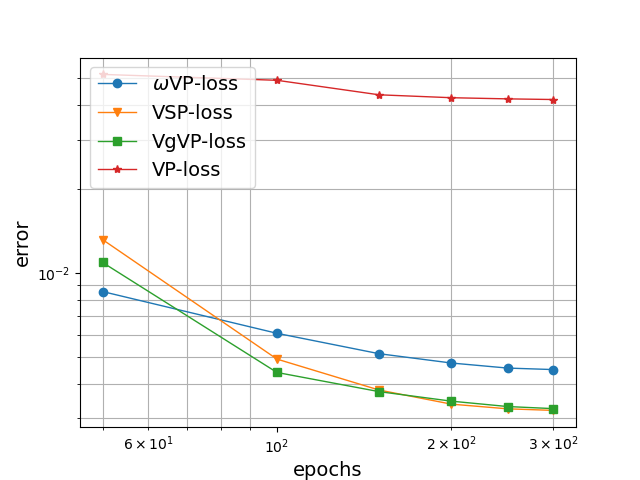}}
	\subfigure[$Err({ p})$]{\includegraphics[scale=0.3]{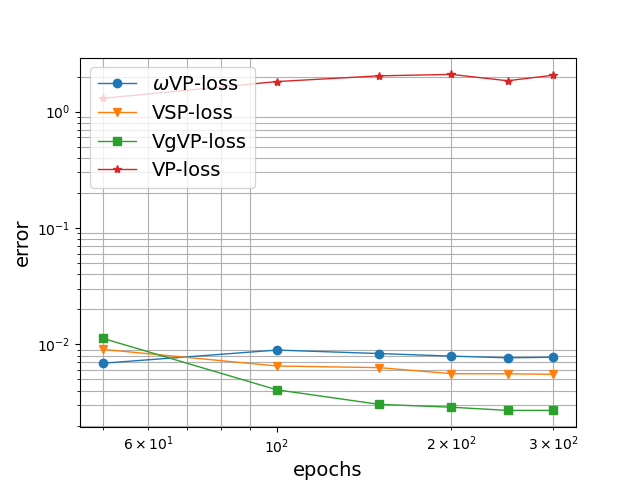}}
	\vspace{-10pt}
	\caption{Normal DNN with different loss functions.}%
	\label{example1-1}%
\end{figure}

\begin{figure}[ht!]
	\center
	\subfigure[loss]{\includegraphics[scale=0.3]{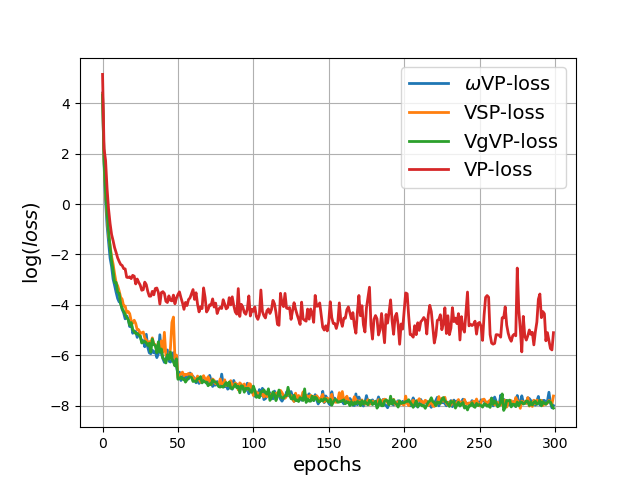}}
	\subfigure[$Err({\bs u})$]{\includegraphics[scale=0.3]{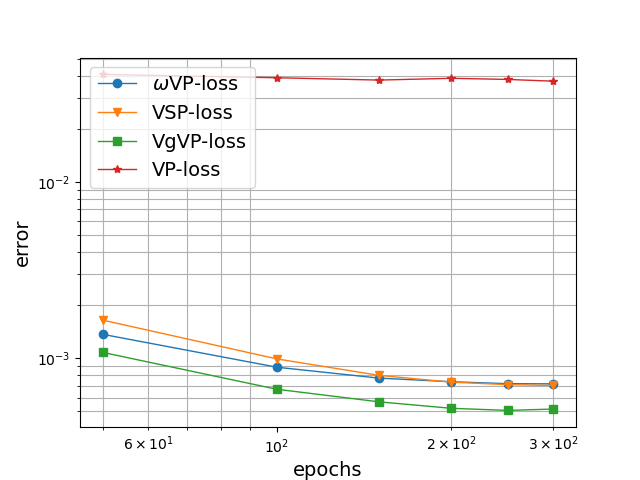}}
	\subfigure[$Err({ p})$]{\includegraphics[scale=0.3]{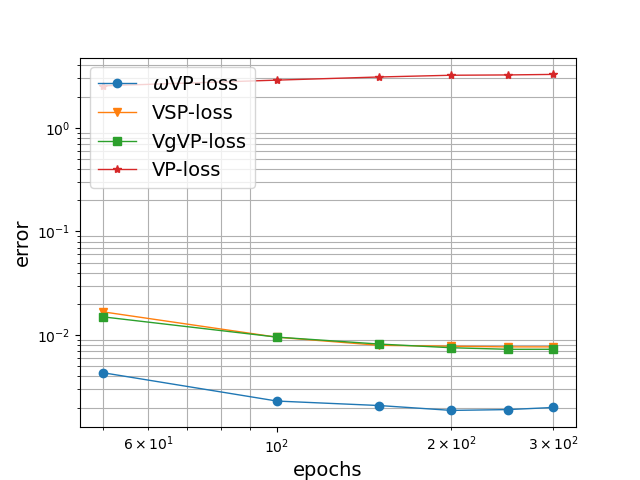}}
	\vspace{-10pt}
	\caption{MscaleDNN with different loss functions.}%
	\label{example1-2}%
\end{figure}

The DNN solutions obtained by minimizing different loss functions in \eqref{firstorderloss}-\eqref{secondorderloss} are compared in Fig.  \ref{example1-1}-\ref{example1-2}.
The results show that both fully connected DNN and MscaleDNNs converge in 300 epochs with any one of the loss functions in \eqref{firstorderloss}. However, the simple VP-loss in \eqref{secondorderloss} has a very poor performance no matter if the fully connected DNN or the MscaleDNNs is used. In particular, both fully connected DNN and MscaleDNNs can not produce reasonable results within 300 epochs if the VP-loss function is used.
\begin{figure}[ht!]
	\center
	\subfigure[loss]{\includegraphics[scale=0.3]{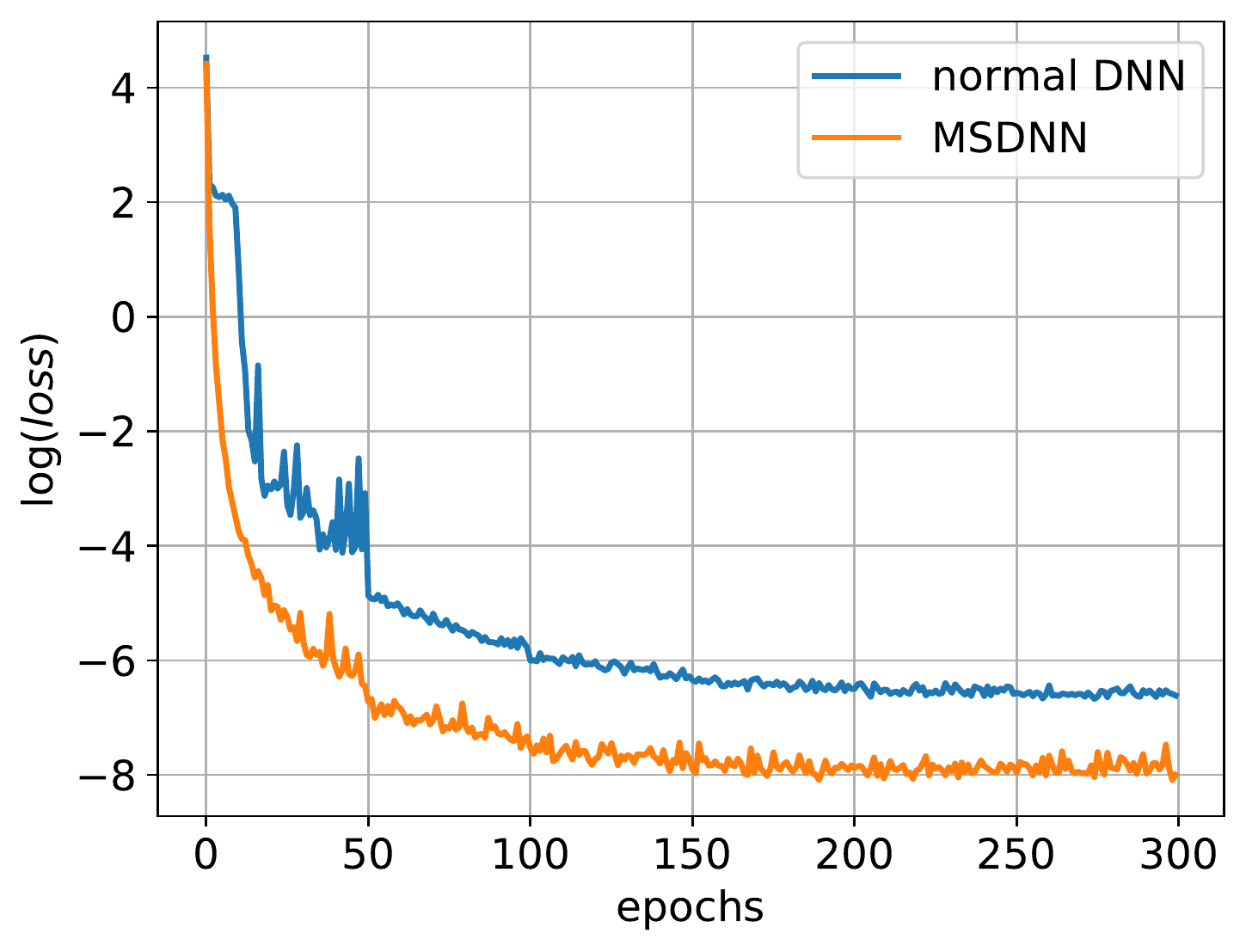}}
	\subfigure[$Err({\bs u})$]{\includegraphics[scale=0.3]{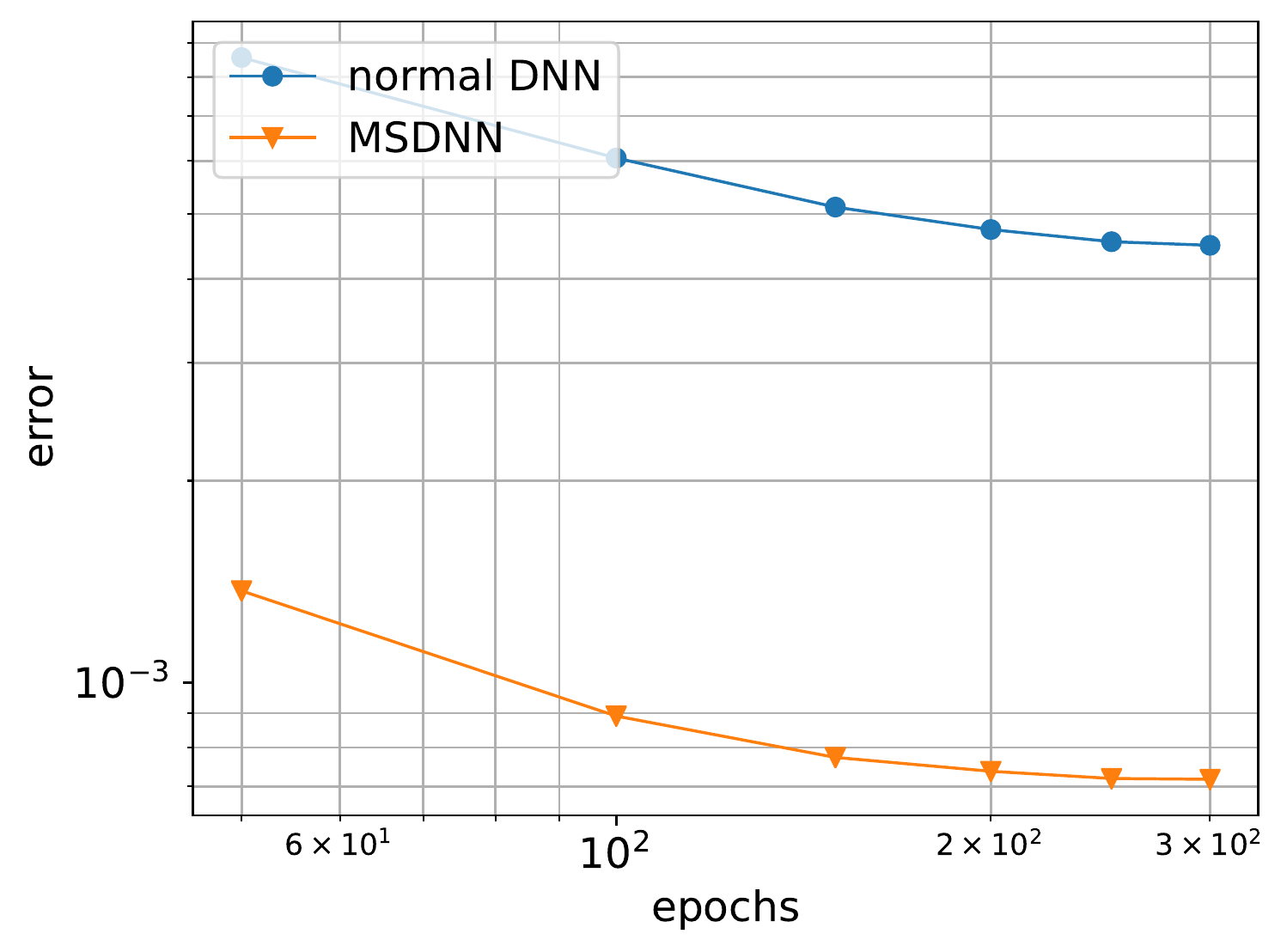}}
	\subfigure[$Err({p})$]{\includegraphics[scale=0.3]{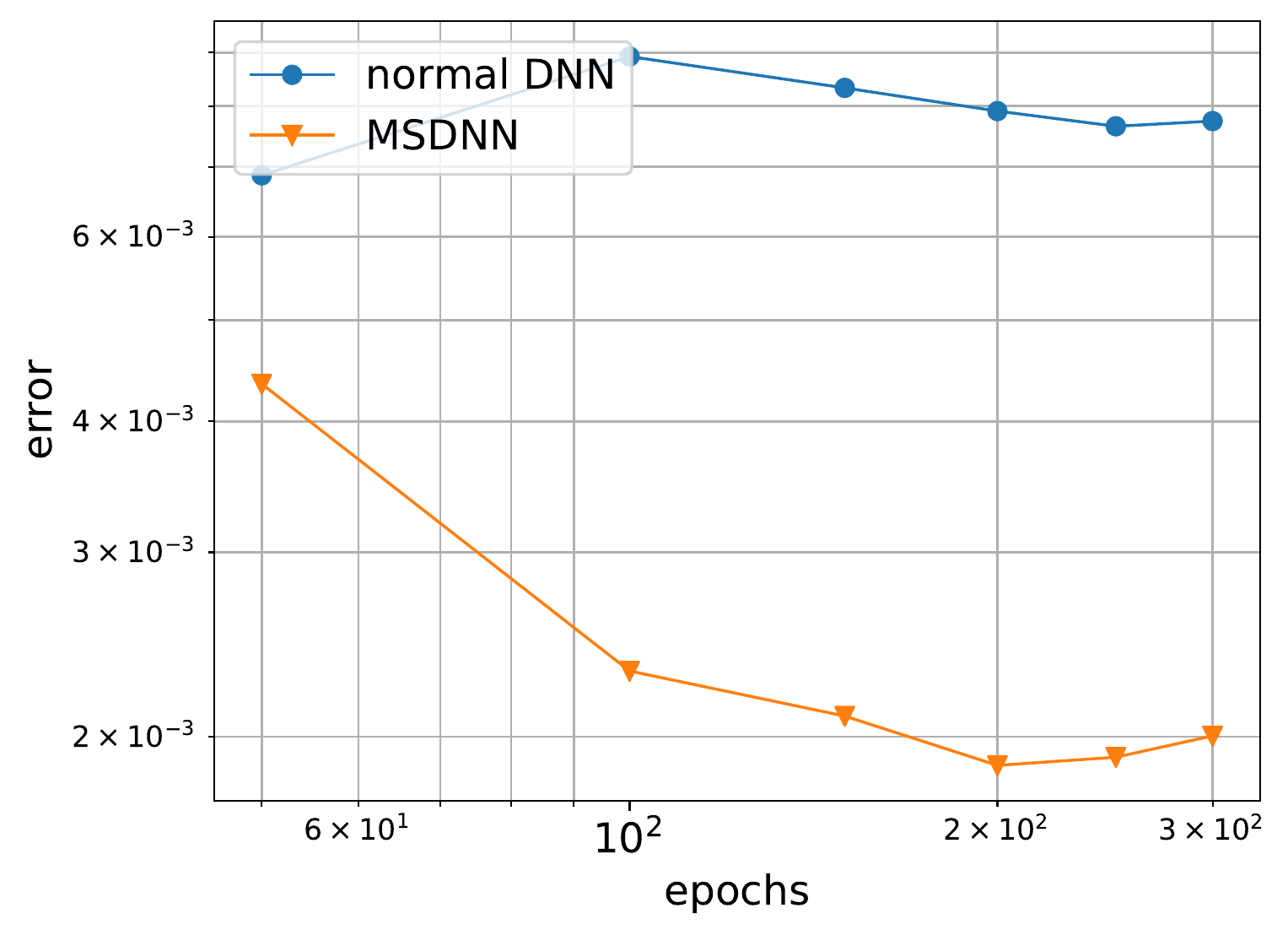}}
	\vspace{-10pt}
	\caption{Normal DNN and MscaleDNN with loss function $\bs L_{\mbox{$\omega$}VP}(\theta_{\bs u},\theta_p, \theta_{\bs{\omega}},\theta_{\bs q})$.}%
	\label{example1-3}%
\end{figure}
\begin{figure}[ht!]
	\center
	\subfigure[loss]{\includegraphics[scale=0.3]{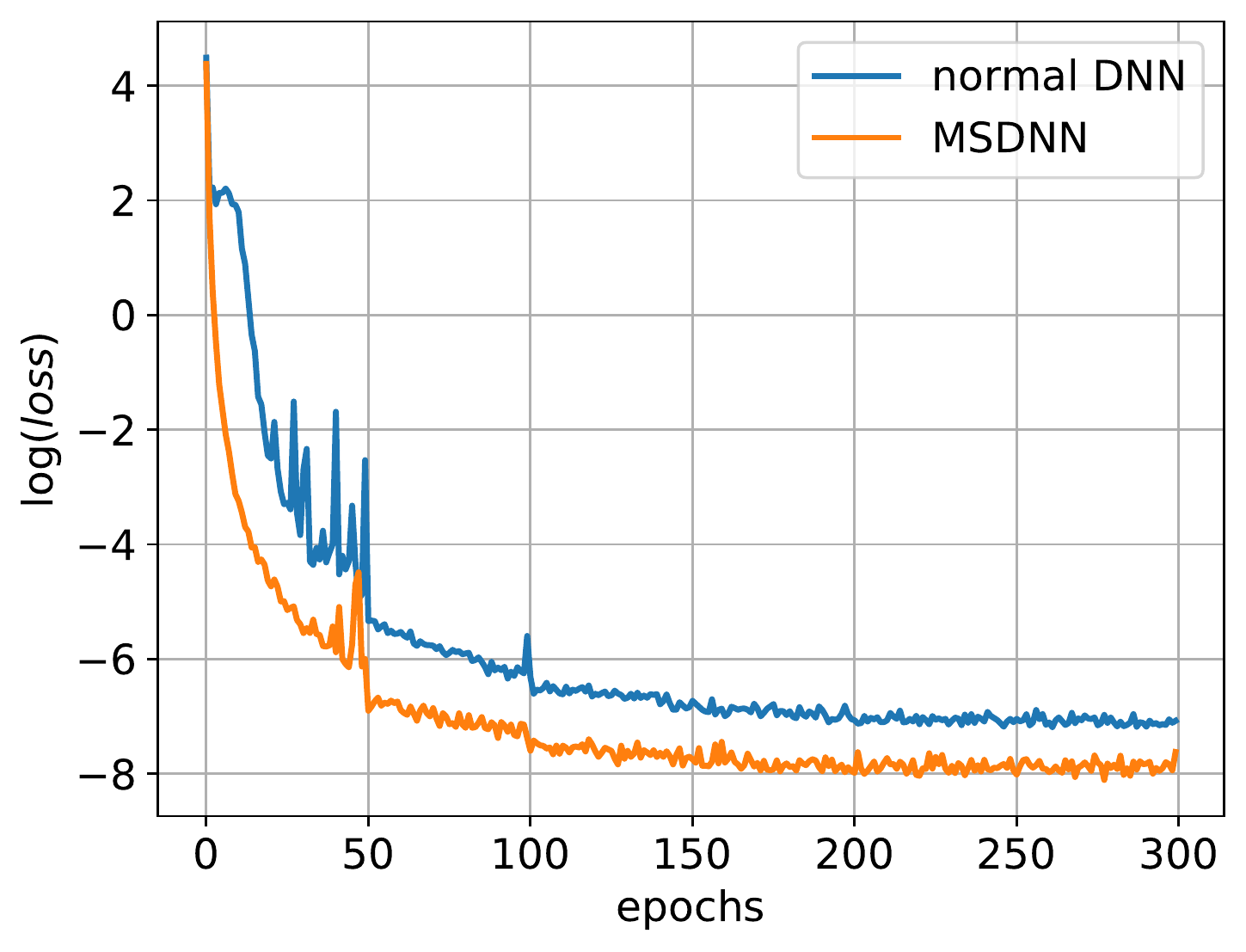}}
	\subfigure[$Err({\bs u})$]{\includegraphics[scale=0.3]{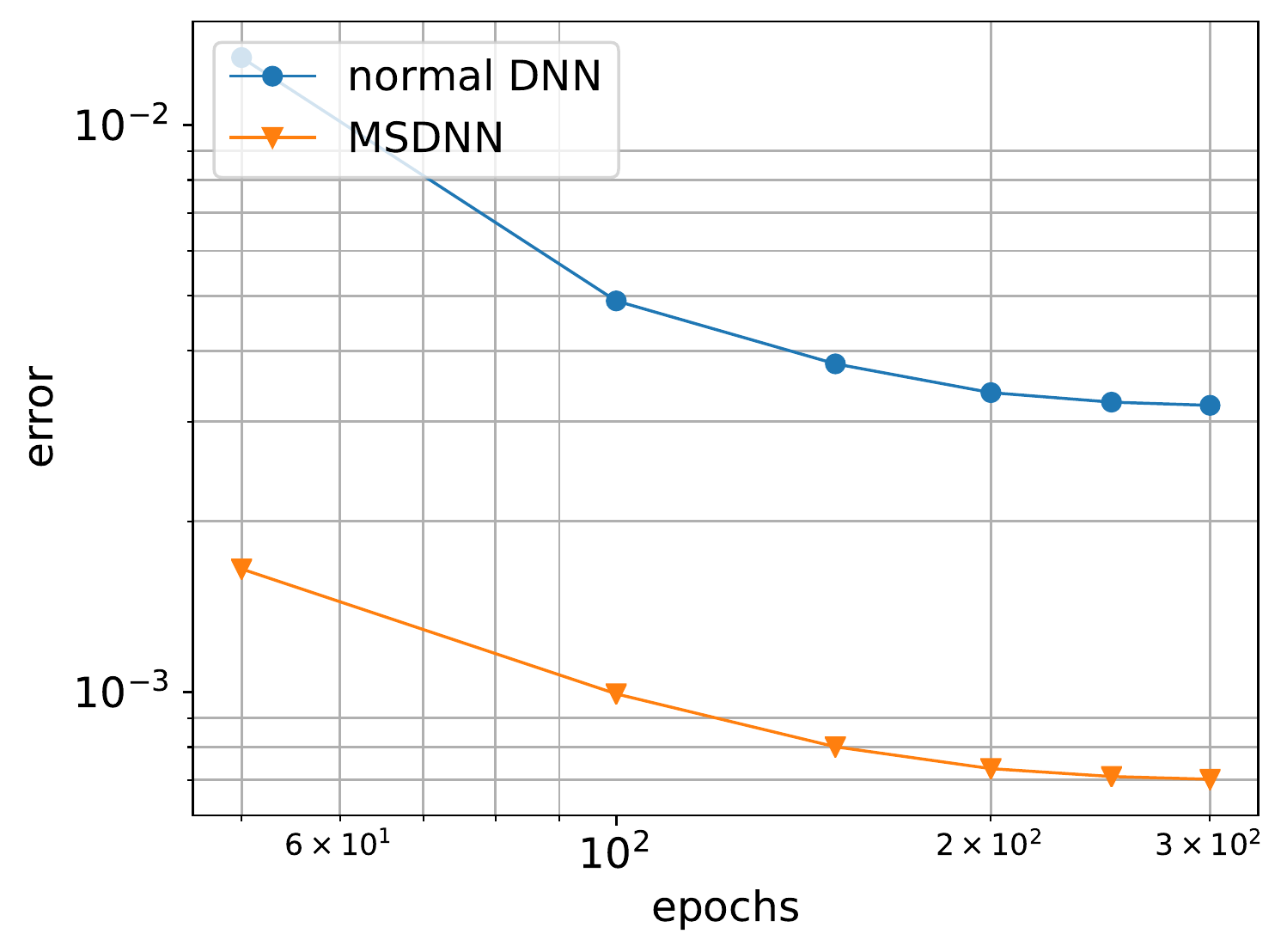}}
	\subfigure[$Err(p)$]{\includegraphics[scale=0.3]{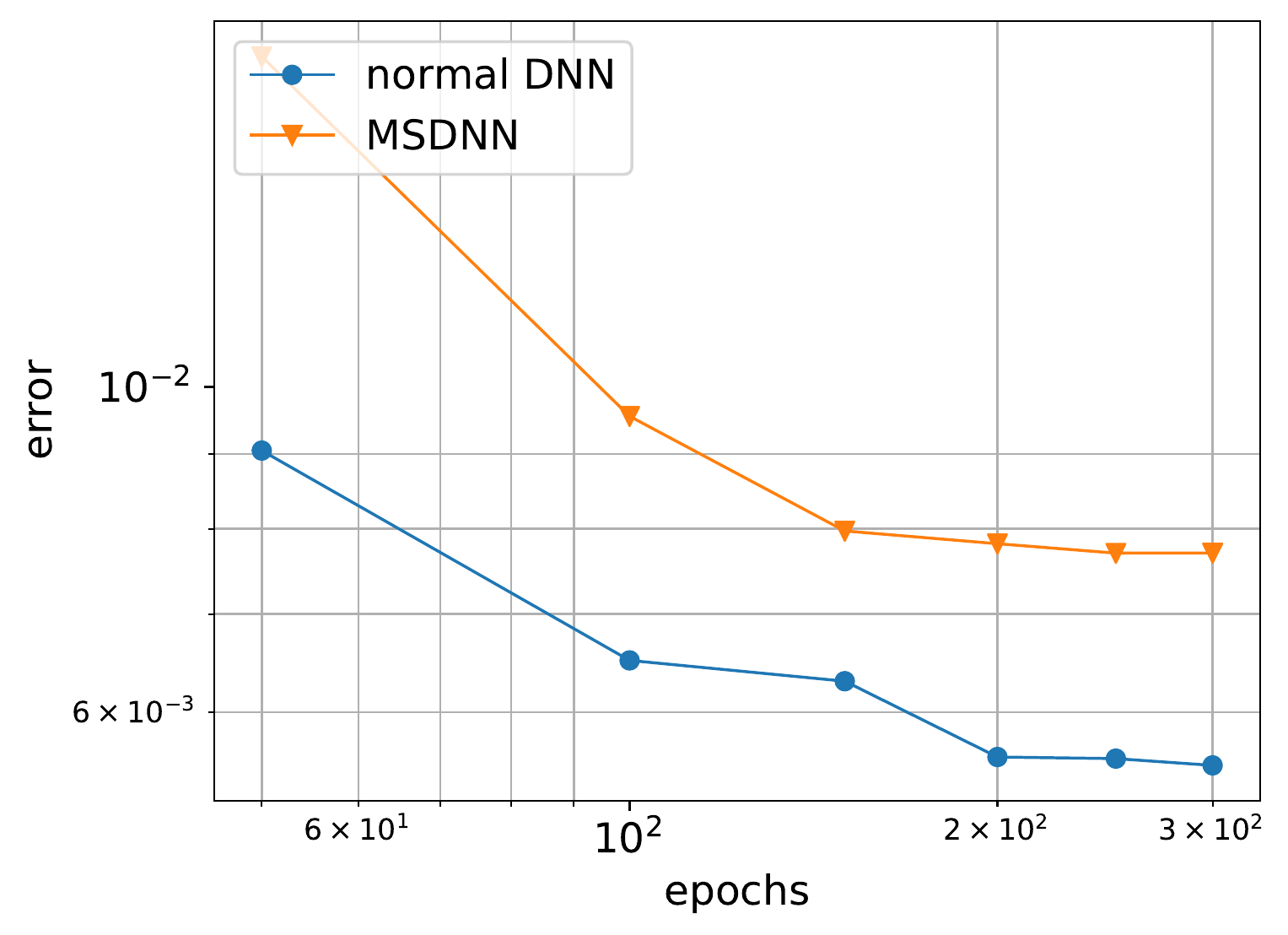}}
	\vspace{-10pt}
	\caption{Normal DNN and MscaleDNN with loss function $\bs L_{VSP}(\theta_{\bs u},\theta_p, \theta_{\bs T},\theta_{\bs q})$.}%
	\label{example1-4}%
\end{figure}
\begin{figure}[ht!]
	\center
	\subfigure[loss]{\includegraphics[scale=0.3]{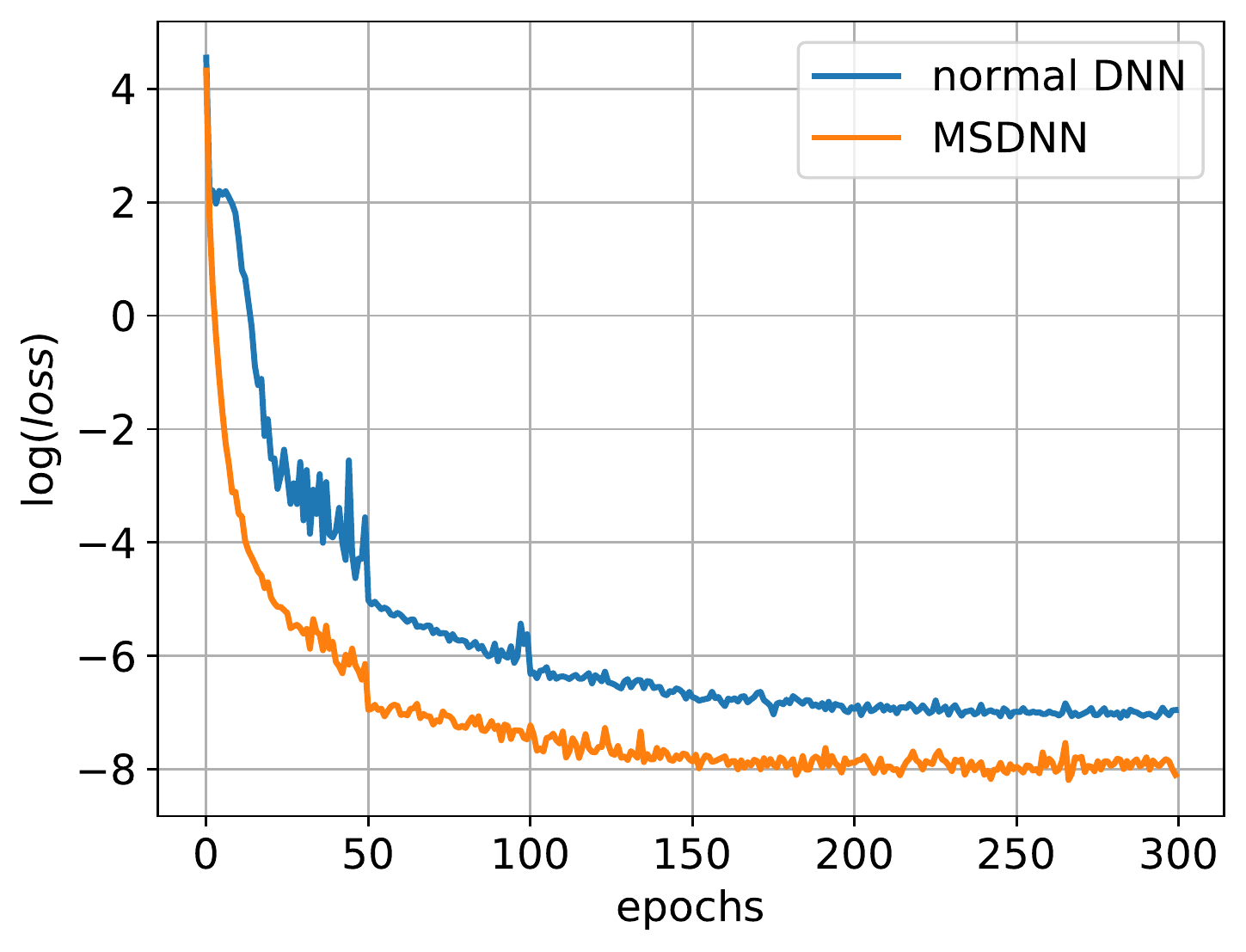}}
	\subfigure[$Err({\bs u})$]{\includegraphics[scale=0.3]{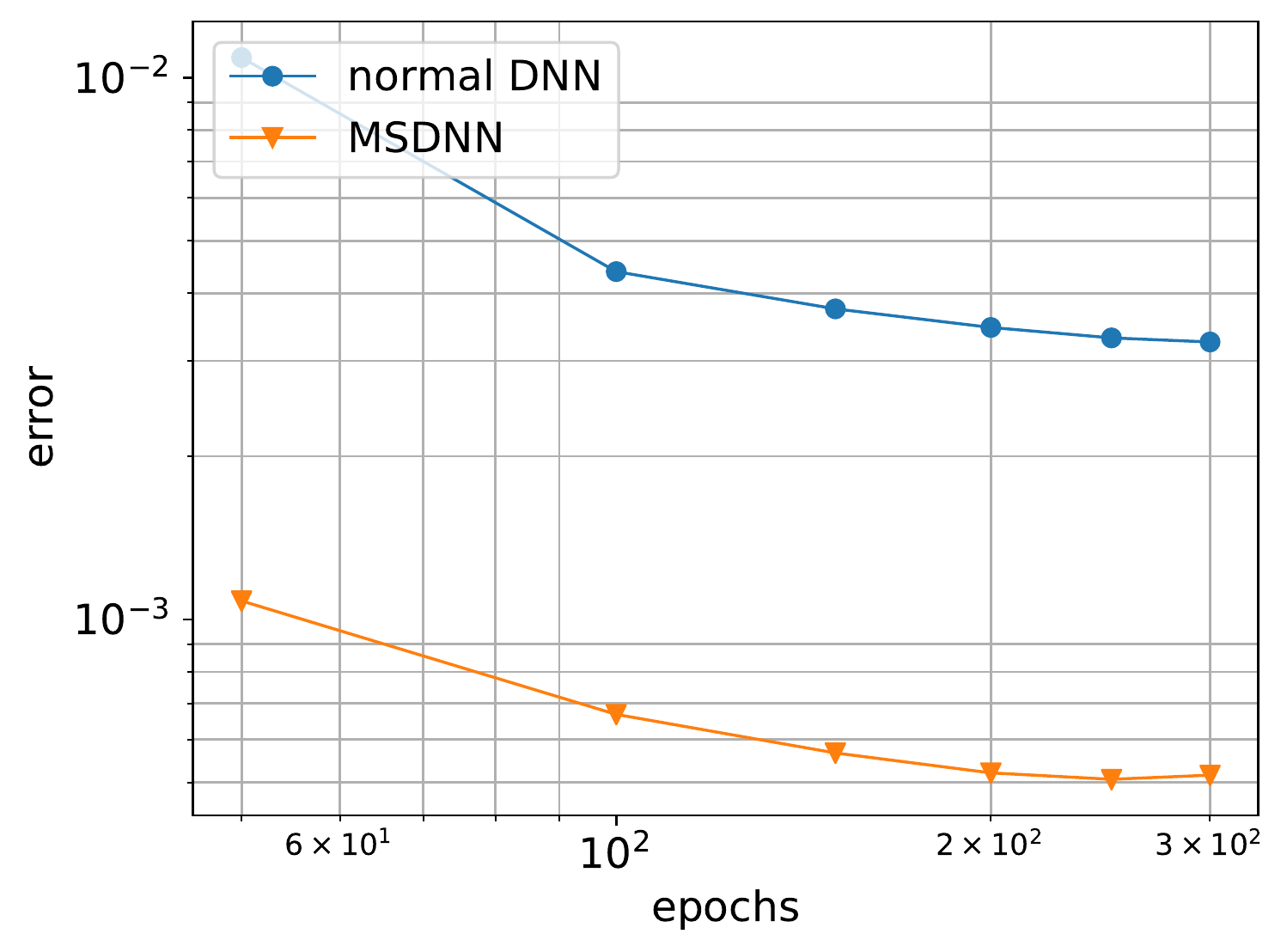}}
	\subfigure[$Err(p)$]{\includegraphics[scale=0.3]{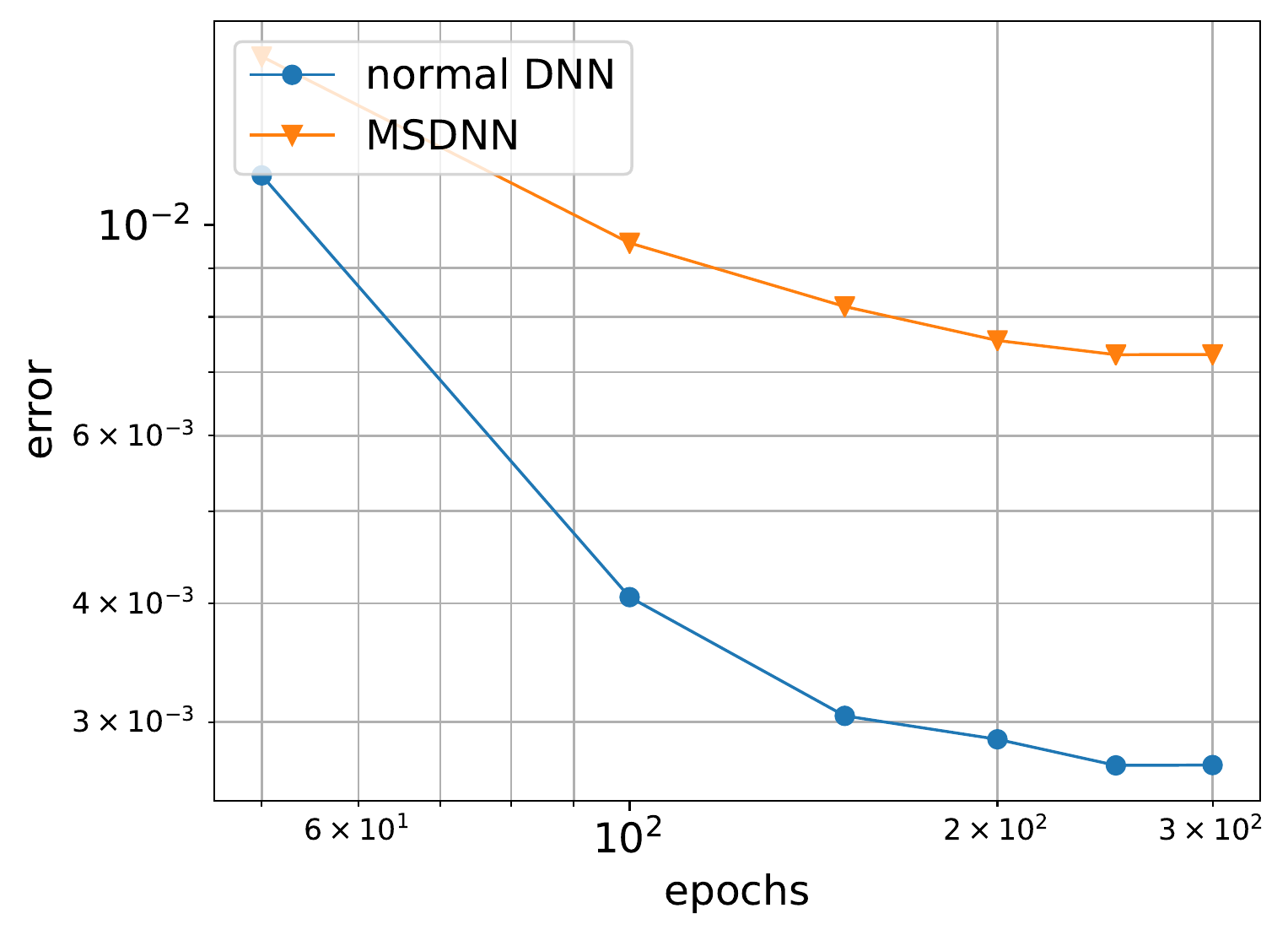}}
	\vspace{-10pt}
	\caption{Normal DNN and MscaleDNN with loss function $\bs L_{VgVP}(\theta_{\bs u},\theta_p, \theta_{\bs U},\theta_{\bs q})$.}%
	\label{example1-5}%
\end{figure}
\begin{figure}[ht!]
	\center
	\subfigure[$u_x$]{\includegraphics[scale=0.3]{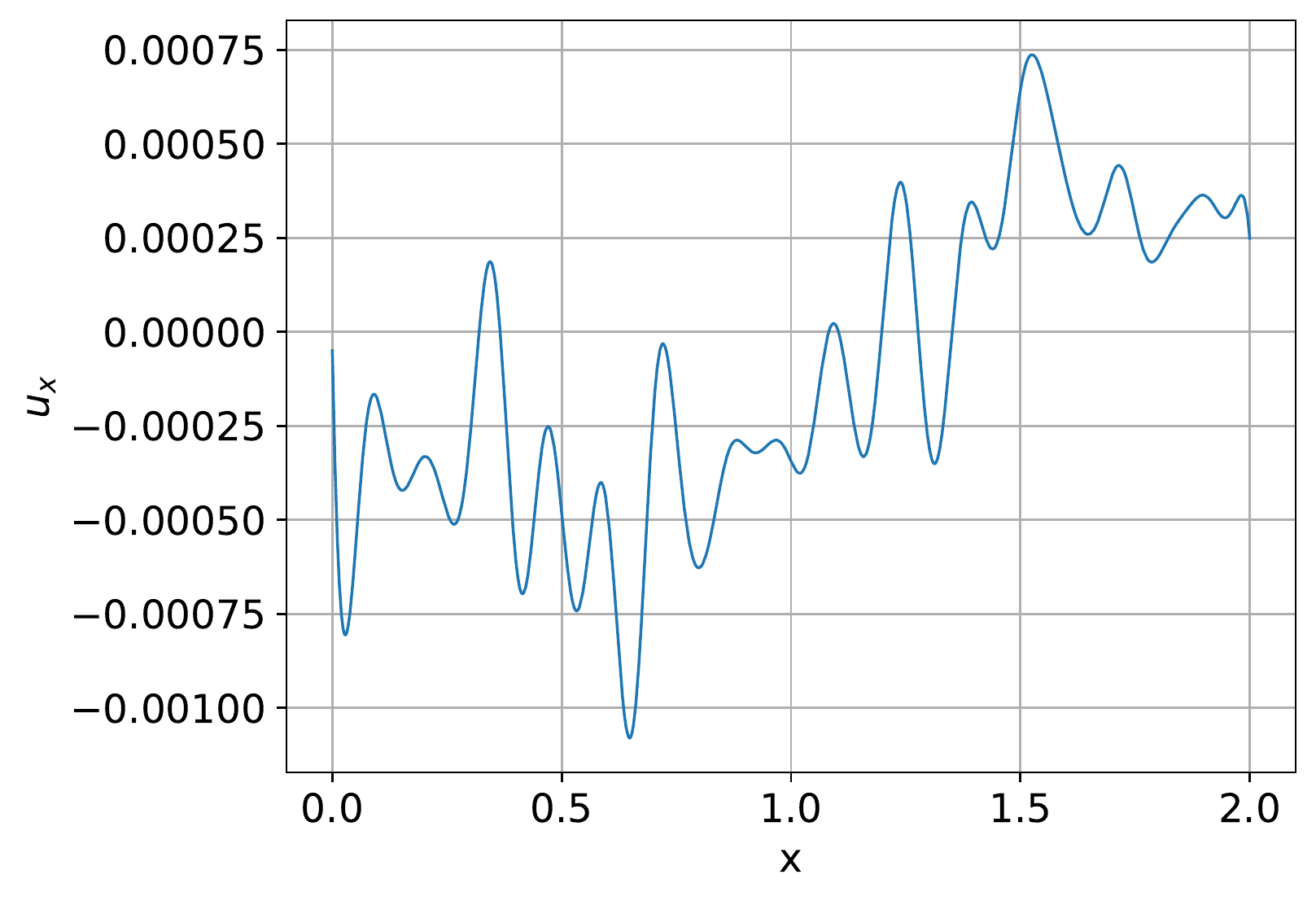}}
	\subfigure[$u_y$]{\includegraphics[scale=0.3]{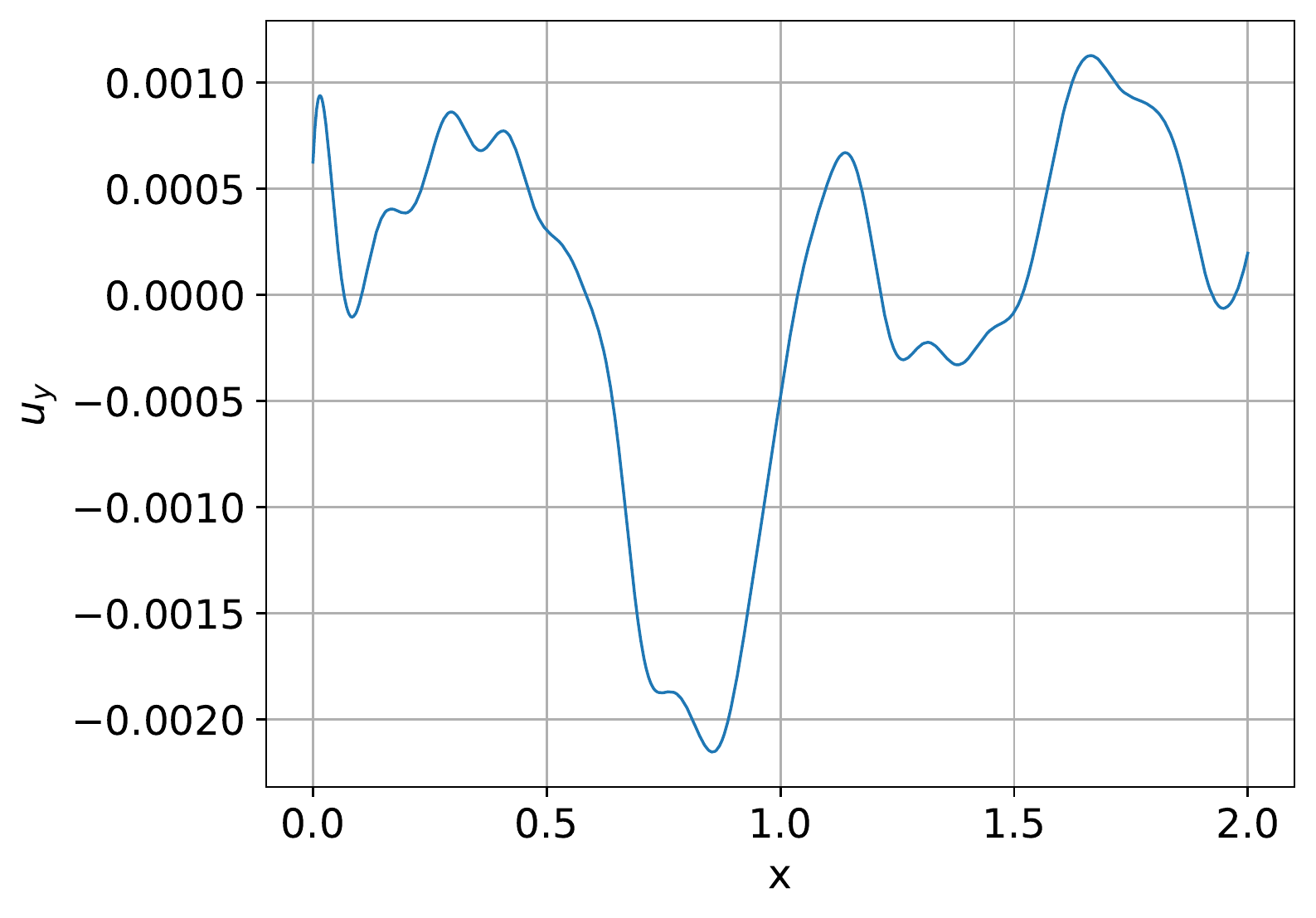}}
	\subfigure[$p$]{\includegraphics[scale=0.3]{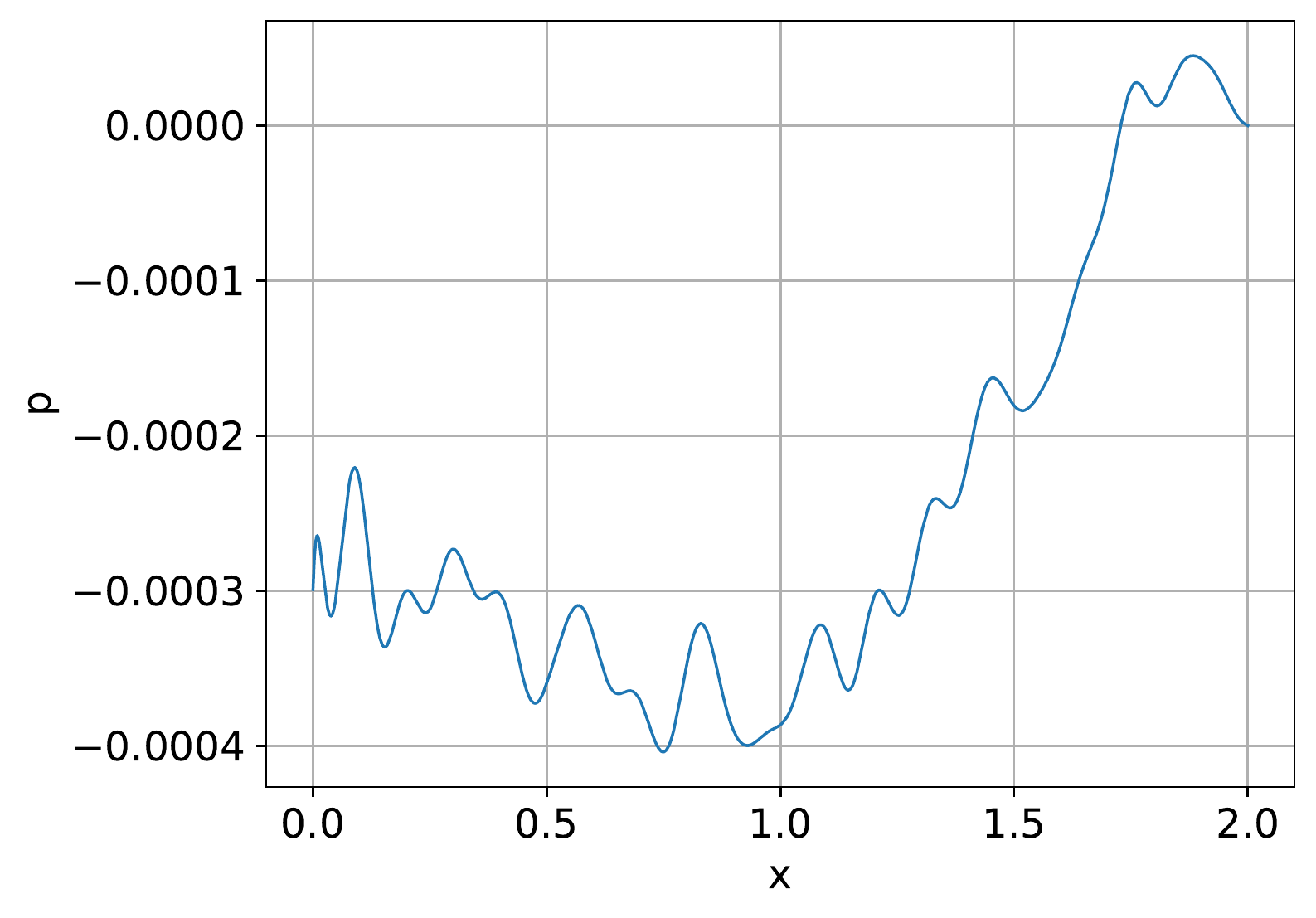}}
	\vspace{-10pt}
	\caption{Error of MscaleDNN solutions at epoch $300$ with loss function $\bs L_{{\bs\omega} VP}(\theta_{\bs u},\theta_p, \theta_{\bs \omega},\theta_{\bs q})$.}%
	\label{example1-6}%
\end{figure}

More detailed difference can been seen from the comparison of loss and errors between the normal DNN and MscaleDNN for the three loss functions in Fig. \ref{example1-3}-\ref{example1-5}. The results show that the MscaleDNNs have much faster convergence no matter which loss function is used. In fact, MscaleDNNs can achieve much better accuracy than normal DNN as we can see in Fig. \ref{example1-3}(b)-\ref{example1-5}(b). In particular,  the MscaleDNN solutions obtained by minimizing the $\mbox{$\omega$}VP$-loss are compared with exact solution along the line $y=0.7$ in Fig. \ref{example1-6}. It is clear that the MscaleDNN solutions match very well with the exact solutions.

\section{Oscillatory Kovasznay flows in a domain with multiple cylindrical voids}
\begin{figure}[htbp]
\centering
\subfigure[computational domain]{\includegraphics[width=0.4\linewidth]{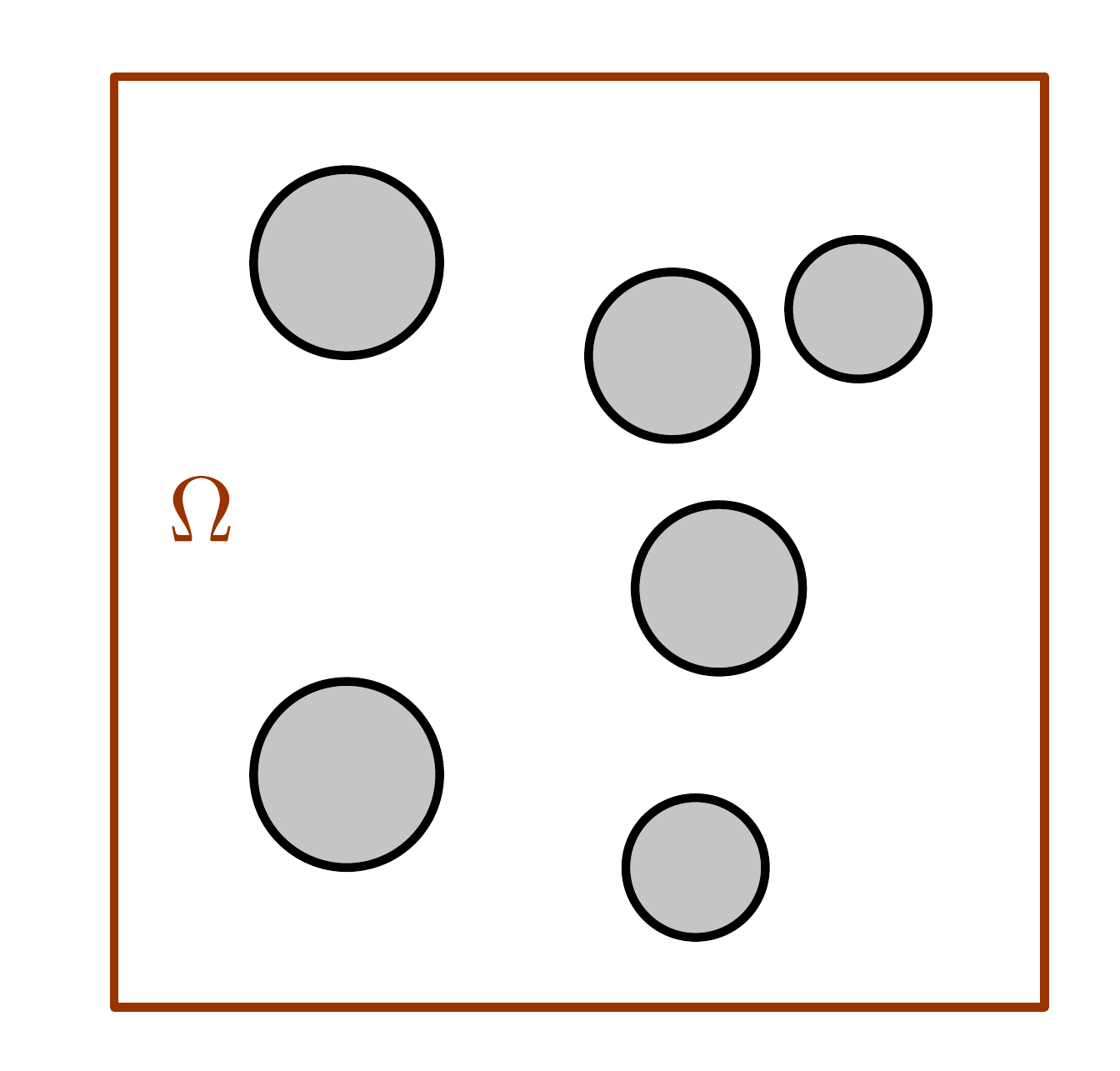}}\qquad
\subfigure[an example of $u_1$ in \eqref{highoscillatedflowmixed}]{\includegraphics[width=0.5\linewidth]{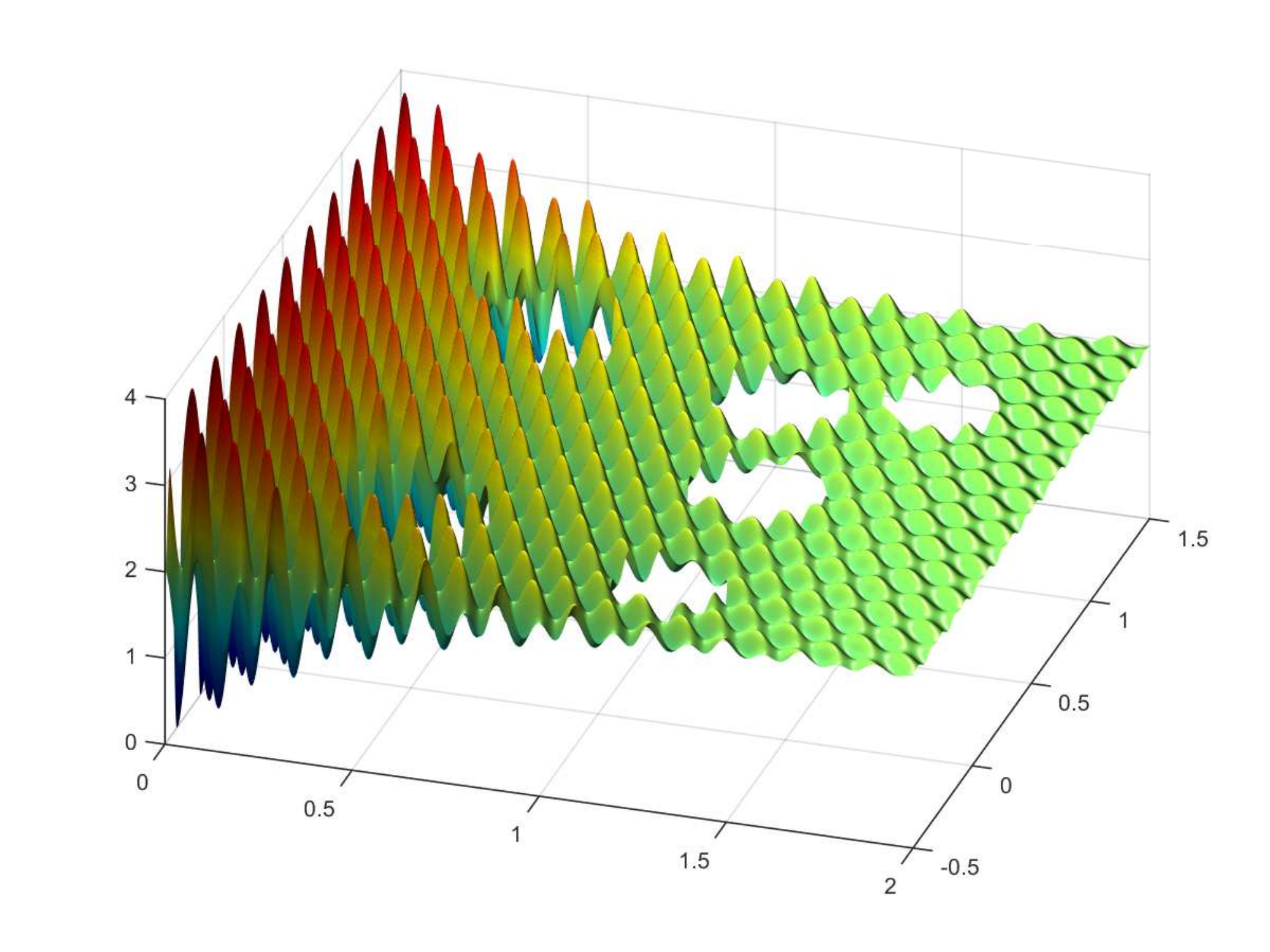}}
\vspace{-10pt}
\caption{A oscillatory solution over a domain with six cylindrical voids.}
\label{domain6holes}
\end{figure}

The MscaleDNN is more powerful than a normal DNN due to the former's capability on solving complicate problems with oscillatory solutions. Here, we consider the Stokes flow in the domain $\Omega=[0, 2]\times[-0.5,1.5]$ with $6$ cylindrical holes (refer to Fig. \ref{domain6holes}) centered at
$$(0.5, 0.0), \;\;(1.25,-0.2),\;\; (1.3,0.4), \;\; ( 0.5, 1.1), \;\; (1.2, 0.9),\;\; (1.6, 1),$$
inside the domain. The radius of the cylinders are set to be $0.2, 0.15, 0.18, 0.2, 0.18, 0.15$, respectively.  We will test two  exact solutions with highly oscillatory velocity fields. All examples are set to run $1500$ epochs using Adam.

The adaptive learning rates technique will be used in the numerical tests below, where the learning rate of the first $500$ epochs is set to be $0.001$, then, the learning rate will be reduced by a factor of 10 after each $500$ epochs. The change of learning rate can be seen clearly in the history of losses later. In the loss functions, we fix the penalty parameter $\beta=100$ and set an initial penalty parameter $\alpha=2000$. Every 50 epochs, we check the errors $Err(\bs u)$ and $Err(p)$ and adjust parameter $\alpha$ as follows
\begin{itemize}
	\item If $Err(\bs u)>2Err(p)$, $\alpha=\alpha+500$;
	\item If $Err(p)>2Err(\bs u)$ and $\alpha>500$, $\alpha=\alpha-500$.
\end{itemize}

In the results below, the $\ell^2$-errors defined in \eqref{ell2err} are computed again with 34,072 randomly selected points in the computational domain.

\subsection{Two frequency solution}

The first case has an exact solution given by
\begin{equation}\label{highoscillatedflow}
\begin{split}
u_1=&1-e^{\lambda x_1}\cos(2n\pi x_1+2m\pi x_2)),\\
u_2=&\frac{\lambda}{2m\pi}e^{\lambda x_1}\sin(2n\pi x_1+2m\pi x_2)+\frac{n}{m}e^{\lambda x_1}\cos(2n\pi x_1+2m\pi x_2),\\
p=&\frac{1}{2}e^{2\lambda x_1},\quad \lambda=\frac{Re}{2}-\sqrt{\frac{Re^2}{4}+4\pi^2},\quad Re=\frac{1}{\nu},
\end{split}
\end{equation}
with frequencies $n=50, m=55$. In the simulations of this example, the MscaleDNNs for $\bs u$, $\bs \omega$, $\bs T$ and $\bs U$ are set to have $11$
scales: $\{\bs x, 2\bs x, \cdots, 2^{10}\bs x\}$ and the embedded fully connected DNN for each scale is set to have $8$ hidden layers and $150$ neurons in each hidden layer. As the pressure does not have high oscillations, the MscaleDNNs for $p$ and $\bs q$ are set to have $6$
scales: $\{\bs x, 2\bs x, \cdots, 2^{5}\bs x\}$ and the embeded fully connected DNN for each scale is set to have $8$ hidden layers and $50$ neurons in each hidden layer. We randomly sample 850621 points inside $\Omega$ and 140000 points on the boundary for learning. In the learning process, we set batch size equal to 10000 points inside the domain and randomly pick 2000 points on the boundary for each step.

\begin{figure}[ht!]
	\center
	\includegraphics[scale=0.28]{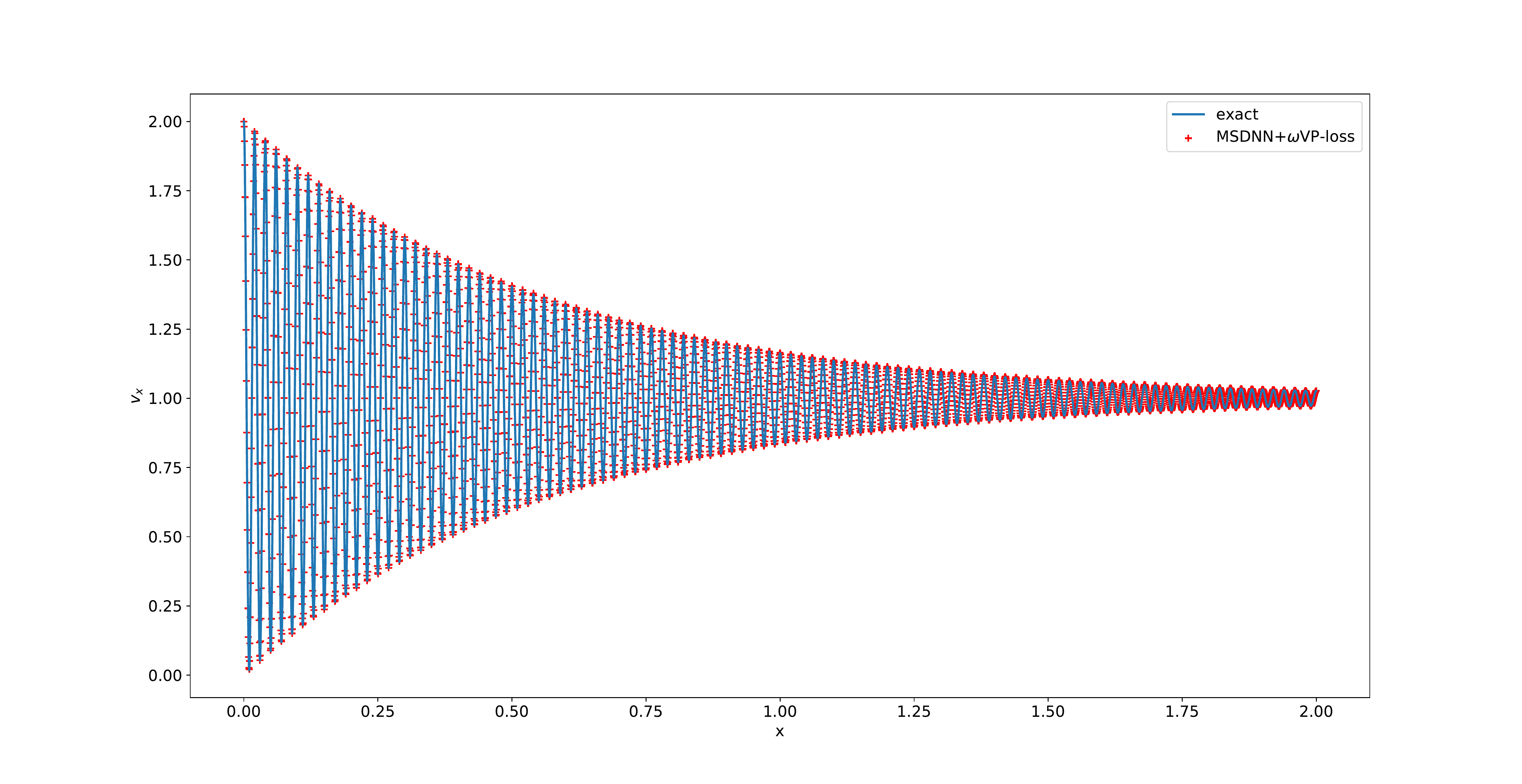}
	\vspace{-10pt}
	\caption{Exact $u_1$ and its MscaleDNN approximation with $\omega$VP-loss $\bs L_{\mbox{$\omega$}VP}(\theta_{\bs u},\theta_p, \theta_{\bs \omega},\theta_{\bs q})$.}
	\label{example2-1}%
\end{figure}

\begin{figure}[ht!]
	\center
	\includegraphics[scale=0.28]{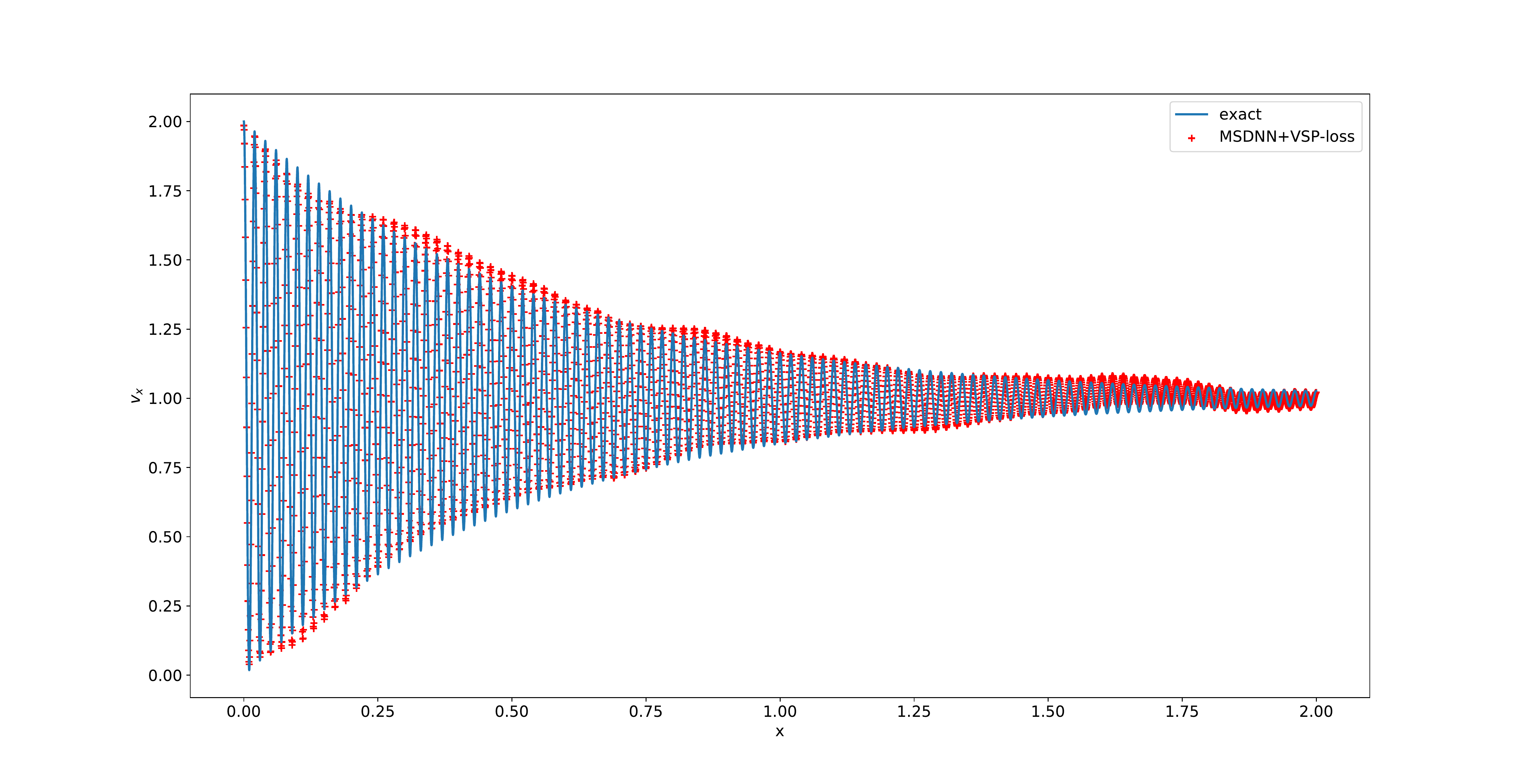}
	\vspace{-10pt}
	\caption{Exact $u_1$ and its MscaleDNN approximation with VSP-loss $\bs L_{VSP}(\theta_{\bs u},\theta_p, \theta_{\bs T},\theta_{\bs q})$.}
	\label{example2-2}%
\end{figure}

\begin{figure}[ht!]
	\center
	\includegraphics[scale=0.28]{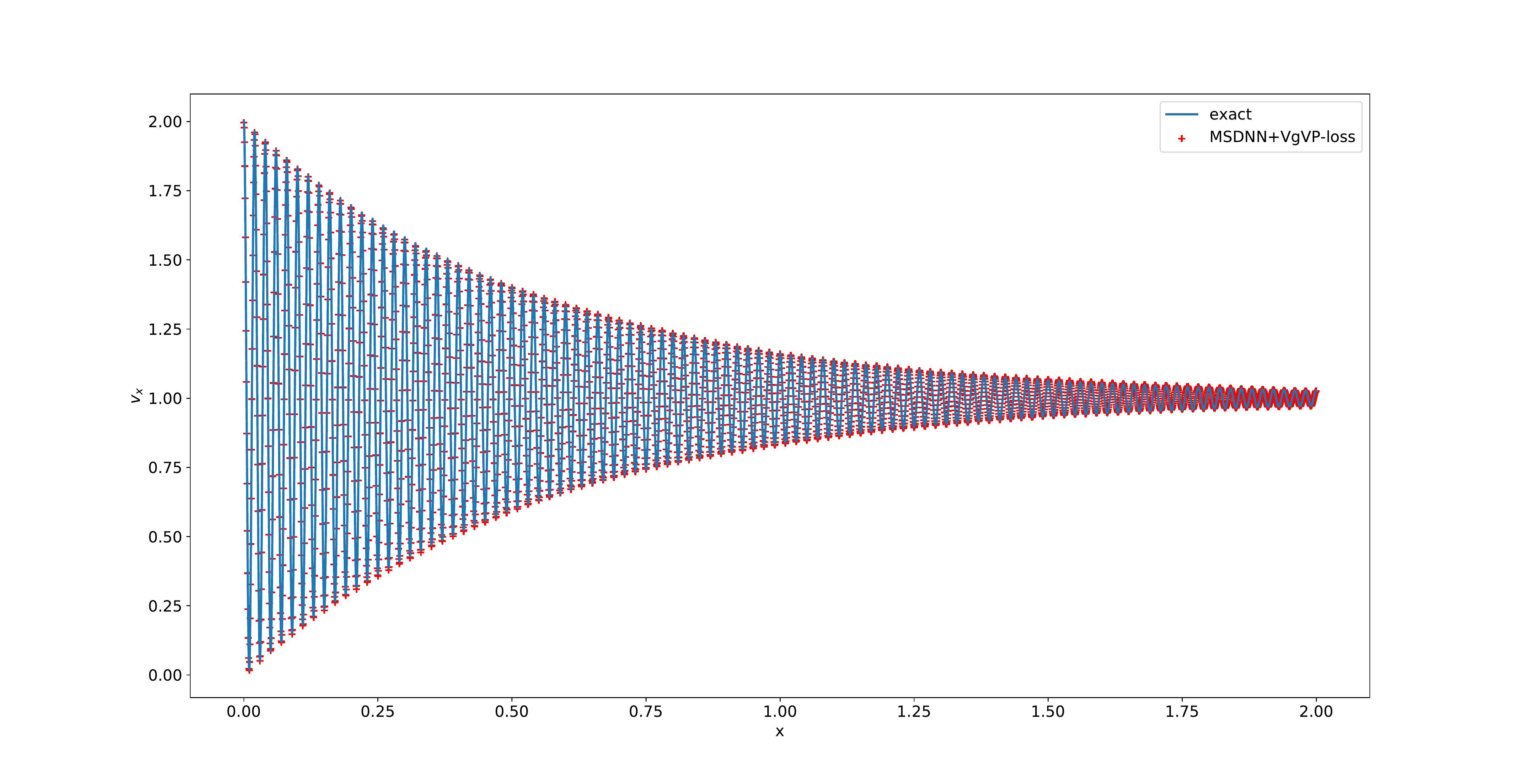}
	\vspace{-10pt}
	\caption{Exact $u_1$ and its MscaleDNN approximation with VgVP-loss $\bs L_{VgVP}(\theta_{\bs u},\theta_p, \theta_{\bs U},\theta_{\bs q})$.}
	\label{example2-3}%
\end{figure}

\begin{figure}[ht!]
	\center
	\subfigure[$Err(\bs u)$]{\includegraphics[scale=0.4]{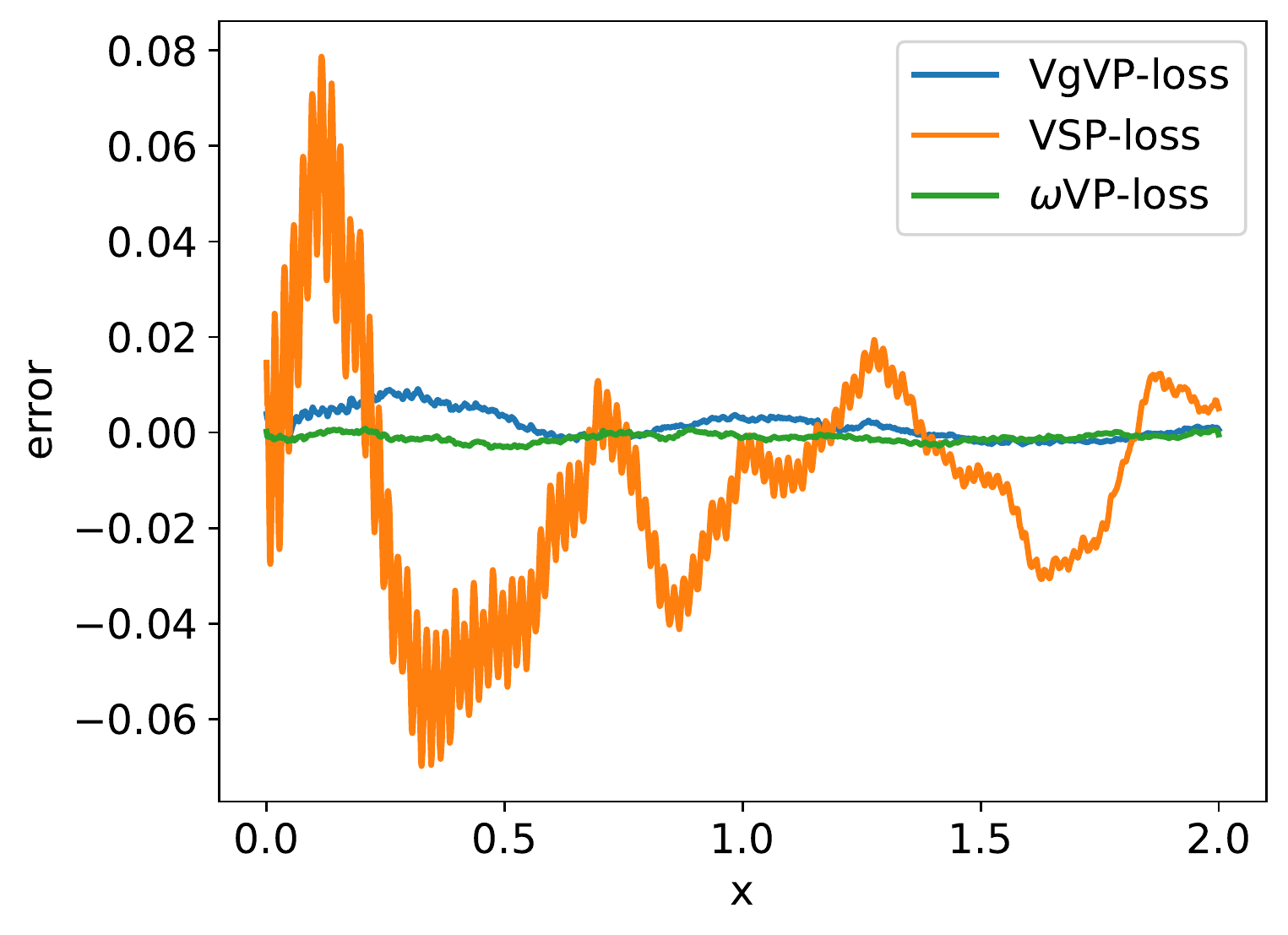}}\qquad
	\subfigure[$Err(p)$]{\includegraphics[scale=0.4]{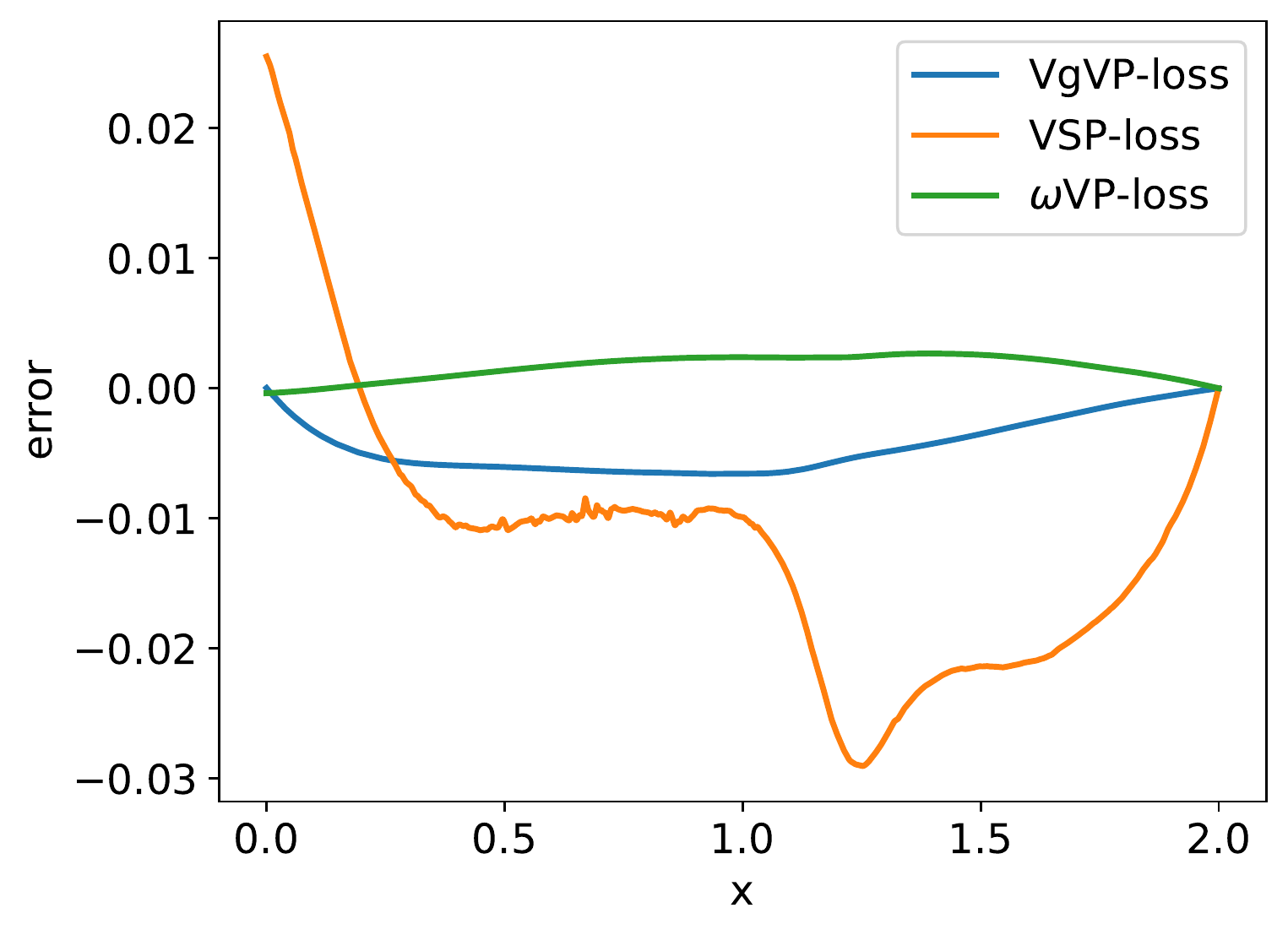}}
	\vspace{-10pt}
	\caption{Errors of MscaleDNN approximations using different loss functions.}
	\label{example2-4}%
\end{figure}

The MscaleDNN solutions of $u_1$ are compared with the exact $u_1$ in Fig. \ref{example2-1}-\ref{example2-3}. Here, we plot the solutions along the line $y=0.7$ which does not cross any of cylinders inside the domain. Errors of the MscaleDNN approximations for $\bs u$ and $p$ using different losses are depicted in Fig. \ref{example2-4}.
We can see that the $\omega$VP-loss or $VgVP$-loss with MscaleDNN can produce very accurate solutions in just 1500 epochs while the VSP-loss needs more learning to achieve similar accuracy.
\vspace{-10pt}
\vspace{-20pt}

\begin{figure}[ht!]
	\center
	\subfigure[loss]{\includegraphics[scale=0.3]{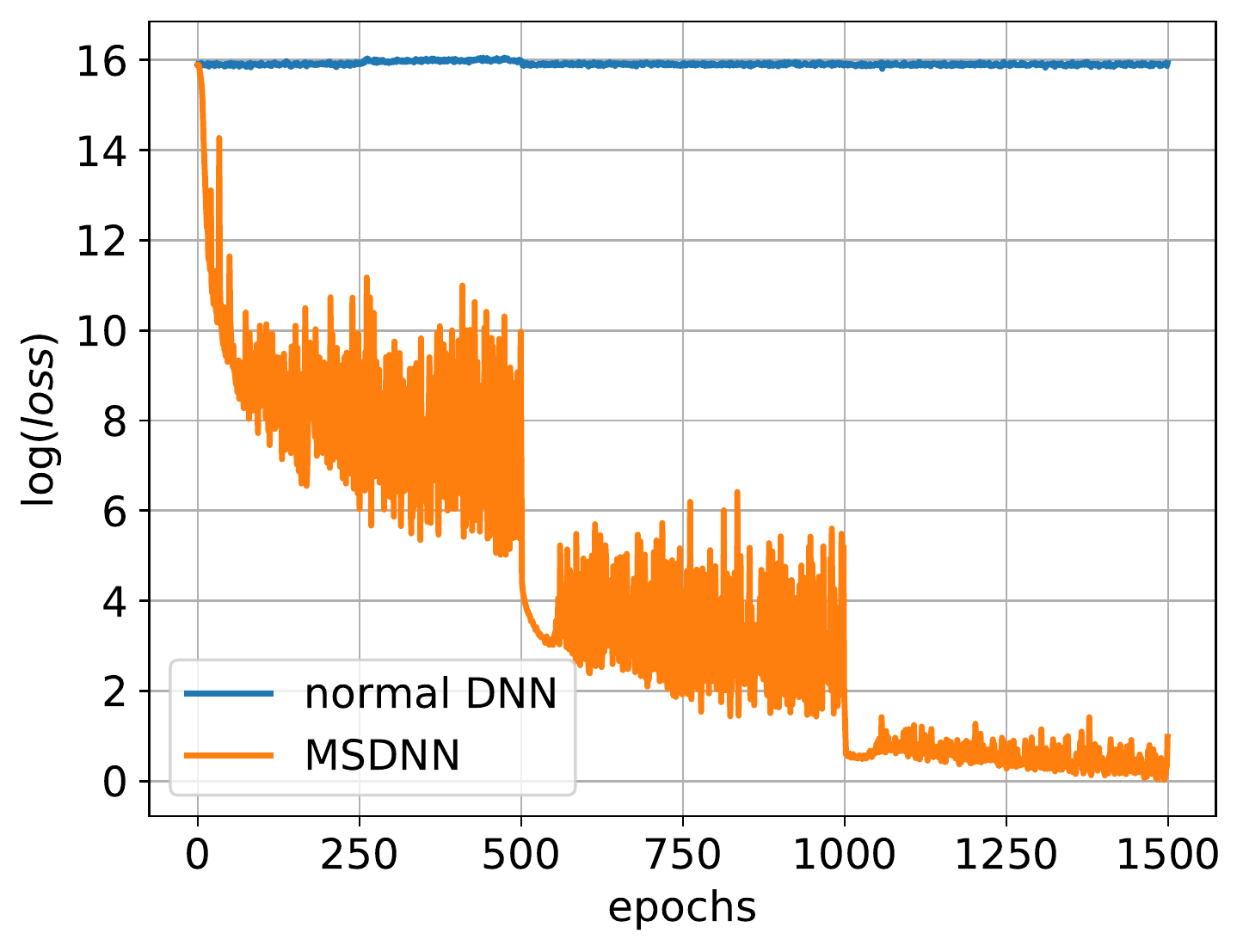}}
	\subfigure[$Err(\bs u)$]{\includegraphics[scale=0.3]{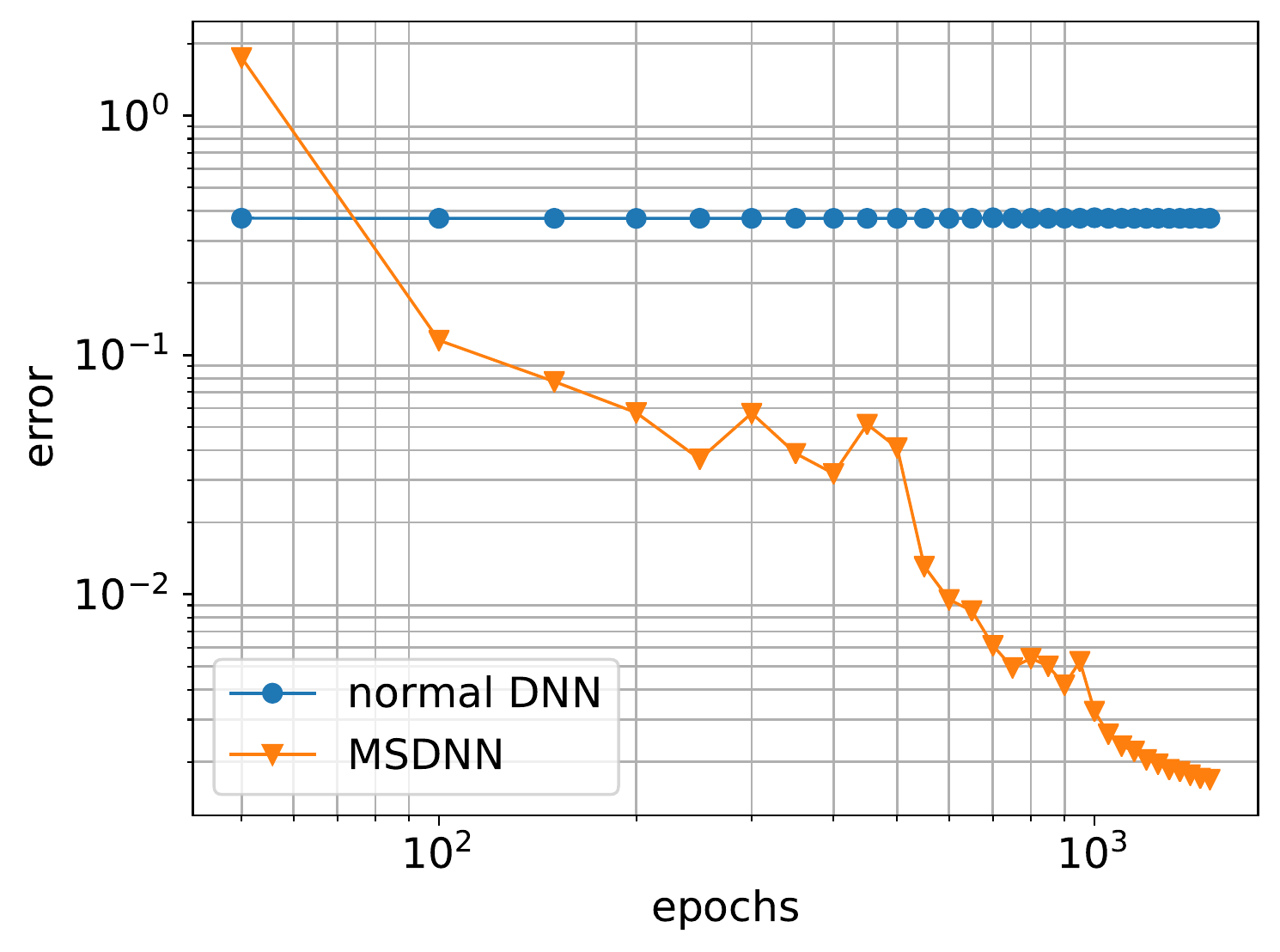}}
	\subfigure[$Err(p)$]{\includegraphics[scale=0.3]{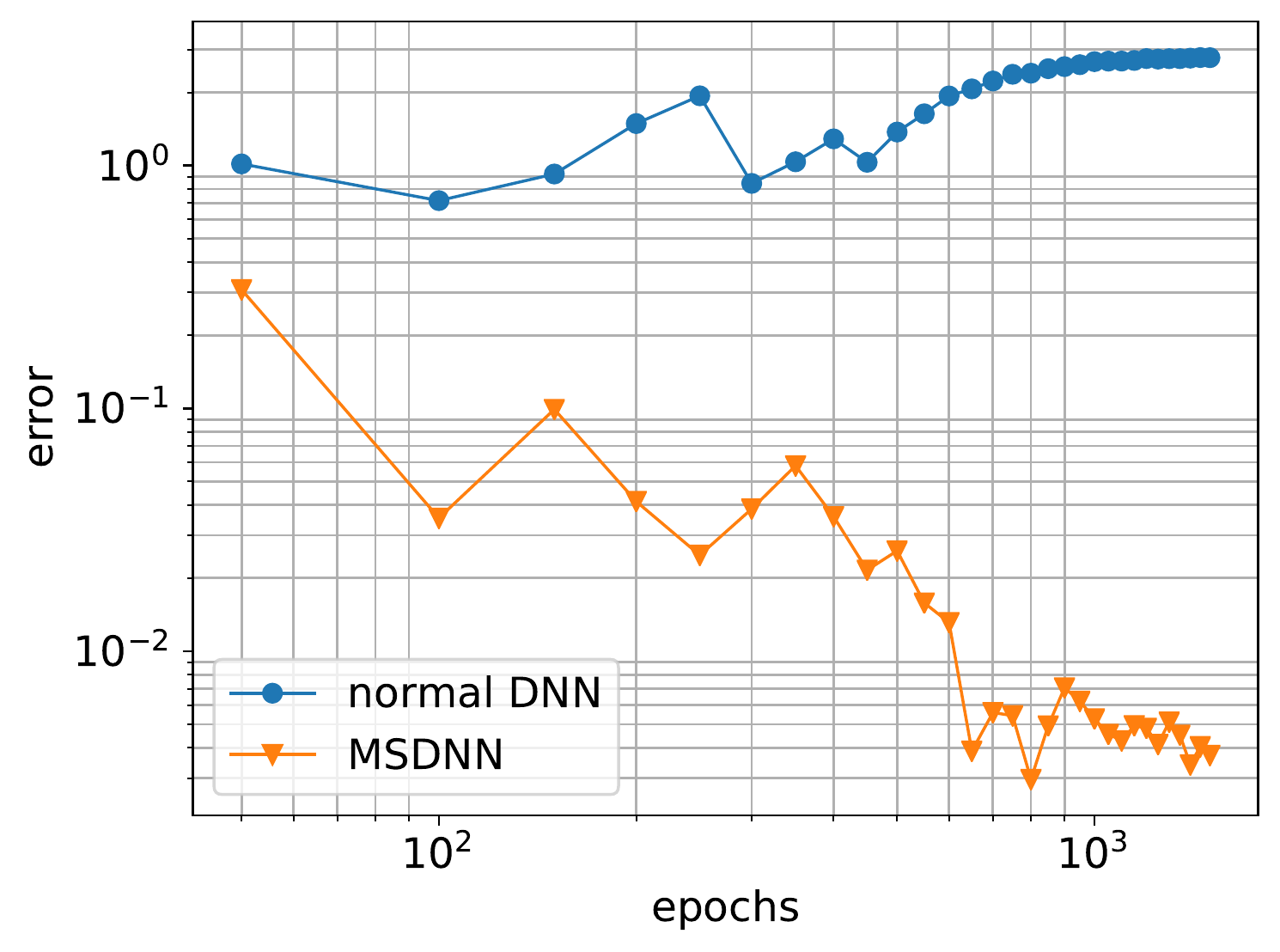}}
	\vspace{-10pt}
	\caption{Comparison of a normal DNN and the MscaleDNN with loss function  $\bs L_{\mbox{$\omega$}VP}(\theta_{\bs u},\theta_p, \theta_{\bs \omega},\theta_{\bs q})$.}
	\label{example2-5}%
\end{figure}

\begin{figure}[ht!]
	\center
	\subfigure[loss]{\includegraphics[scale=0.3]{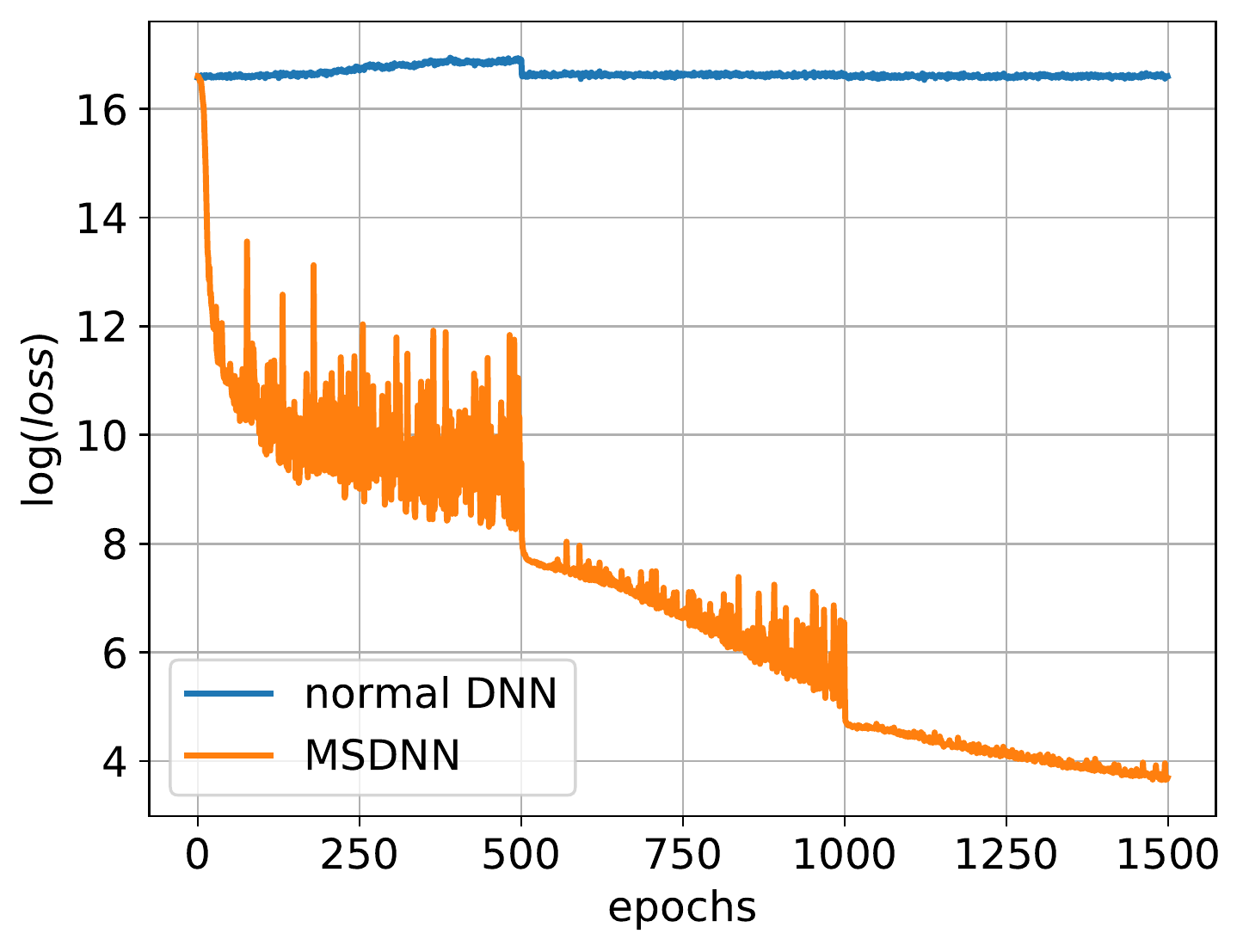}}
	\subfigure[$Err(\bs u)$]{\includegraphics[scale=0.3]{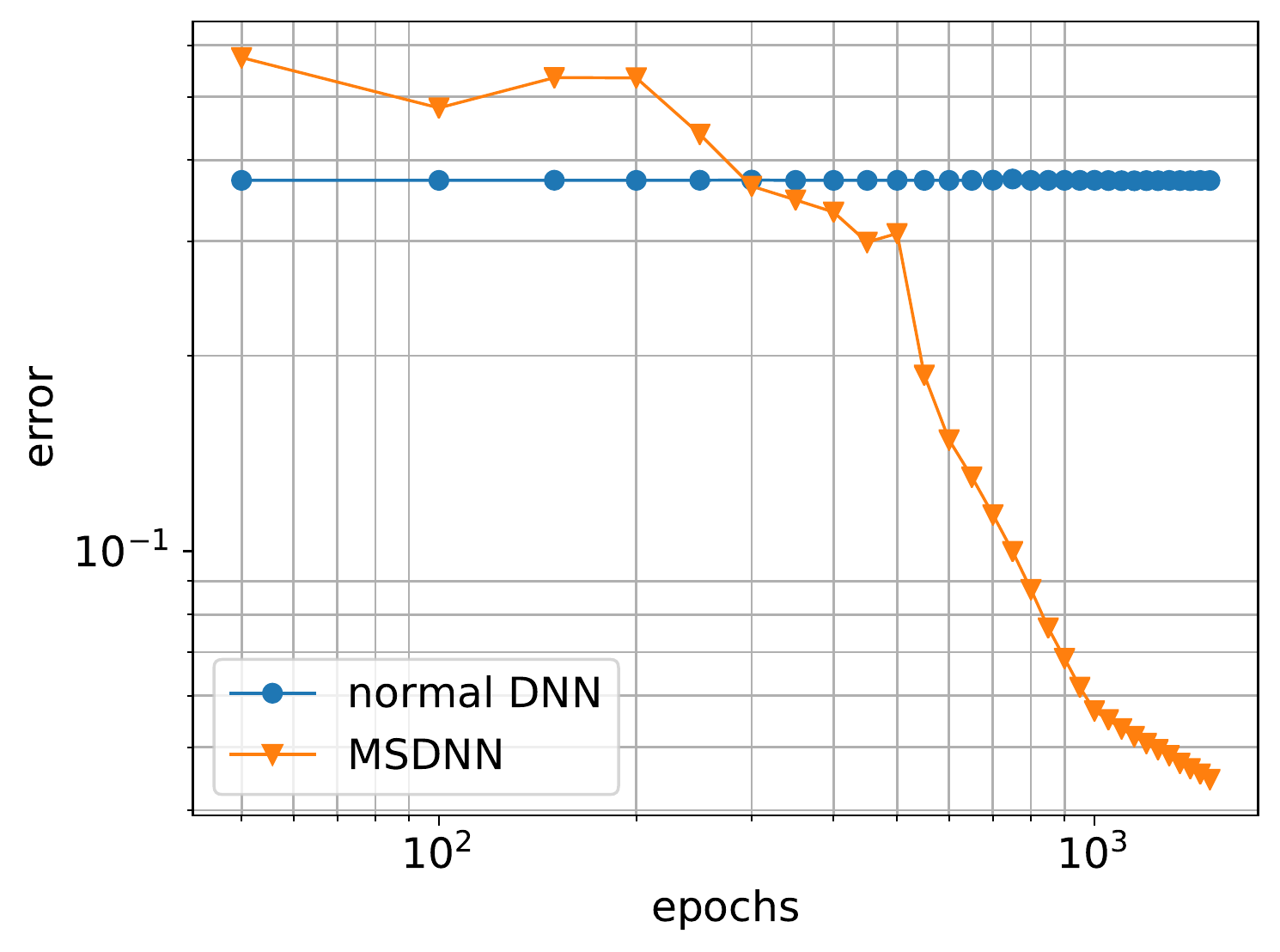}}
	\subfigure[$Err(p)$]{\includegraphics[scale=0.3]{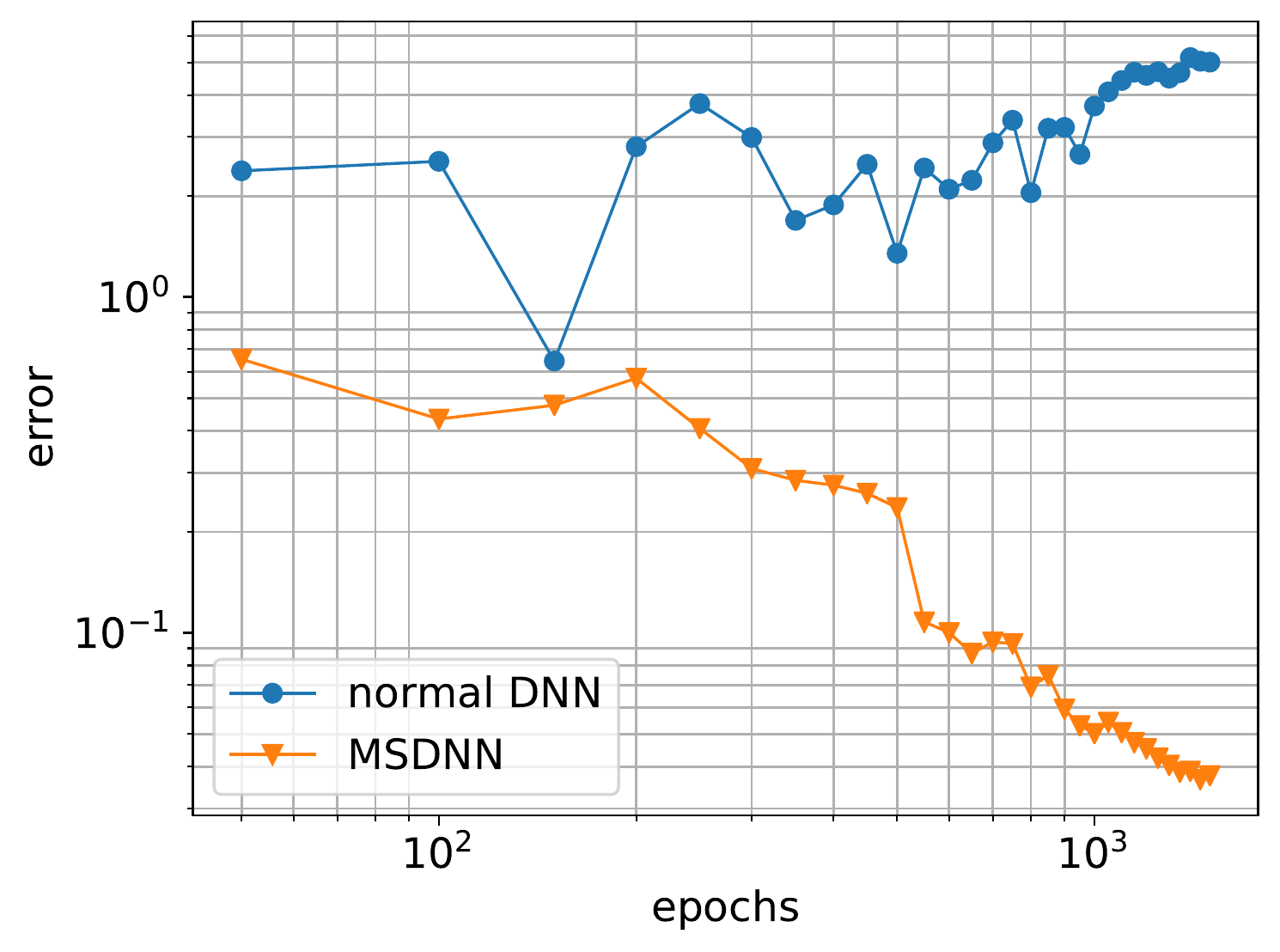}}
	\vspace{-10pt}
	\caption{Comparison of a normal DNN and the MscaleDNN with loss function  $\bs L_{VSP}(\theta_{\bs u},\theta_p, \theta_{\bs T},\theta_{\bs q})$.}
	\label{example2-6}%
\end{figure}

\begin{figure}[ht!]
	\center
	\subfigure[loss]{\includegraphics[scale=0.3]{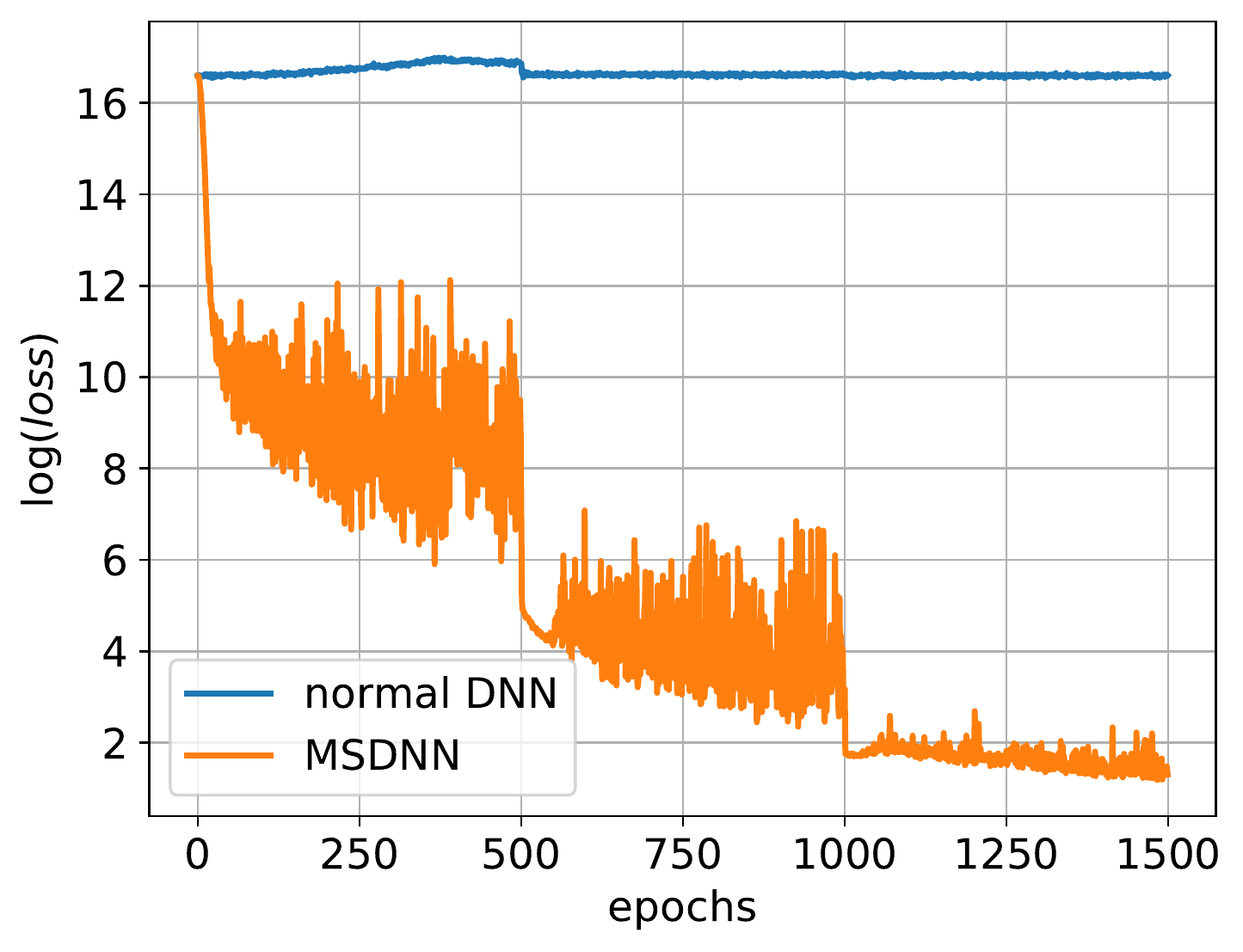}}
	\subfigure[$Err(\bs u)$]{\includegraphics[scale=0.3]{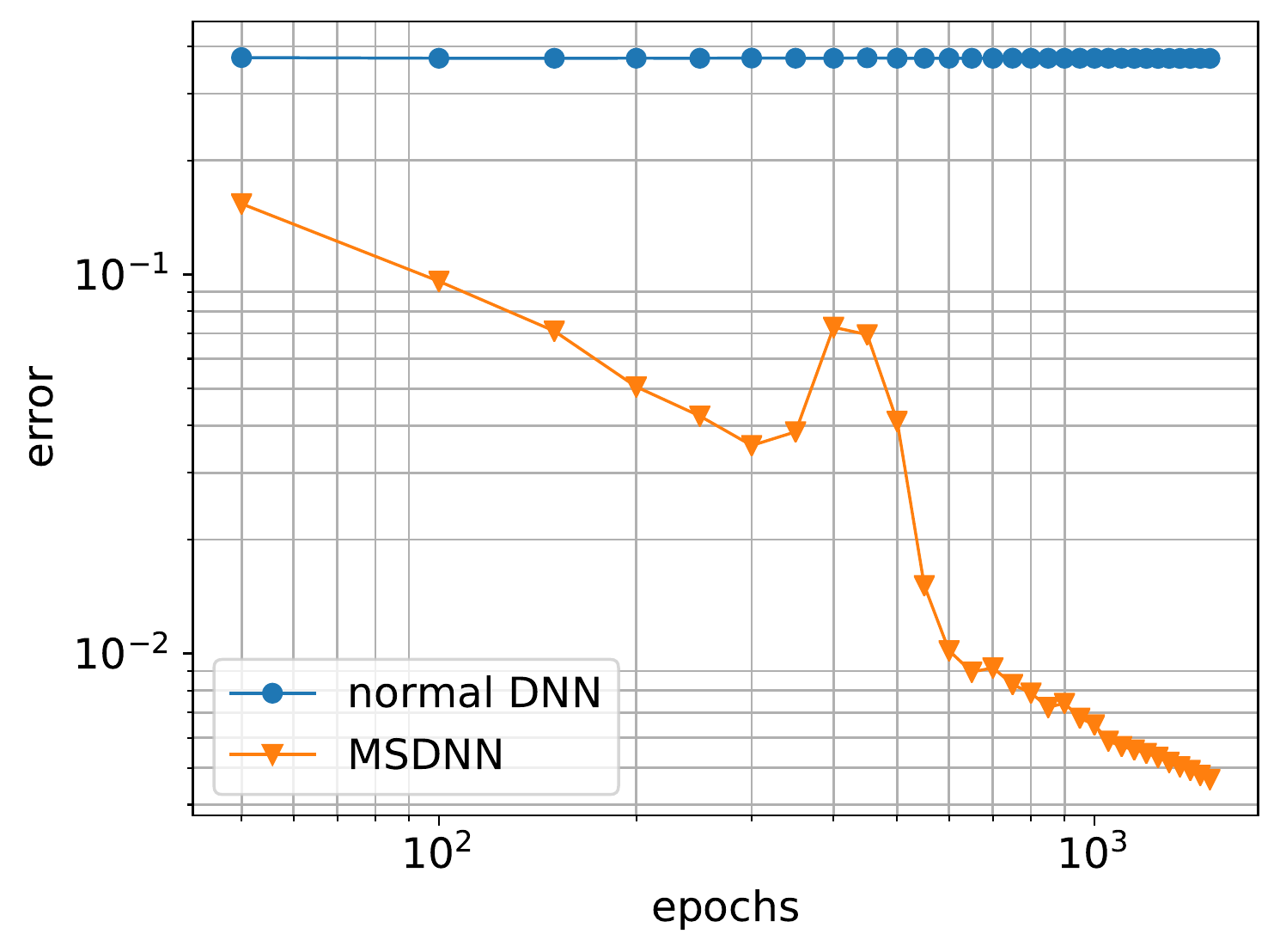}}
	\subfigure[$Err(p)$]{\includegraphics[scale=0.3]{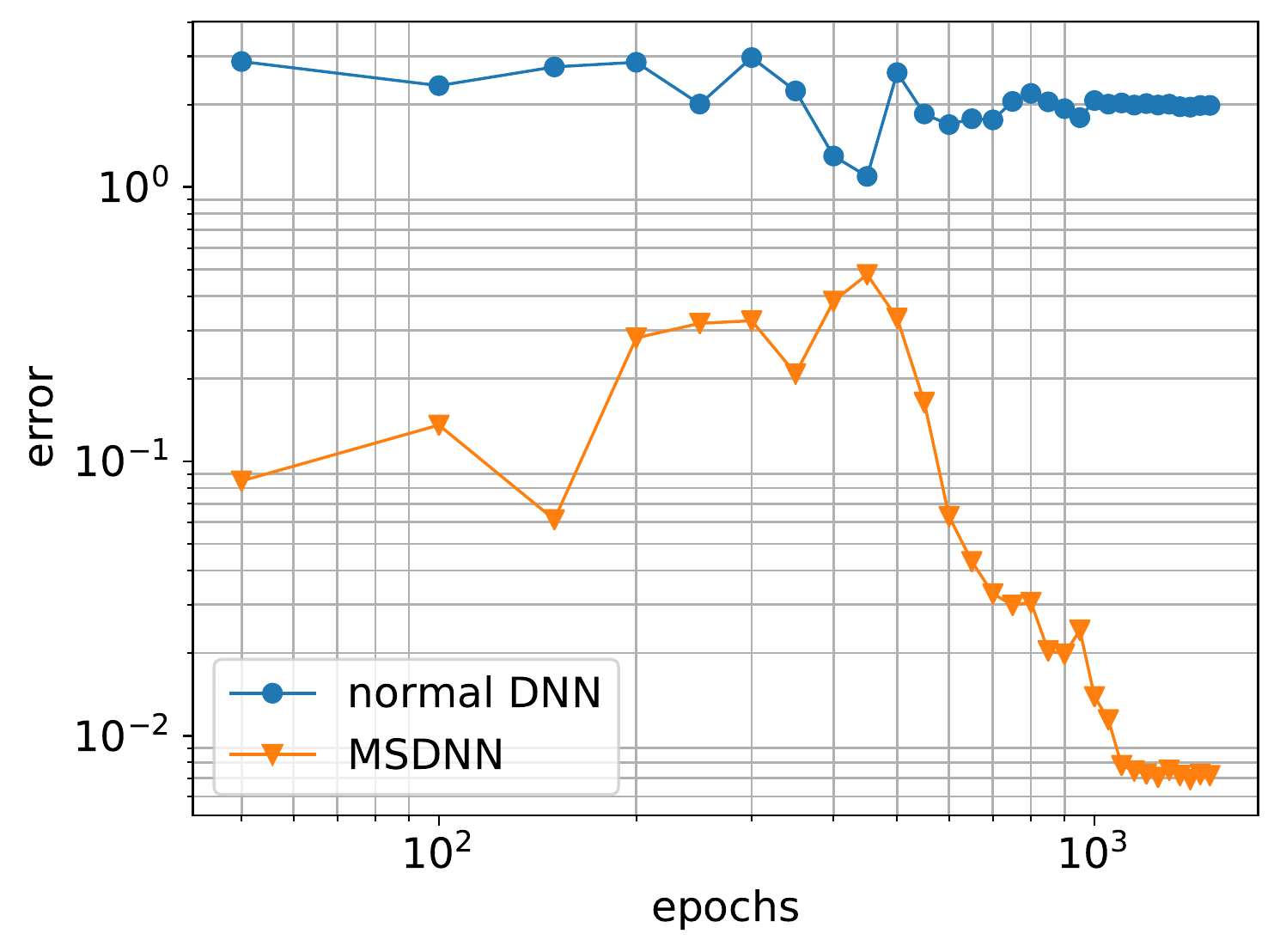}}
	\vspace{-10pt}
	\caption{Comparison of a normal DNN and the MscaleDNN with loss function  $\bs L_{VgVP}(\theta_{\bs u},\theta_p, \theta_{\bs U},\theta_{\bs q})$.}
	\label{example2-7}%
\end{figure}

For comparison, we also test the DNN-based algorithm only using fully connected DNNs. For $\bs u$ and intermediate variables $\bs\omega$, $\bs T$ and $\bs U$, we use fully connected DNNs with $8$ hidden layers and $1650$ neurons in each hidden layer. For $p$ and $\bs q$, we use fully connected DNNs with $8$ hidden layers and $300$ neurons in each hidden layer. Therefore, the total number of neurons in the fully connected DNNs and the MscaleDNNs are the same. The losses and $\ell^2$-errors obtained by minimizing different loss functions in \eqref{firstorderloss} are compared in Fig.  \ref{example2-5}-\ref{example2-7}. For this highly oscillatory solution, algorithms using fully connected DNNs can not learn anything within 1500 epochs. However, the ones using MscaleDNNs converge very fast within 1500 epochs.

\subsection{Multiple frequency solution}
\label{section_multiF}

Our second test problem will be a case where the velocity field has multiple high frequencies as follows,
\begin{equation}\label{highoscillatedflowmixed}
\begin{split}
u_1=&2-e^{\lambda x_1}\cos(70\pi x_1+60\pi x_2))-e^{\lambda x_1}\cos(80\pi x_1+90\pi x_2)),\\
u_2=&\frac{\lambda}{60\pi}e^{\lambda x_1}\sin(70\pi x_1+60\pi x_2)+\frac{7}{6}e^{\lambda x_1}\cos(70\pi x_1+60\pi x_2)\\
+&\frac{\lambda}{90\pi}e^{\lambda x_1}\sin(80\pi x_1+90\pi x_2)+\frac{8}{9}e^{\lambda x_1}\cos(80\pi x_1+90\pi x_2),\\
p=&\frac{1}{2}e^{2\lambda x_1},\quad \lambda=\frac{Re}{2}-\sqrt{\frac{Re^2}{4}+4\pi^2},\quad Re=\frac{1}{\nu}.
\end{split}
\end{equation}

For this test, the MscaleDNNs for $\bs u$, $\bs \omega$, $\bs T$ and $\bs U$ are set to have $10$ scales: $\{\bs x, 2\bs x, \cdots, 2^{9}\bs x\}$ and the embedded fully connected DNN for each scale is set to have $8$ hidden layers and $120$ neurons in each hidden layer. As in the last numerical test, the MscaleDNNs for $p$ and $\bs q$ are set to have $6$ scales: $\{\bs x, 2\bs x, \cdots, 2^{5}\bs x\}$ and the embeded fully connected DNN for each scale is set to have $8$ hidden layers and $50$ neurons in each hidden layer. We randomly sample 425290 points inside $\Omega$ and 140000 points on the boundary for learning. In the learning process, we set batch size equal to 5000 points inside the domain and randomly select 2000 points on the boundary for each step.

\begin{figure}[ht!]
	\center
	\includegraphics[scale=0.28]{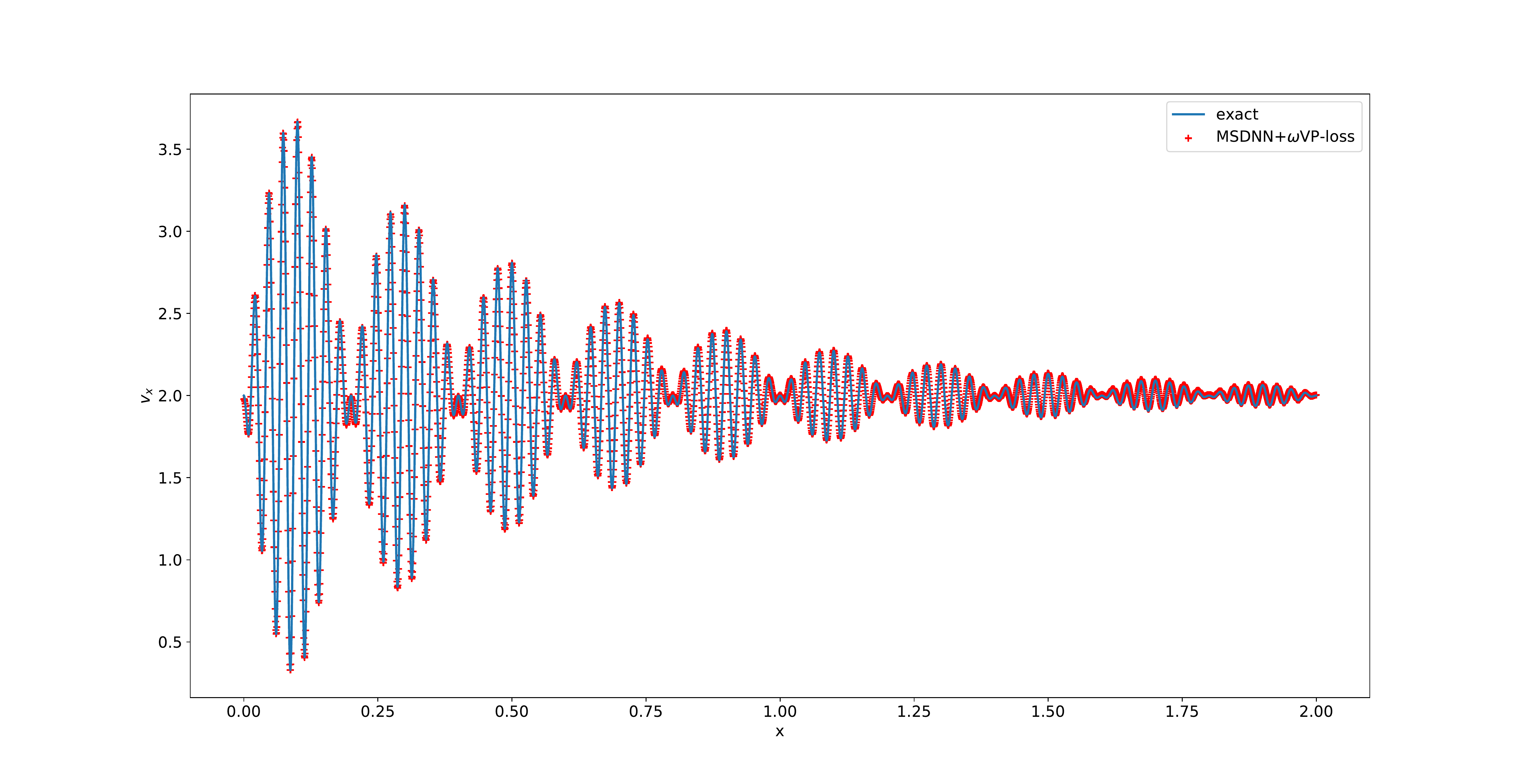}
	\vspace{-10pt}
	\caption{Exact $u_1$ and its MscaleDNN approximation with loss function $\bs L_{\mbox{$\omega$}VP}(\theta_{\bs u},\theta_p, \theta_{\bs \omega},\theta_{\bs q})$.}
	\label{example3-4}%
\end{figure}

\begin{figure}[ht!]
	\center
	\includegraphics[scale=0.28]{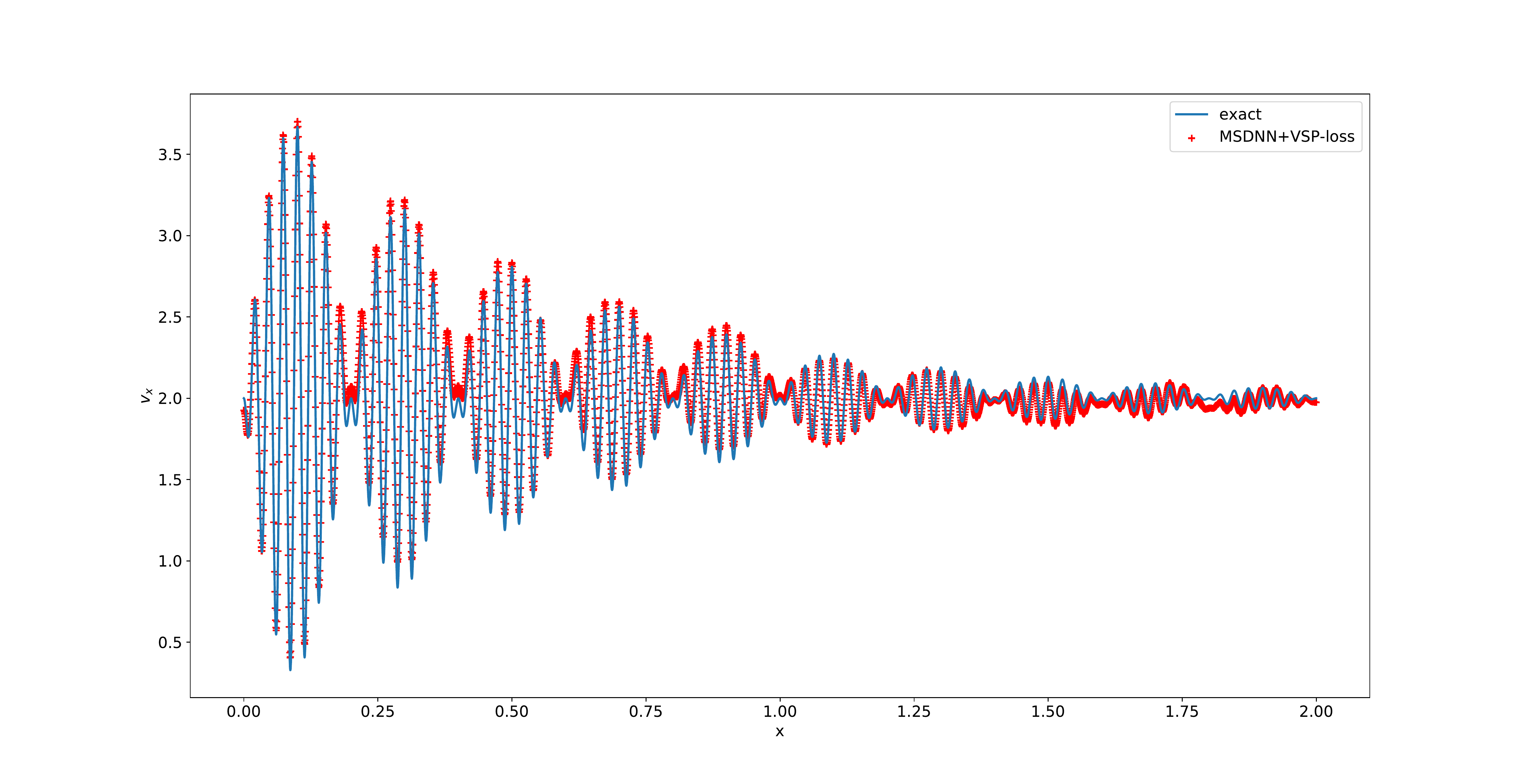}
	\vspace{-10pt}
	\caption{Exact $u_1$ and its MscaleDNN approximation with loss function $\bs L_{VSP}(\theta_{\bs u},\theta_p, \theta_{\bs T},\theta_{\bs q})$.}
	\label{example3-5}%
\end{figure}

\begin{figure}[ht!]
	\center
	\includegraphics[scale=0.28]{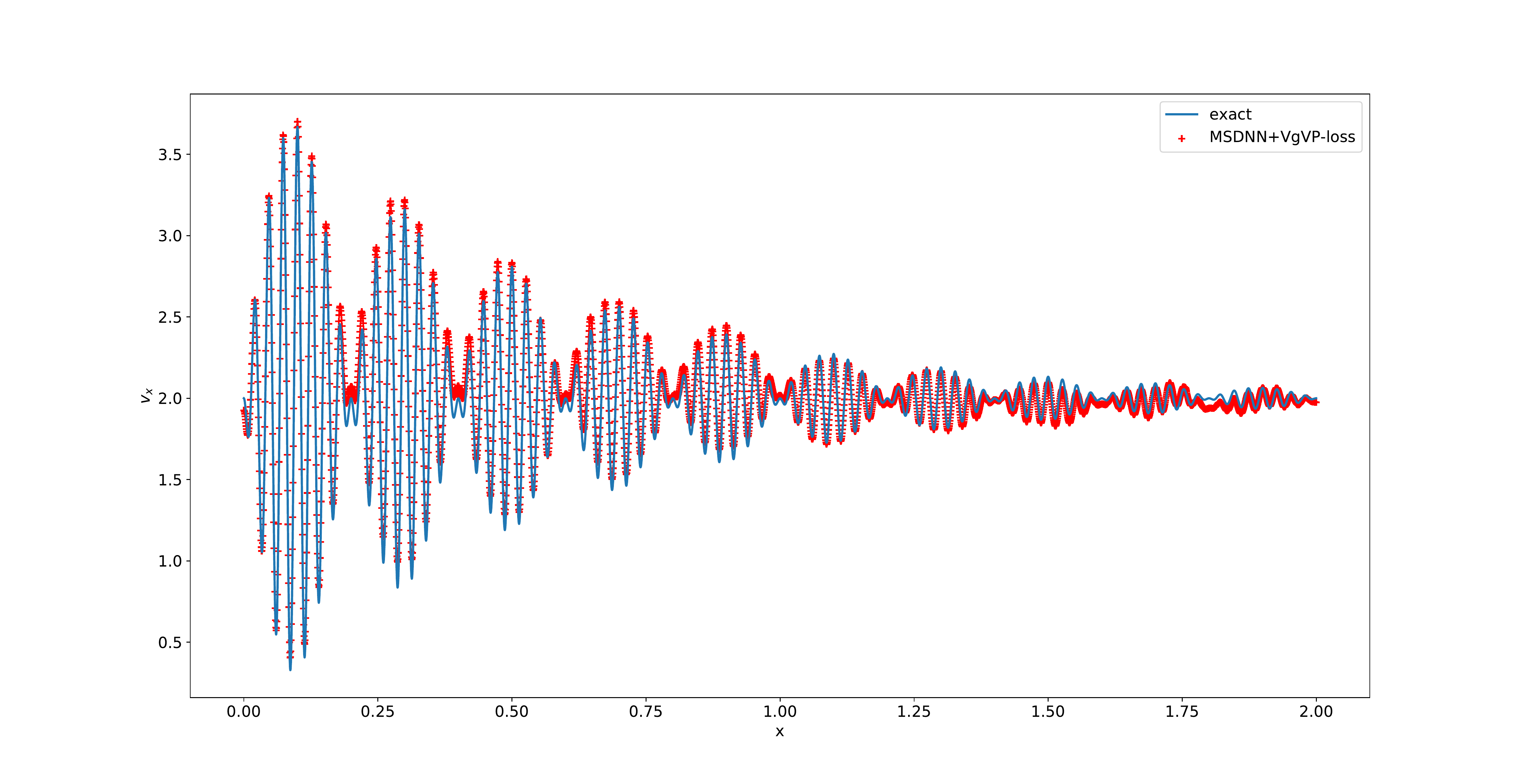}
	\vspace{-10pt}
	\caption{Exact $u_1$ and its MscaleDNN approximation with loss function $\bs L_{VgVP}(\theta_{\bs u},\theta_p, \theta_{\bs U},\theta_{\bs q})$.}
	\label{example3-6}%
\end{figure}

\begin{figure}[ht!]
	\center
	\subfigure[$Err(\bs u)$]{\includegraphics[scale=0.4]{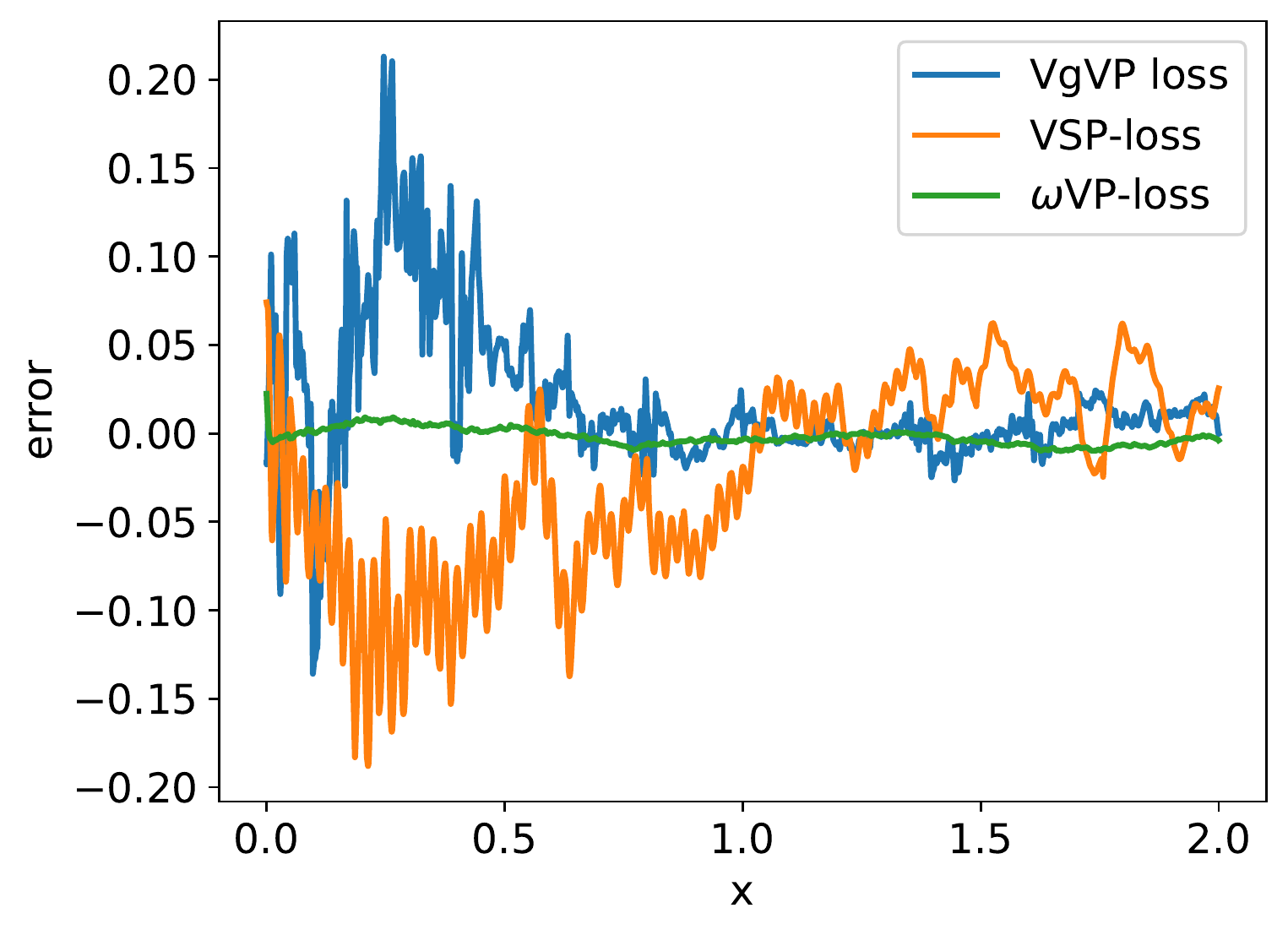}}\qquad
	\subfigure[$Err(p)$]{\includegraphics[scale=0.4]{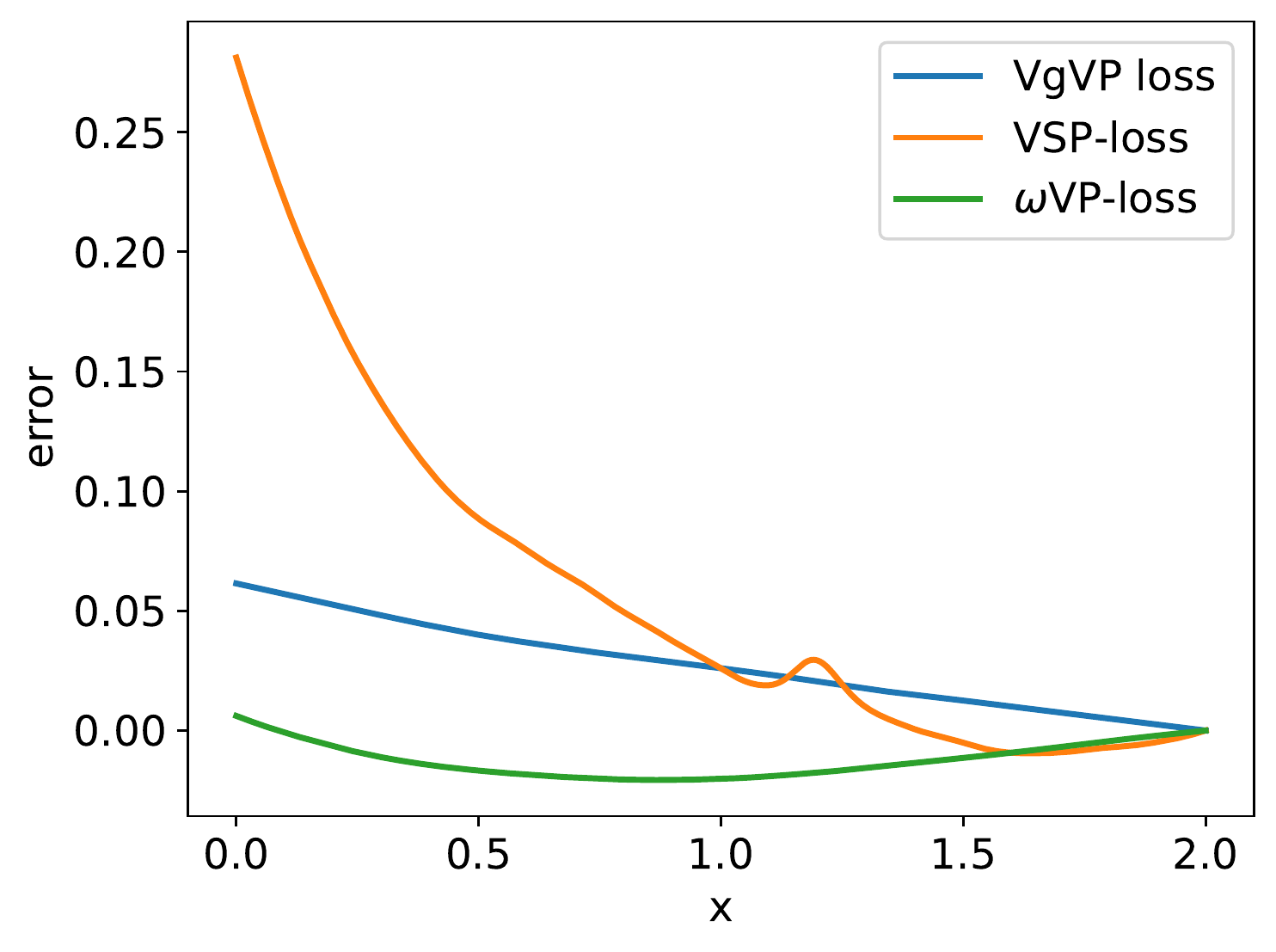}}
	\vspace{-10pt}
	\caption{Error of MscaleDNN approximations using different loss functions.}
	\label{example3-7}%
\end{figure}

The MscaleDNN solutions of $u_1$ are compared with the exact $u_1$ in Fig. \ref{example3-4}-\ref{example3-6}.  Here, we again plot the solutions along the line $y=0.7$ which does not cross any of cylinders inside the domain. Errors of the MscaleDNN approximations for $\bs u$ and $p$ using different losses are depicted in Fig. \ref{example3-7}.  We can see that the $\omega$VP-loss or $VgVP$-loss with the MscaleDNN can obtain very accurate solutions within 1500 epochs. Again, the VSP-loss need more learning to achieve similar accuracy.

%For comparison, we test algorithms using only fully connected DNNs. For $\bs u$ and variables $\bs\omega$, $\bs T$ and $\bs U$, we use fully connected DNNs with $8$ hidden layers and $1200$ neurons in each hidden layer. For $p$ and $\bs q$, we use fully connected DNNs with $8$ hidden layers and $300$ neurons in each hidden layer. Obviously, the total number of neurons in the fully connected DNNs and the MscaleDNNs are the same.  The losses and $\ell^2$-errors obtained by minimizing different loss functions in \eqref{firstorderloss} are compared in Fig.  \ref{example3-1}-\ref{example3-3}.


\begin{figure}[ht!]
	\center
	\subfigure[loss]{\includegraphics[scale=0.3]{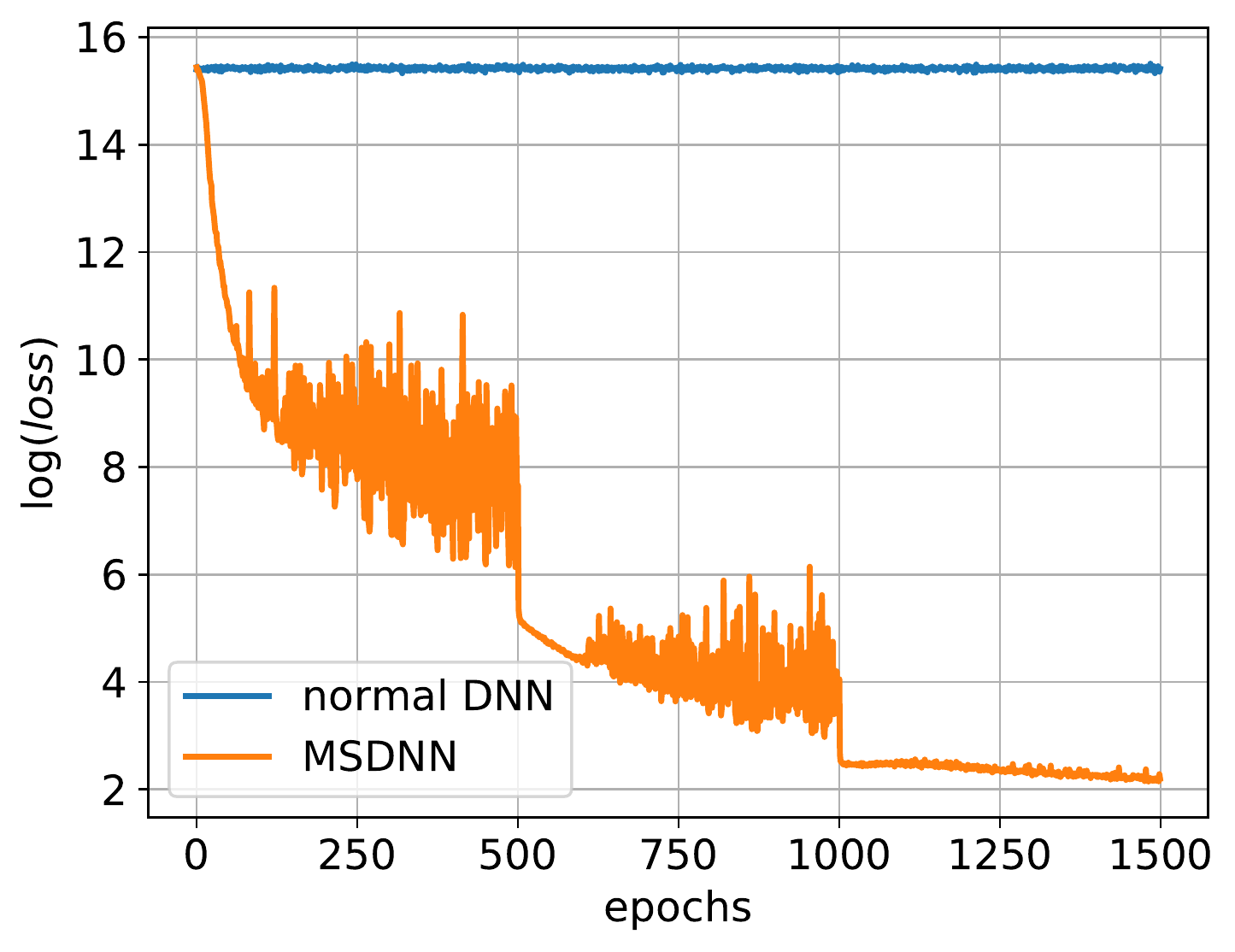}}
	\subfigure[$Err(\bs u)$]{\includegraphics[scale=0.3]{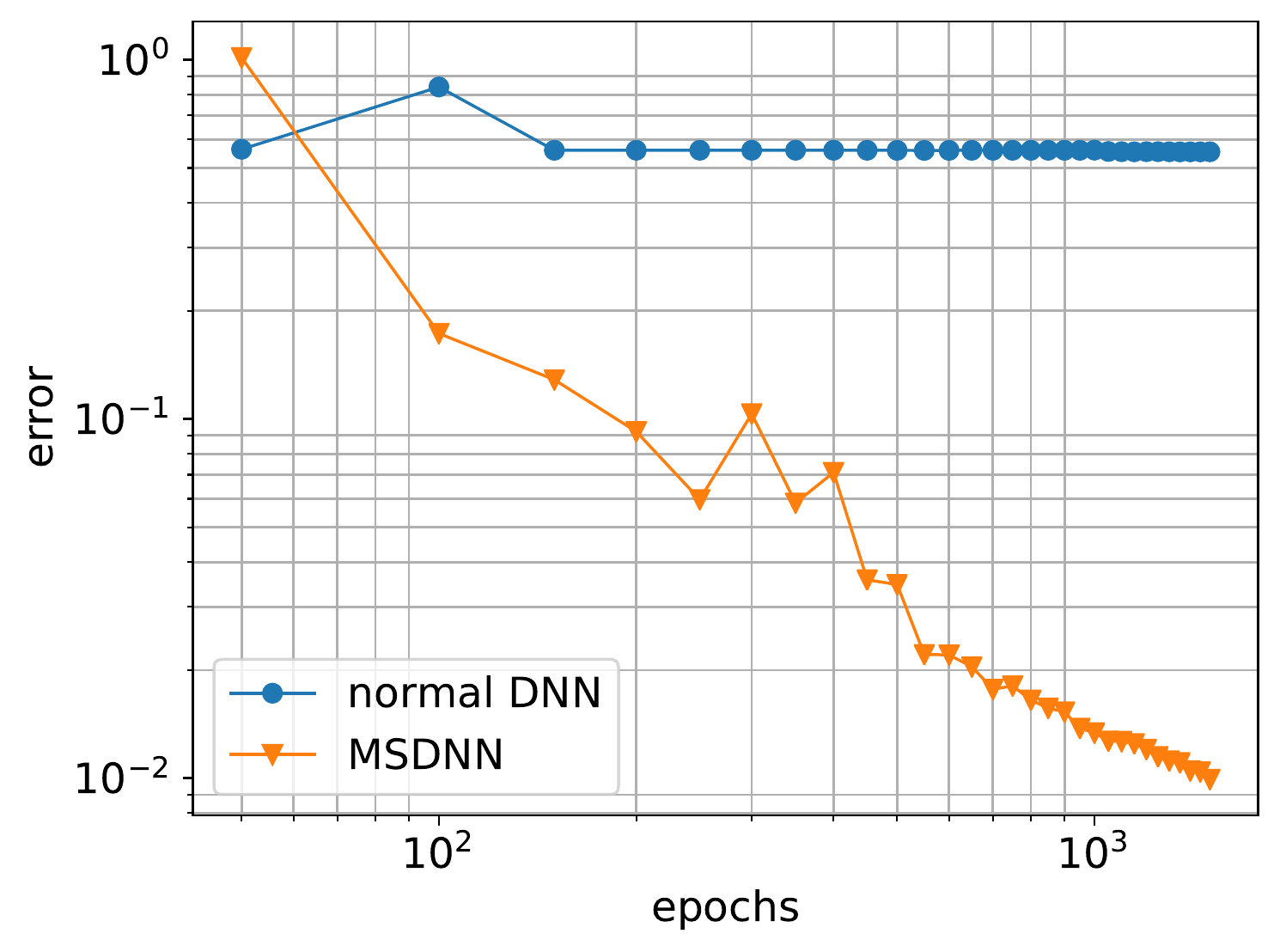}}
	\subfigure[$Err(p)$]{\includegraphics[scale=0.3]{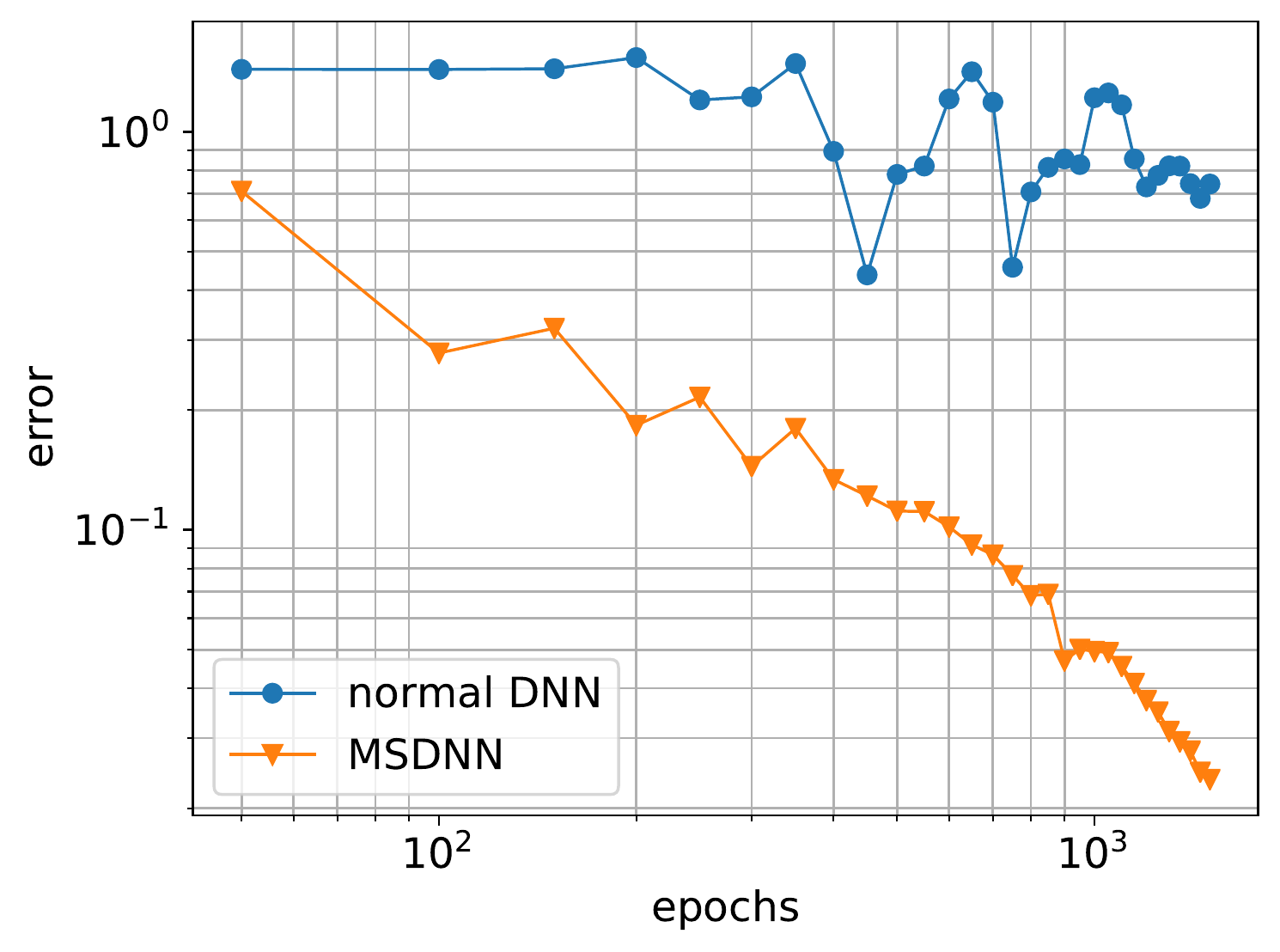}}
	\vspace{-10pt}
	\caption{Comparison of a normal DNN and the MscaleDNN with loss function  $\bs L_{\mbox{$\omega$}VP}(\theta_{\bs u},\theta_p, \theta_{\bs \omega},\theta_{\bs q})$.}
	\label{example3-1}%
\end{figure}

\begin{figure}[ht!]
	\center
	\subfigure[loss]{\includegraphics[scale=0.3]{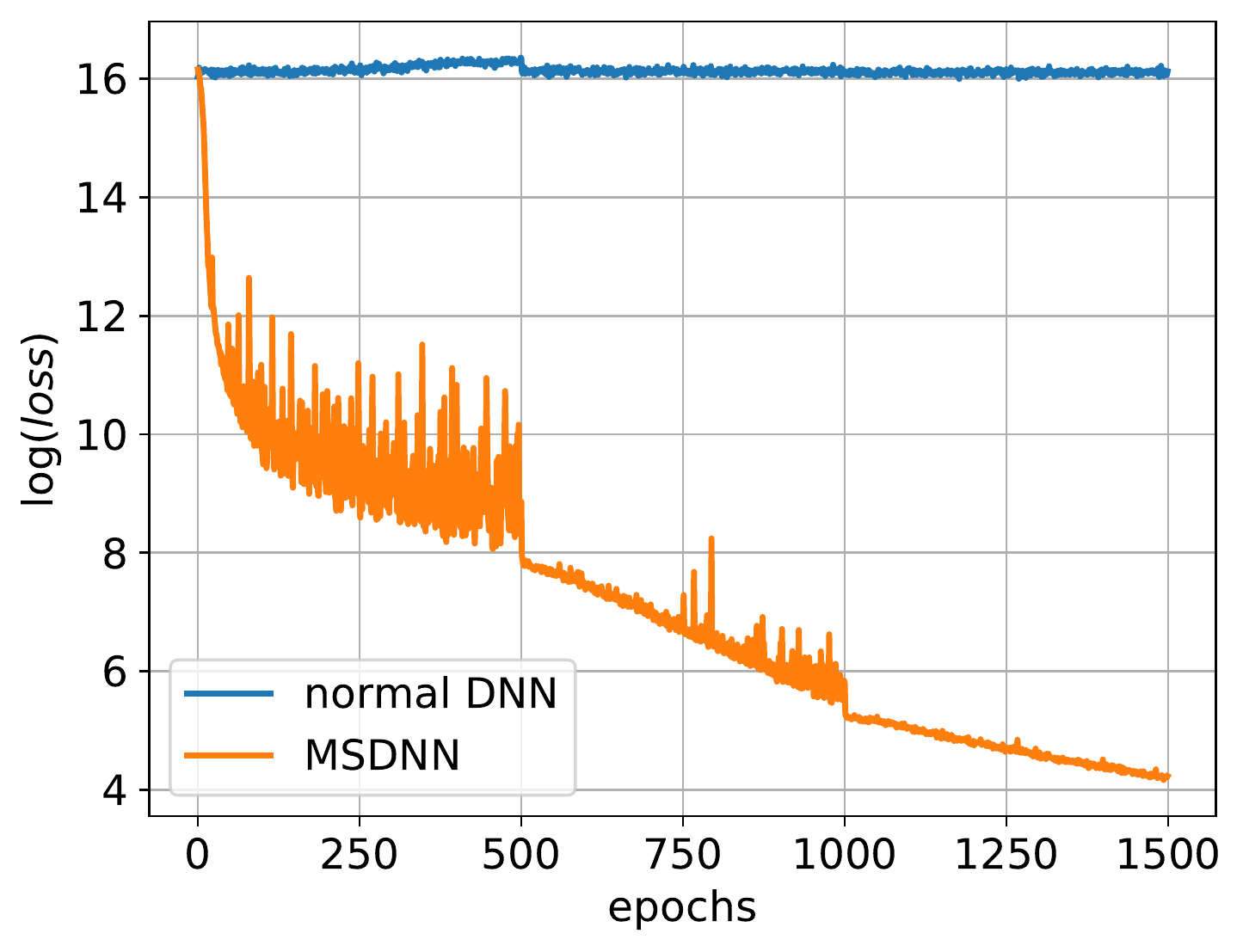}}
	\subfigure[$Err(\bs u)$]{\includegraphics[scale=0.3]{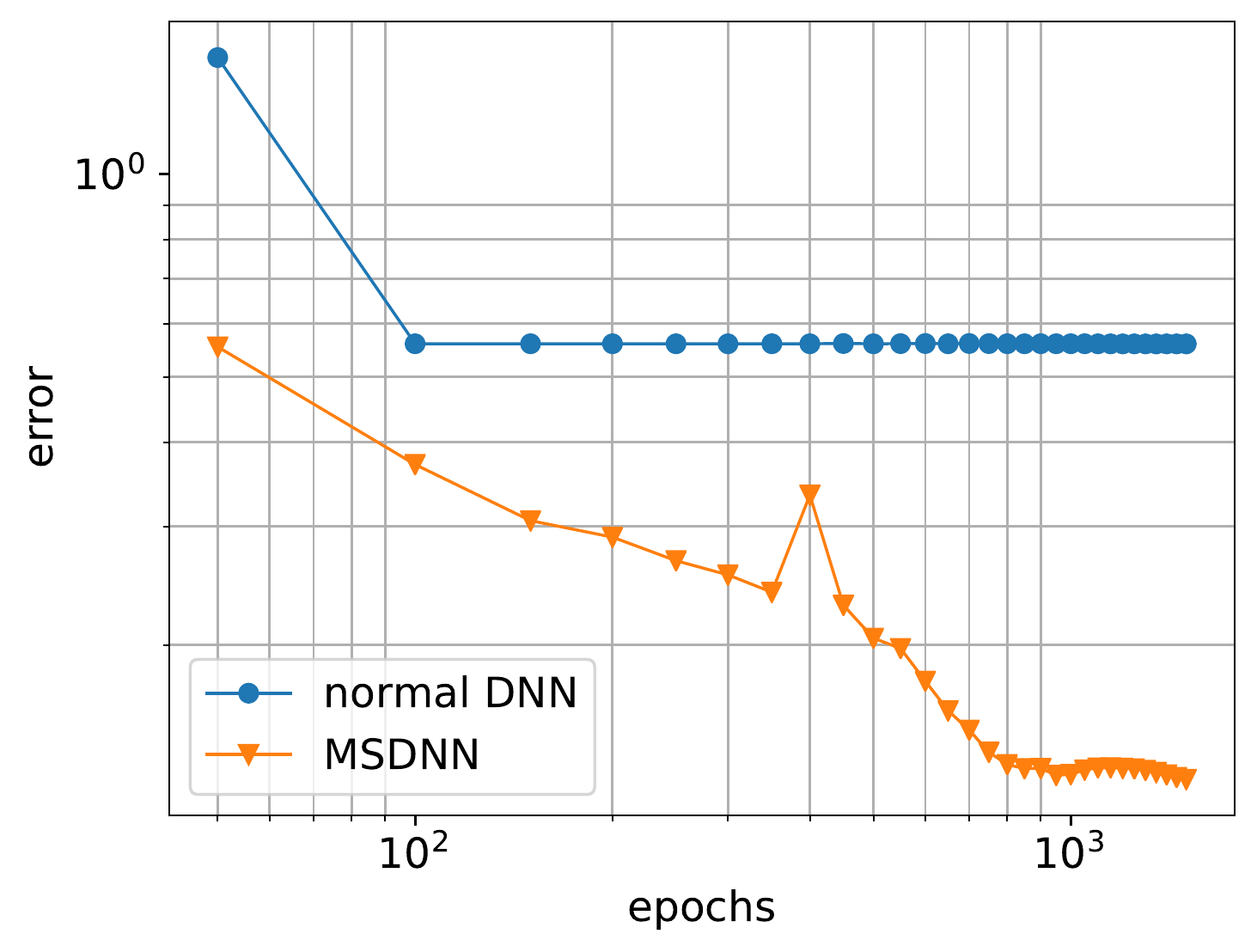}}
	\subfigure[$Err(p)$]{\includegraphics[scale=0.3]{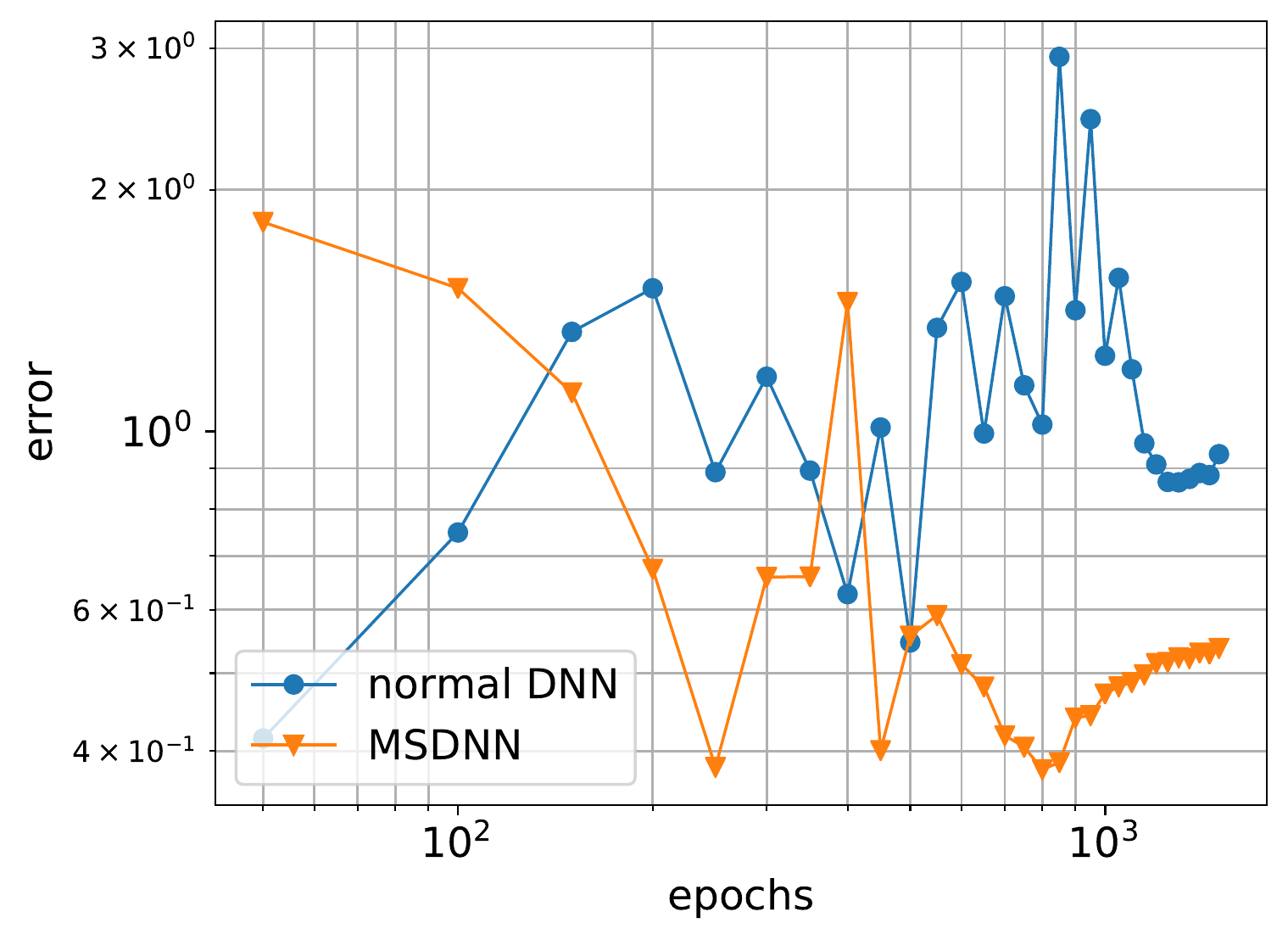}}
	\vspace{-10pt}
	\caption{Comparison of a normal DNN and the MscaleDNN with loss function  $\bs L_{VSP}(\theta_{\bs u},\theta_p, \theta_{\bs T},\theta_{\bs q})$.}
	\label{example3-2}%
\end{figure}

\begin{figure}[ht!]
	\center
	\subfigure[loss]{\includegraphics[scale=0.3]{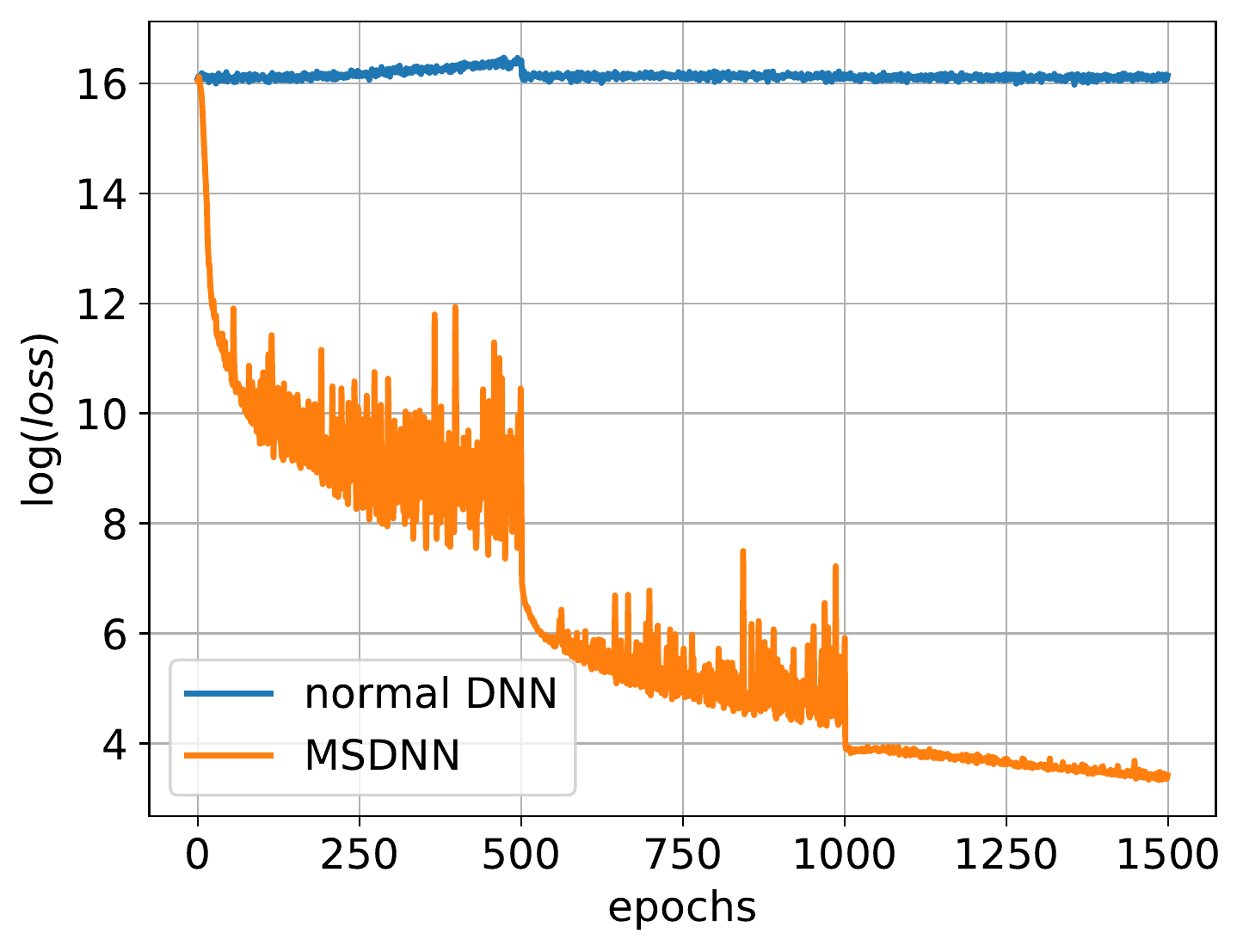}}
	\subfigure[$Err(\bs u)$]{\includegraphics[scale=0.3]{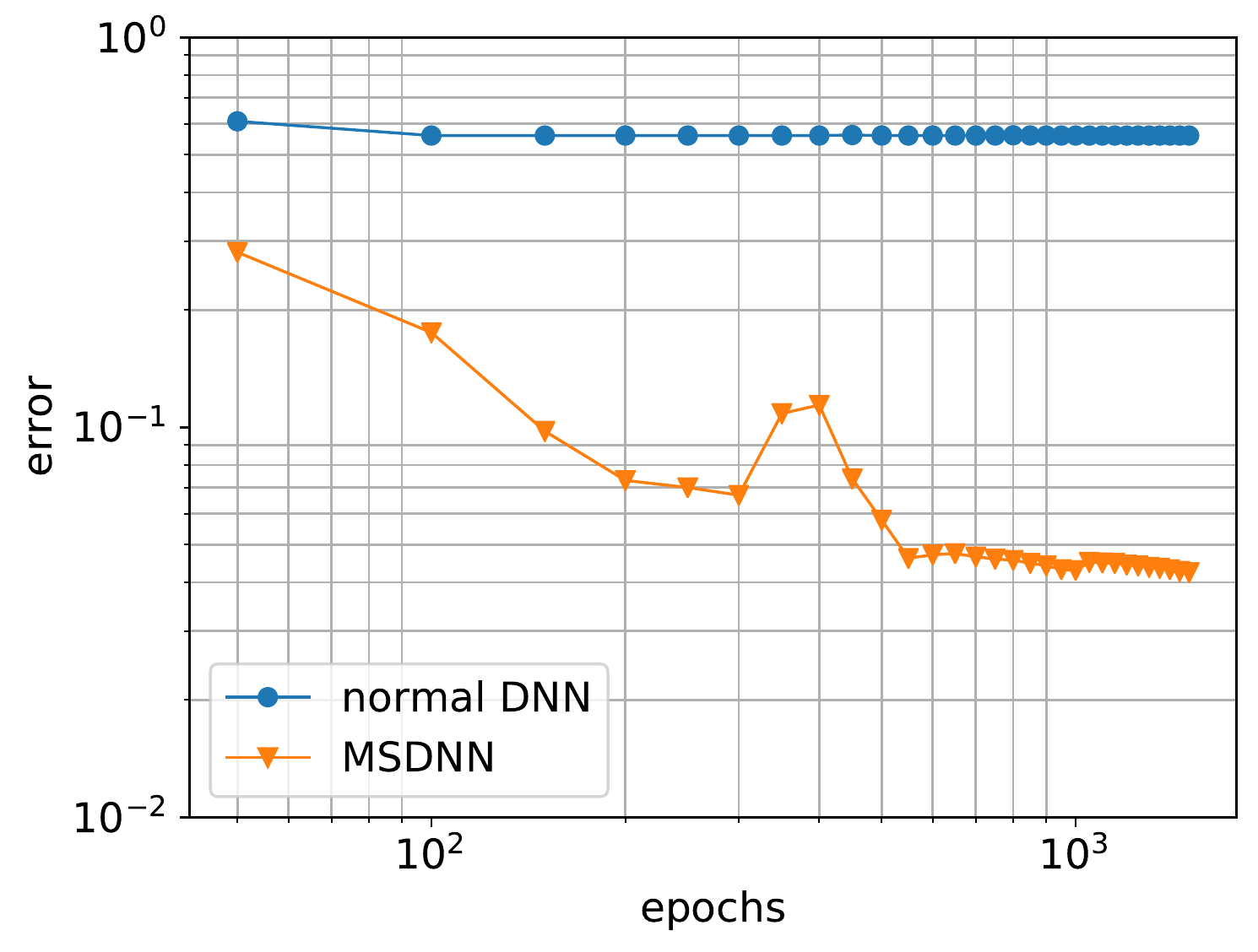}}
	\subfigure[$Err(p)$]{\includegraphics[scale=0.3]{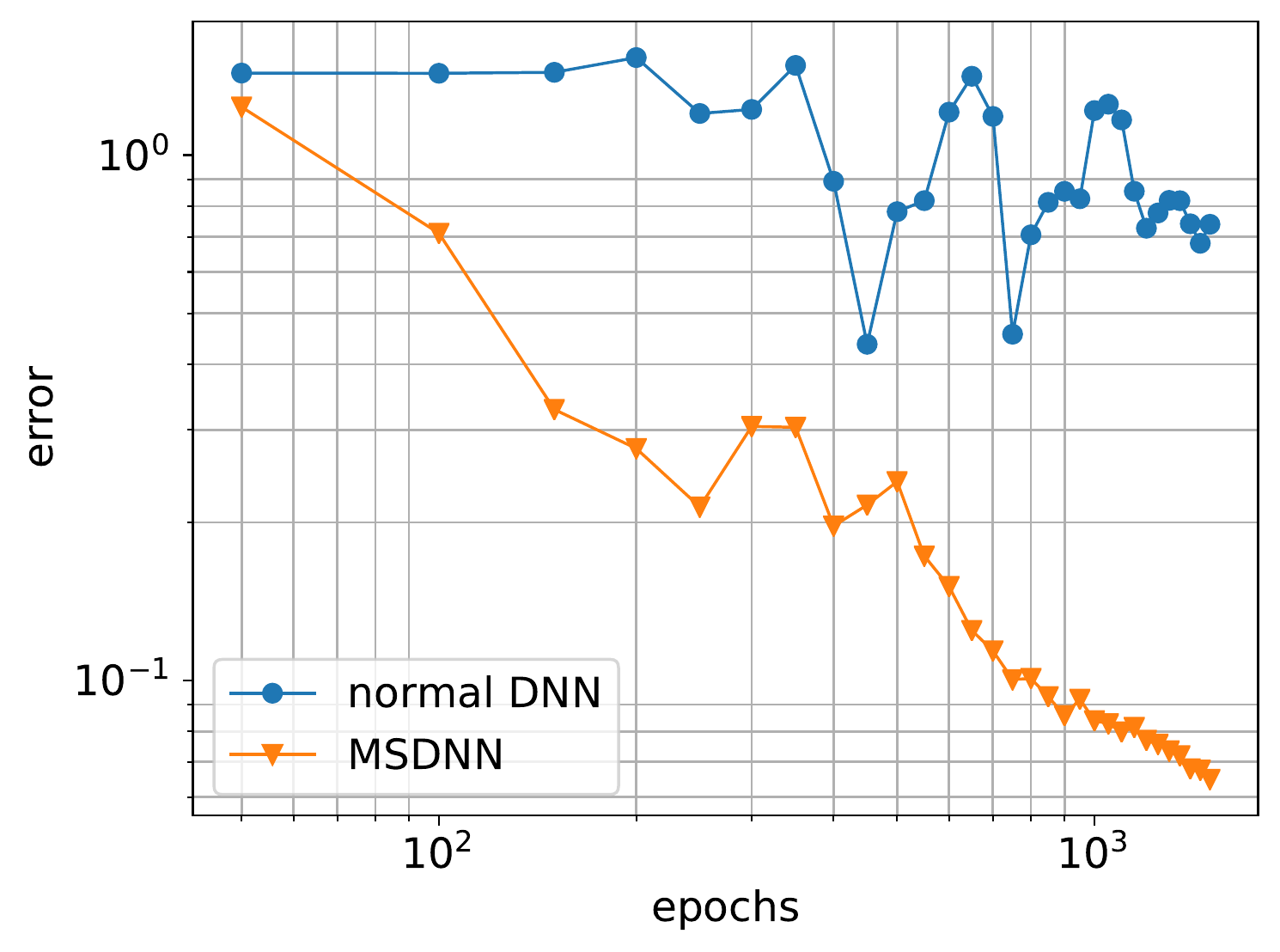}}
	\vspace{-10pt}
	\caption{Comparison of a normal DNN and the MscaleDNN with loss function  $\bs L_{VgVP}(\theta_{\bs u},\theta_p, \theta_{\bs U},\theta_{\bs q})$.}
	\label{example3-3}%
\end{figure}

For comparison, we test algorithms using only fully connected DNNs. For $\bs u$ and variables $\bs\omega$, $\bs T$ and $\bs U$, we use fully connected DNNs with $8$ hidden layers and $1200$ neurons in each hidden layer. For $p$ and $\bs q$, we use fully connected DNNs with $8$ hidden layers and $300$ neurons in each hidden layer. Again, the total number of neurons in the fully connected DNNs and the MscaleDNNs are the same.  The losses and $\ell^2$-errors obtained by minimizing different loss functions in \eqref{firstorderloss} are compared in Fig.  \ref{example3-1}-\ref{example3-3}, which clearly show the fast convergence of the MscaleDNNs when the normal fully connected DNNs fail to converge at all.

\subsection{Error in pressure $p$ and Poisson equation}
\label{section_P}

\begin{figure}[ht!]
	\center
	\subfigure[loss]{\includegraphics[scale=0.24]{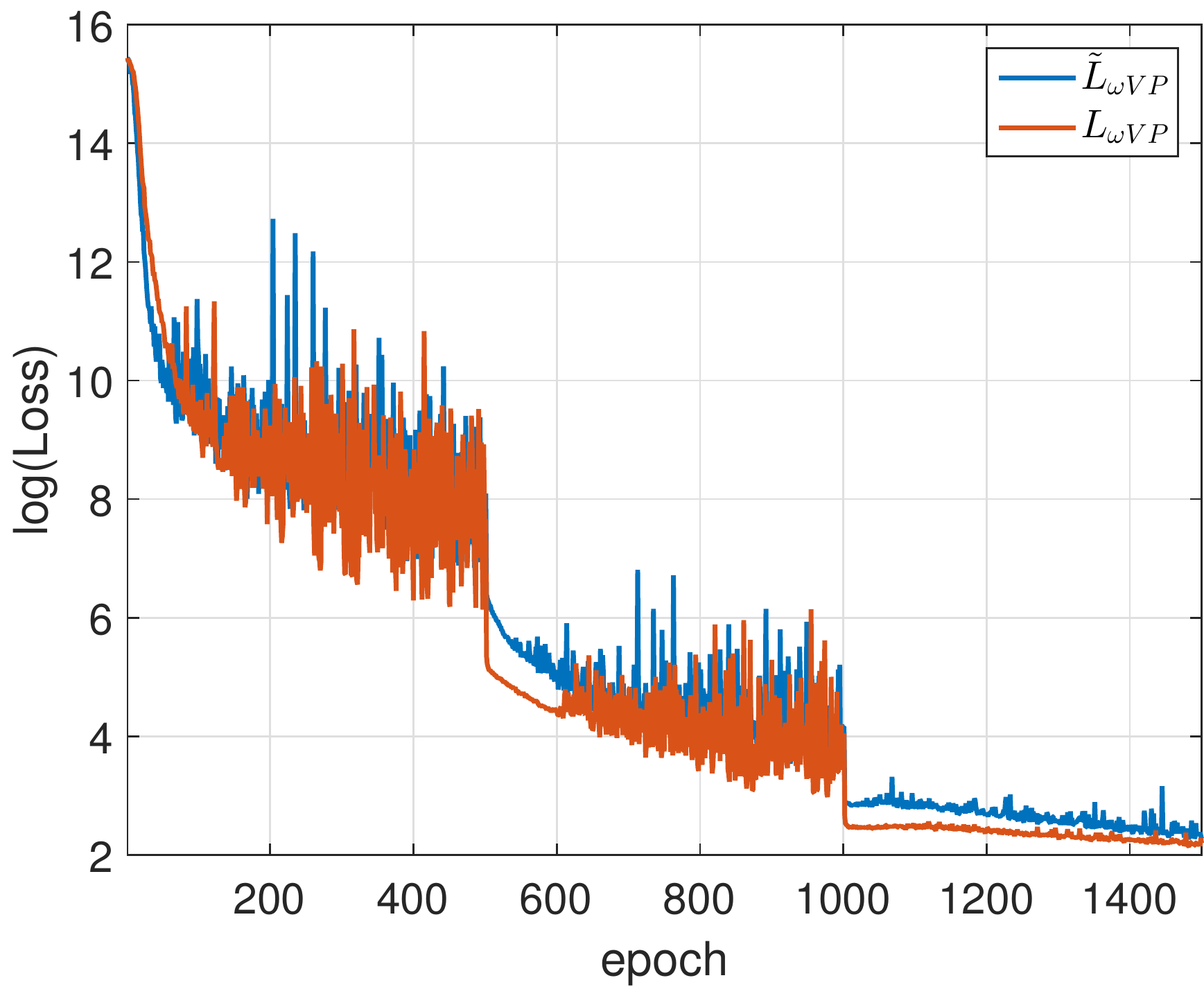}}
	\subfigure[$Err(\bs u)$]{\includegraphics[scale=0.24]{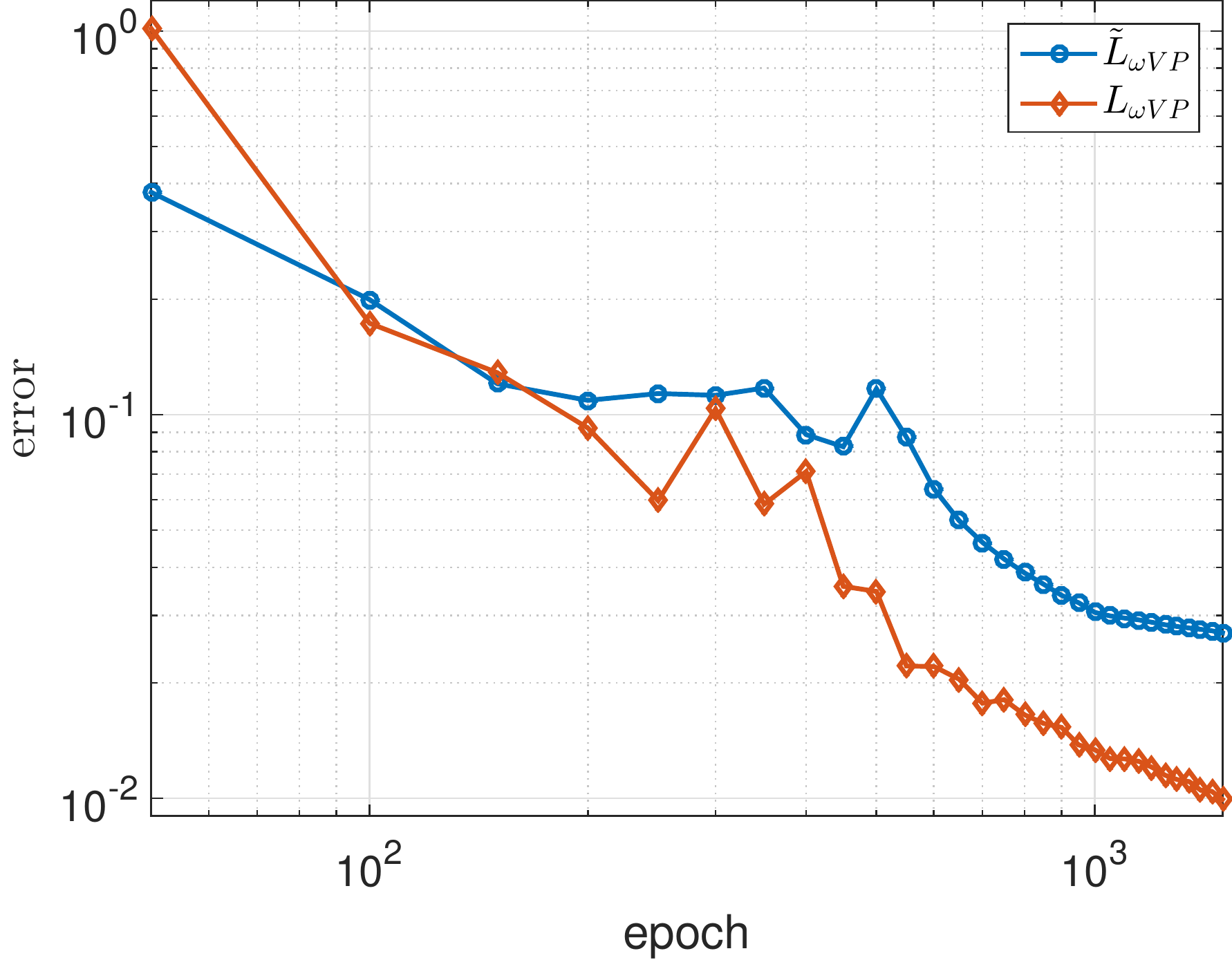}}
	\subfigure[$Err(p)$]{\includegraphics[scale=0.24]{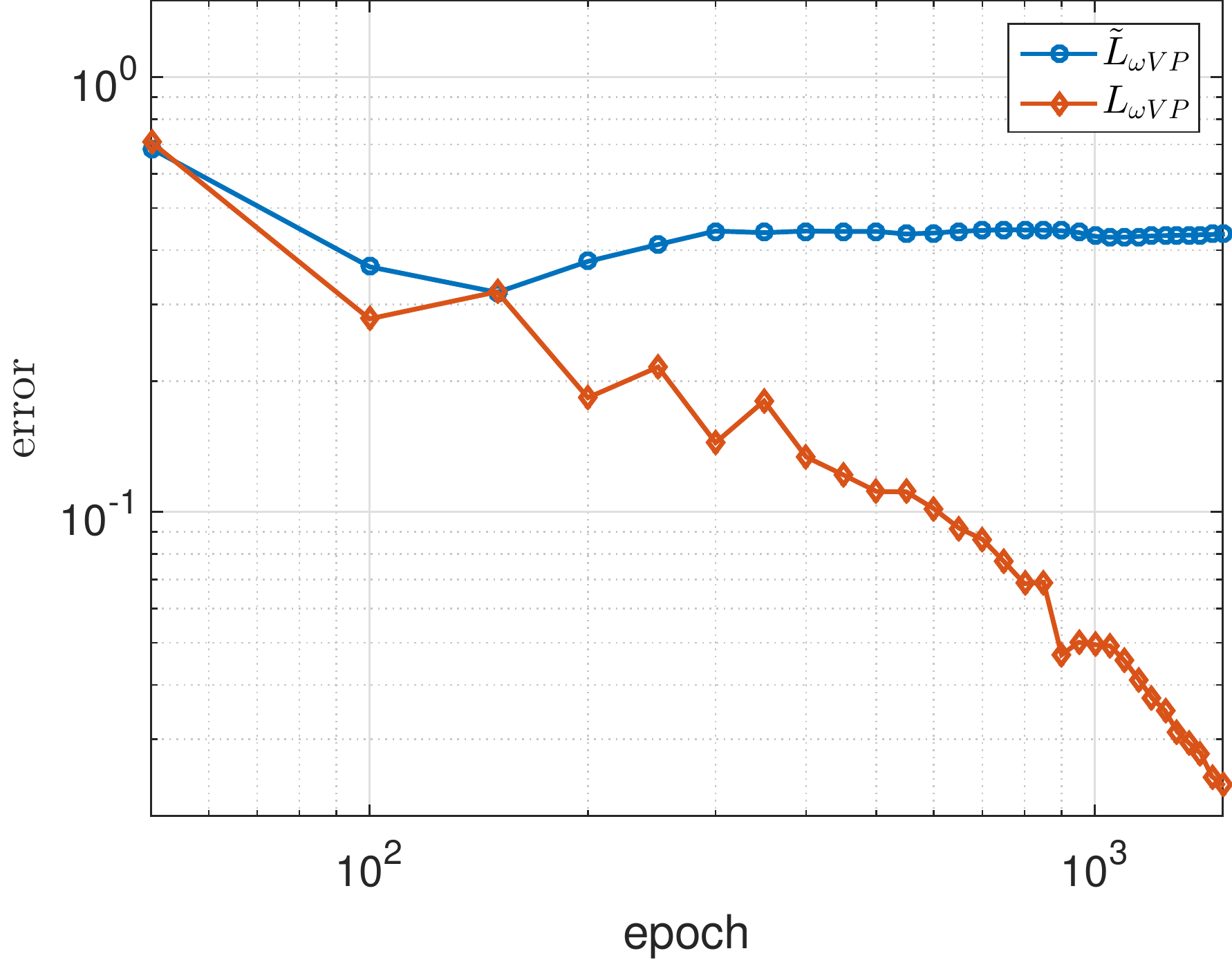}}
	\vspace{-10pt}
	\caption{Effect of Poisson equation in the loss function: Loss functions $\bs L_{\mbox{$\omega$}VP}(\theta_{\bs u},\theta_p, \theta_{\bs \omega},\theta_{\bs q})$ (red lines and diamonds) vs. $\tilde{\bs L}_{\mbox{$\omega$}VP}(\theta_{\bs u},\theta_p, \theta_{\bs \omega})$ (blue lines and circles).}
	\label{example3-8}%
\end{figure}

It is a well-known fact that the traditional projection methods for incompressible flow may experience an error degeneration for pressure near the boundaries depending on the types of pressure boundary conditions used for the Poisson equation \eqref{poissonp} \cite{gk91}. To show the importance of the pressure's Poisson equation in the DNN-based approaches for the Stokes problem, we will study a loss function without explicitly including the residual of the Poisson equation. Here, we consider a modification of the loss function ${\bs L}_{\mbox{$\omega$}VP}(\theta_{\bs u}, \theta_p, \theta_{\bs\omega},\theta_{\bs q})$ given by
\begin{equation}
	\widetilde{\bs L}_{\mbox{$\omega$}VP}(\theta_{\bs u}, \theta_p, \theta_{\bs\omega}):=\|\nu\nabla\times\bs \omega+\nabla p-\bs f\|^2_{\Omega}+\|\nabla\times \bs u-\bs \omega\|^2_{\Omega}+\|\nabla\cdot\bs u\|^2_{\Omega}+\beta\|\bs u-\bs g\|_{\partial\Omega}^2.
\end{equation}
The input data, size of the MscaleDNNs and other settings are exactly the same as we have used in the numerical tests in Section \ref{section_multiF}. The loss, errors of the MscaleDNN solutions are compared with the algorithm using loss function $\bs L_{\bs{\omega}VP}$ in Fig. \ref{example3-8}. We can see that the loss and $Err(\bs u)$ are compatible. However, $Err(p)$ is significantly improved if the loss function with a Poisson equation residual is used.

\section{Conclusion and future work}
In this paper, we have studied the MscaleDNN methods for solving highly oscillatory Stokes flow in complex domains and demonstrated the capability of the MscaleDNN as a meshless and high resolution numerical method for simulating flows in complex domains. Several least square formulations of the Stokes equations using different forms of first order systems are used to construct the loss functions for the MscaleDNN learning. The numerical results have clearly demonstrated the increased resolution power of the MscaleDNN to capture the fine structures in the flow fields when the normal fully connected network with the same overall sizes fail to converge at all. The MscaleDNN shows the potential of DNN machine learning as a practical alternative numerical method to traditional finite element methods. The DNN-based methods have an obvious advantage of no need for expensive mesh generations and matrix solvers as for traditional mesh-based numerical methods nor the delicate treatment of pressure boundary conditions and incompressibility constrains of the flow field.

There are many unresolved issues for solving Navier-Stokes equation, among them the most important one is to understand the convergence property of the MscaleDNN learning. A related issue is to find adaptive strategies to dynamically selecting the penalty constants for various terms in the loss functions, which are sensitive for the performance of DNN based machine learning PDE algorithms. It should also be mentioned that the structure of MscaleDNN is amendable to adaptive selections of scales by either adding or removing a scale dynamically during learning, future work will be done to explore this feature as well as to apply the MscaleDNN to 3-D time-dependent incompressible flows.

\section*{Acknowledgments}
W.C. is supported by the U.S. Army Research Office (grant W911NF-17-1-0368). B. W. acknowledges the financial support provided by NSFC (grant 11771137,12022104).

%%Vancouver style references.
\bibliographystyle{plain}

\end{document}